

\input epsf.tex

\def\2{{1\over 2}}

\def\d{\delta}
\def\a{\alpha}
\def\b{\beta}
\def\g{\gamma}

\def\e{\epsilon}
\def\l{\lambda}
\def\o{\omega}
\def\D{\Delta}

\def\fun#1#2#3{#1\colon #2\rightarrow #3}

\def\abs#1{\vert #1 \vert}
\def\frac#1#2{{{#1} \over {#2}}}

\def\sqr{\sqrt}
\def\st{\;\colon\;}
\def\tends{\rightarrow}
\def\weak{\rightharpoonup}

\def\dx{\hbox{{\rm d}$x$}}

\def\dr{ {\rm d} }
\def\dt{\hbox{{\rm d}$t$}}

\def\R{{\bf R}}
\def\N{{\bf N}}
\def\Z{{\bf Z}}

\def\T{{\bf T}}

\def\M{{\cal M}}

\def\thm#1{\vskip 1 pc\noindent{\bf Theorem #1.\quad}\sl}
\def\lem#1{\vskip 1 pc\noindent{\bf Lemma #1.\quad}\sl}
\def\prop#1{\vskip 1 pc\noindent{\bf Proposition #1.\quad}\sl}

\def\proof{\rm\vskip 1 pc\noindent{\bf Proof.\quad}}
\def\fin{\par\hfill $\backslash\backslash\backslash$\vskip 1 pc}
\def\txt#1{\quad\hbox{#1}\quad}
\def\m{\mu}
\def\L{{\cal L}}

\def\dcal{{\cal D}}
\def\G{\Gamma}

\def\o{\omega}
\def\r{\rho}

\def\cin#1{\2\abs{{#1}}^2}
\def\cinn#1{{n\over 2}\abs{{#1}}^2}
\def\cinh#1{{1\over{2h}}\abs{{#1}}^2}

\def\2{\frac{1}{2}}
\def\inn#1#2{{\langle #1 ,#2\rangle}}
\def\Mt{{\cal M}_1(\T^p)}

\def\part{{\partial_{x}}}
\def\div{{\rm div}}

\def\pprime{{{}^\prime{}^\prime}}

\def\Mcal{{\cal M}}



\baselineskip= 17.2pt plus 0.6pt
\font\titlefont=cmr17
\centerline{\titlefont A time-step approximation scheme}
\vskip 1 pc
\centerline{\titlefont for a viscous version of the Vlasov equation}
\vskip 4pc
\font\titlefont=cmr12
\centerline{         \titlefont {Ugo Bessi}\footnote*{{\rm 
Dipartimento di Matematica, Universit\`a\ Roma Tre, Largo S. 
Leonardo Murialdo, 00146 Roma, Italy.}}   }{}\footnote{}{
{{\tt email:} {\tt bessi@matrm3.mat.uniroma3.it}Work partially supported by the PRIN2009 grant "Critical Point Theory and Perturbative Methods for Nonlinear Differential Equations}} 
\vskip 0.5 pc
 
\par
\vskip 2pc
\centerline{\bf Abstract}
Gomes and Valdinoci have introduced a time-step approximation scheme for a viscous version of Aubry-Mather theory; this scheme is a variant of that of Jordan, Kinderlehrer and Otto.  Gangbo and Tudorascu have shown that the Vlasov equation can be seen as an extension of Aubry-Mather theory, in which the configuration space is the space of probability measures, i. e. the different distributions of infinitely many particles on a manifold. Putting the two things together, we show that Gomes and Valdinoci's theorem carries over to a viscous version of the Vlasov equation. In this way, we shall recover a theorem of J. Feng and T. Nguyen, but by  a different and more "elementary" proof.

\vskip 2 pc
\centerline{\bf  Introduction}
\vskip 1 pc

The Vlasov equation models a group of particles governed by an external potential $V$ and a mutual interaction $W$; we shall always suppose that the particles move on the $p$-dimensional torus $\T^p\colon=\frac{\R^p}{\Z^p}$, that $V$ and $W$ are sufficiently regular and that $V$ depends periodically on time. More precisely,

\noindent 1) $V\in C^4(\T\times\T^p)$ and

\noindent 2) $W\in C^4(\T^p)$; moreover $W$, seen as a periodic potential on $\R^p$, is even: $W(x)=W(-x)$. Up to adding a constant, we can suppose that $W(0)=0$. 

Let $\Mcal_1(\T^p\times\R^p)$ denote the space of Borel probability measures on $\T^p\times\R^p$; we say that a continuous curve $\fun{\eta}{\R}{\Mcal_1(\T^p\times\R^p)}$ solves the Vlasov equation if it satisfies, in the weak sense, the continuity equation
$$\partial_t\eta_t+
\div_{(x,v)}(\eta_t\cdot(v,\partial_x P^{\eta_t}(x)))=0   \eqno (CE)$$
where $(x,v)$ are the position and velocity coordinates on 
$\T^p\times\R^p$, 
$$P^{\eta_t}(t,x)=V(t,x)+W^{\eta_t}(x)$$
and
$$W^{\eta_t}(x)=\int_{\T^p\times\R^p}W(x-y)\dr\eta_t(y,v)  .  $$
An idea underlying several papers (see for instance [1], [9], [10], [12]) is to consider the Vlasov equation as a Hamiltonian system with infinitely many particles, i. e. as a Hamiltonian system on the space $\Mt$ of probability measures on $\T^p$; in particular, one can define, on $\Mt$, both the Hopf-Lax semigroup and the Hamilton-Jacobi equation.

In this paper, we follow [13] adding a viscosity term to the Hopf-Lax semigroup; we want to check two things. The first one (theorem 1 below) is that the minimal characteristics are solutions of a Fokker-Planck equation whose drift is determined by Hamilton-Jacobi, exactly as in the case without viscosity. The second check we want to do is about the time-discretization method of [14], which was developed for a final condition linear on measures, say
$$U_f(\mu)=\int_{\T^p}f\dr\mu  .  $$
We would like to see if it survives when the final condition $U$ is merely differentiable. Our definition of differentiability will be a little different from the usual one: indeed, we shall approximate minimal characteristics through "discrete characteristics"; since we shall see that the latter always have a density, we shall differentiate $U$ as a function on $L^1(\T^p)$, i. e. $U^\prime (\mu)$ will be a scalar function, an element of $L^\infty(\T^p)$.

We are going to consider a Lagrangian on 
$\R\times\T^p\times\R^p$ given by 
$$L^{\g_t}(t,q,\dot q)=
\2|\dot q|^2-P^{\g_t}(t,q)    $$
whose Legendre transform is 
$$H^{\g_t}(t,q,p)=\2|p|^2+P^{\g_t}(t,q) . $$

\thm{1} Let $\fun{U}{\Mt}{\R}$ be Lipschitz for the 1-Wasserstein distance and differentiable in the sense of section 4 below; let 
$\L^p$ denote the Lebesgue measure on $\T^p$. Then, the following three points hold.

\noindent 1) For every $\mu\in\Mt$ and every $m\in\N$, the $\inf$ below is a minimum.
$$(\Lambda^mU)(\mu)\colon =\inf\left\{
\int_{-m}^0\dt\int_{\T^p}L^{\2\r}(t,x,Y(t,x))\r(t,x)\dx+
U(\r(0)\L^p)    \right\}   .   \eqno (1) $$
In the formula above, the $\inf$ is taken over all the Lipschitz vector fields $Y$; the curve of measures $\r$ is a weak solution of the Fokker-Planck equation
$$\left\{
\eqalign{
\partial_t\r_t-\D\r_t+\div(\r_t\cdot Y)&=0,\quad t\in[-m,0]\cr
\r_{-m}&=\mu      .
}     \right. \eqno (FP)_{-m,Y,\mu}    $$

\noindent 2) The operator $\Lambda^m$ defined in point 1) has the semigroup property
$$\Lambda^{m+n}U=\Lambda^m\circ\Lambda^n U\quad
\forall m,n\in\N   .   $$

\noindent 3) There is a vector field $Y$ minimal in (1); it is given by  
$Y=c-\partial_x u$, where $u$ solves the Hamilton-Jacobi equation with time reversed
$$\left\{
\eqalign{
\partial_t u +\Delta u-H^\r(t,x,-\partial_xu)&=0,\quad
t\in(-m,0)\cr
u(0,x)&=f        
}     \right.   \eqno(HJ)_{0,\r,f}   $$
for a suitable $f\in L^\infty(\T^p)$.

\rm

\vskip 1pc

Note that [8] contains a stronger version of this theorem; in a sense, the aim of this paper is to show that it is possible to prove part of [8] using the technique of [14]. 

We briefly expand on this technique: roughly speaking, the difference with [15] is that the entropy term is embedded in the kinetic energy. Let us be more precise and describe the time-step, which is backwards in time. Given a continuous function $U$ on $\Mt$, we are going to define 
$$U(-\frac{1}{n},\mu)=\min
\left\{
\int_{\T^p\times\R^p}[
\frac{1}{n}L^{\2\mu}(\frac{-1}{n},x,nv)+\log\g(x,v)
]    \g(x,v)   \dr\mu(x)\dr v  +U(\mu\ast\g)
\right\}   -\log\left(\frac{n}{2\pi}\right)^\frac{p}{2}$$
where the minimum is over all the functions $\g$ on 
$\T^p\times\R^p$ such that $\g(x,\cdot)$ is a probability density on $\R^p$ for all $x$. One should look at $\g$ as at the probability distribution of the velocities: a particle starting at $x$ has velocity $nv$ with probability $\g(x,v)$. Since $U$ is non linear there is some work to do in order to show that the minimal $\g$ exists; we shall prove this in section 1 below. In section 2, we prove a bound on the $L^\infty$ norm of the minimal; in section 3, we shall iterate backward the formula above, getting the "discrete value function" $U(\frac{j}{n},\mu)$ for $j\le 0$; naturally, we shall also get a discrete characteristic 
$\mu_\frac{j}{n},\mu_\frac{j+1}{n},\dots,\mu_0$. We shall show that the discrete value functions is bounded as the time-step tends to zero. In section 4 we reduce to the linear case  expressing the minima of section 1 in terms of the differential of $U$ at the endpoint of the discrete characteristic. In section 5, we discuss the regularity of the linear problem.  Thanks to this regularity, in section 6 we can prove that the discrete characteristics converge to a solution of the Fokker-Planck equation and that the discrete value function converges to a solution of Hamilton-Jacobi; this will end the proof of theorem 1.

\vskip 1pc

\vskip 2pc

\centerline{\bf \S 1}

\centerline{\bf The time-step: existence of the minimal}

\vskip 1pc

We begin with a few standard definitions.  

\vskip 1pc

\noindent{\bf Definitions.}  \noindent $\bullet$) We denote by $\Mt$ the space of Borel probability measures on $\T^p$.

\noindent $\bullet$) Let $\tilde x,\tilde y\in\R^p$, and let $x$, $y$ be their projections on $\T^p$. We define
$$|x-y|_{\T^p}\colon=\min_{k\in\Z^p}|\tilde x-\tilde y-k|  .  $$

\noindent $\bullet$) For $\l\ge 1$ and $\mu_1,\mu_2\in\Mt$, we set
$$d_\l(\mu_1,\mu_2)^\l\colon=\min
\int_{\T^p\times\T^p}|x-y|_{\T^p}^\l\dr\G(x,y)$$
where the minimum is over all the measures $\G$ on 
$\T^p\times\T^p$ whose first and second marginals are $\mu_1$ and $\mu_2$ respectively; we recall from [2] that 
$(\Mt,d_\l)$ is a complete metric space whose topology is equivalent to the weak$\ast$ one.

\vskip 1pc

The term on the right in the formula above is a minimum by a standard theorem ([2], [16]); a useful characterization of $d_1$ is the dual one, i. e. 
$$d_1(\mu_1,\mu_2)=\sup\left\{
\int_{\T^p}f\dr\mu_1-\int_{\T^p}f\dr\mu_2
\right\}  \eqno (1.1)$$
where the $\sup$ is taken over all the functions $f\in C(\T^p)$ such that
$$|f(x)-f(y)|\le|x-y|_{\T^p}\qquad \forall x,y\in\T^p   .  $$

We need to adapt a few definitions of [14] to our situation.

\vskip 1pc

\noindent{\bf Definitions.} $\bullet$) Let $\mu\in\Mt$. We define 
$\dcal_\mu$ as the set of all the Borel functions 
$\fun{\g}{\T^p\times\R^p}{[0,+\infty)}$ such that
$$\int_{\R^p}\g(x,v)\dr v=1\txt{for $\mu$ a. e. } x\in\T^p  .  
\eqno (1.2)$$

\noindent $\bullet$) We denote by 
$$\fun{\pi_{\T^p}}{\T^p\times\R^p}{\T^p},\qquad
\fun{\pi_{\R^p}}{\T^p\times\R^p}{\R^p},\qquad
\fun{\pi_{cover}}{\R^p}{\T^p}$$
the natural projections, and define 
$\fun{\tilde\pi}{\T^p\times\R^p}{\T^p}$ by 
$\tilde\pi=\pi_{cover}\circ\pi_{\R^p}$.

\noindent $\bullet$) If $\mu\in\Mt$ and $\g\in\dcal_\mu$, we define a measure on $\T^p$ by 
$$\mu\ast\g=
(\pi_{\T^p}-\tilde\pi)_\sharp(\mu\otimes(\g(x,\cdot)\L^p))$$
where $\L^p$ denotes the Lebesgue measure on $\R^p$; the sharp sign denotes, as usual, the push-forward of a measure. In other words, if $f\in C(\T^p)$, then
$$\int_{\T^p}f(z)\dr(\mu\ast\g)(z)=
\int_{\T^p\times\R^p}f(x-v)\g(x,v)\dr\mu(x)\dr v  .  $$
Note that, if $\g$ does not depend on $x\in\T^p$, this is the usual convolution of the two measures $\mu$ and $\g\L^p$. One can see $\g$ as the probability, for a particle placed in $x$, to jump to  $x-v$; if the initial distribution of the particles is $\mu$, 
$\mu\ast\g$ is the distribution after one jump.

\noindent $\bullet$) Let now $U\in C(\Mt,\R)$; for $h>0$ and 
$t\in\R$ we define
$$\fun{G^h_tU}{\Mt}{\R}$$
by
$$(G^h_tU)(\mu)=\inf_{\g\in\dcal_\mu}\left\{
\int_{\T^p\times\R^p}[
hL^{\2\mu}(t,x,\frac{1}{h}v)+\log\g(x,v)
]  \g(x,v)\dr\mu(x)\dr v+U(\mu\ast\g)
\right\}     $$
where the Lagrangian $L_c^{\2\mu}$ has been defined in the introduction.

\vskip 1pc

\noindent{\bf Observation.} For $c\in\R^p$, it is natural to consider the Lagrangian
$$L_c^{\g_t}(t,q,\dot q)=\cin{\dot q}-\inn{c}{\dot q}-
P^{\g_t}(t,q)  .  $$
Naturally, it is possible to prove theorem 1 for $L_c^{\g_t}$. Indeed, let 
$$\fun{\tau_c}{\T^p}{\T^p},\qquad
\fun{\tau_c}{x}{x+hc}   $$
and 
$$\hat U(\mu)=U((\tau_{hc})_\sharp\mu)  .  $$
If we set $\tilde\g(x,v)=\g(x,v+hc)$, it is easy to see that
$$\int_{\T^p\times\R^p}[
hL^{\2\mu}_c(t,x,\frac{1}{h}v)+\log\g(x,v)
]  \g(x,v)\dr v\dr\mu(x)+U(\mu\ast\g)=$$
$$\int_{\T^p\times\R^p}[
hL^{\2\mu}_0(t,x,\frac{1}{h}v)+\log\tilde\g(x,v)
]  \tilde\g(x,v)\dr v\dr\mu(x)+\hat U(\mu\ast\tilde\g) -
\frac{h}{2}|c|^2 .  $$
In other words, a simple transformation brings the minima for 
$L^{\2\mu}_c$ into those for $L^{\2\mu}_0$. We have restricted statement and proof of theorem 1 to the case $c=0$ to keep the notation (relatively) simple.

We want to write $G^h_t U$ in a different way. First of all, we define
$$A_h(\g,(x,v))=
\frac{1}{2h}|v|^2\g(x,v)+\g(x,v)\log\g(x,v)  .  $$
If $\g$ does not depend on $x\in\T^p$, we shall call this function 
$A_h(\g,v)$.

Then, we note that the minimal $\g$ does not depend on the potential in $L_c^{\2\mu}$ (though the value function $G^h_tU$ obviously does); indeed, since 
$\g\in{\cal D}_\mu$, if $Z$ is any potential on $\T^p$, we have by Fubini 
$$\int_{\T^p\times\R^p}
Z(x)\g(x,v)\dr\mu(x)\dr v=
\int_{\T^p}Z(x)\dr\mu(x)   .  \eqno (1.3)$$
As a consequence,
$$(G^h_tU)(\mu)=\int_{\T^p}P^{\2\mu}(t,x)\dr\mu(x)+
\inf_{\g\in\dcal_\mu} S(U,\mu,\g)  \eqno (1.4)$$
where the single particle functional $S$ is given by
$$S(U,\mu,\g)=
\int_{\T^p\times\R^p}A_h(\g,(x,v))\dr\mu(x)\dr v+U(\mu\ast\g)  
\eqno (1.5)   $$
and the potential $P^{\2\mu}$ is as in the introduction.

\vskip 1pc

\noindent{\bf Observation.} We must show that the integral in (1.5) is well-defined, though possibly $+\infty$. Indeed, denoting by 
$f^-$ the negative part of a function $f$, we have that
$$\int_{\T^p}\dr\mu(x)\int_{\R^p}A_h^-(\g,(x,v))\dr v=$$
$$\int_{\T^p}\dr\mu(x)\int_{\R^p}\left[
\frac{1}{2h}|v|^2\g(x,v)+\g(x,v)\log\g(x,v)
\right]^-\dr v\ge
-\int_{\T^p}\dr\mu(x)\int_{\R^p}e^{
-1-\frac{1}{2h}|v|^2
}  \dr v=-e^{-1}(2\pi h)^\frac{p}{2}     \eqno (1.6)$$
where the inequality comes from the fact that 
$$\frac{1}{2h}|v|^2x+x\log x\ge -e^{
-1-\frac{1}{2h}|v|^2
}   \qquad\forall x\ge 0  .  $$

\vskip 1pc

We want to prove that the $\inf$ in (1.4) is a minimum; since $U$ is nonlinear, we cannot write the minimum explicitly as in [14]; we shall need a few lemmas, the first of which is an elementary fact on the behaviour of the Gaussian.

\lem{1.1} Let $h,\e>0$. Then, 
$$\min\left\{
\int_{\R^p}A_h(\g,v)\dr v\st \g\ge 0,\quad \int_{\R^p}\g(v)\dr v=1,
\quad
\int_{\R^p}\frac{1}{2h}|v|^2\g(v)\dr v=\frac{p\e}{2}   \right\}=$$
$$\frac{p\e}{2}+\log\frac{1}{(2\pi\e h)^\frac{p}{2}}-
\frac{p}{2} .  \eqno (1.7)  $$

\proof We note that we are minimizing the strictly convex functional
$$\fun{J}{L^1((1+\2|v|^2)\L^p)}{\R\cup{+\infty}},\qquad
\fun{J}{\g}{\int_{\R^p}A_h(\g,v)\dr v}$$
on the closed convex set
$$H=\left\{
\g\in L^1((1+\2|v|^2)\L^p)\st \g\ge 0,\quad\int_{\R^p}\g(v)\dr v=1,\quad
\int_{\R^p}\frac{1}{2h}|v|^2\g(v)\dr v=\frac{p\e}{2}
\right\}    .   $$
It is standard (see for instance the argument of proposition 1 of [14] or proposition 5.6 of chapter 1 of [6]) that, if we find a  density $\g$ and $\chi,\d\in\R$ solving the Lagrange multiplier problem
$$\left\{
\eqalign{
\frac{1}{2h}|v|^2+\log\g(v)+1&=\chi+\frac{\d}{2h}|v|^2\cr
\int_{\R^p}\g(v)\dr v&=1\cr
\int_{\R^p}\frac{1}{2h}|v|^2\g(v)\dr v&=\frac{p\e}{2}    
}
\right.    \eqno (1.8)$$
then $\g$ is the unique minimizer of $J$ restricted to $H$. Thus, solving (1.8) is next in the order of business.

By the first one of (1.8), we see that
$$\g(v)=e^{\chi-1}e^{-\frac{1-\d}{2h}|v|^2}  .  $$
Since we want $\g\in L^1$, eventually we shall have to check that $\d<1$. The constant $\chi$ is the unique one for which the second formula of (1.8) holds, i. e. 
$$e^{\chi-1}=\left(
\frac{1-\d}{2\pi h}
\right)^\frac{p}{2}   .   $$
The constant $\d$ is chosen so that the third one of (1.8) holds:
$$\frac{p\e}{2}=\int_{\R^p}\frac{1}{2h}|v|^2
\left(
\frac{1-\d}{2\pi h}
\right)^\frac{p}{2}
e^{-\frac{1-\d}{2h}|v|^2}\dr v=
\frac{1}{(2\pi)^\frac{p}{2}}\cdot\frac{1}{1-\d}
\int_{\R^p}\2|y|^2e^{-\frac{|y|^2}{2}}\dr y=
\frac{p}{2}\cdot\frac{1}{1-\d}   $$
where we have set $y=\sqrt\frac{1-\d}{h}v$.

From this we get 
$$1-\d=\frac{1}{\e}  .  $$
Since $\e>0$, this implies that $\d<1$, as we wanted. From the last four formulas,
$$\g(v)=\left(
\frac{1}{2\pi\e h}
\right)^\frac{p}{2}  e^{
-\frac{1}{2\e h}|v|^2
}   .  $$
This yields the first equality below, while the second one follows from (1.8).
$$\int_{\R^p}\g(v)\log\g(v)\dr v=
\int_{\R^p}\g(v)[
\log\frac{1}{(2\pi\e h)^\frac{p}{2}}-\frac{1}{2\e h}|v|^2
]   \dr v=
\log\frac{1}{(2\pi\e h)^\frac{p}{2}}-\frac{p}{2}   .   $$
From the formula above and the third one of (1.8), we get the second equality below.
$$\int_{\R^p}A_h(\g,v)\dr v=
\int_{\R^p}\left[
\frac{1}{2h}|v|^2\g(v)+\g(v)\log\g(v)
\right]   \dr v=
\frac{p\e}{2}+\log\frac{1}{(2\pi\e h)^\frac{p}{2}}-\frac{p}{2}  $$
which is (1.7).

\fin

\lem{1.2} Let $\mu\in\Mt$, let $C\in\R$ and let us consider the set 
$E_\mu$ of the functions $\g\in\dcal_\mu$ such that
$$\int_{\T^p\times\R^p}A_h(\g,(x,v))\dr\mu(x)\dr v\le C  .  
\eqno (1.9)$$
Then, 

\noindent 1) $E_\mu$ is uniformly integrable for the measure 
$\mu\otimes\L^p$ on $\T^p\times\R^p$. 

\noindent 2) The set of the measures 
$\{  \mu\otimes\g\L^p  \}$ as $\mu$ varies in $\Mt$ and $\g$ varies in $E_\mu$ is tight on $\T^p\times\R^p$.

\noindent 3) The set $\dcal_\mu$ is weakly closed in 
$L^1(\mu\otimes\L^p)$.

\proof We begin with point 1). We fix $a>1$ and consider 
$\g\in E_\mu$; the first inequality below is (1.9), the second one follows from Fubini, (1.6) and the fact that $\log\g\ge 0$ if 
$\g\ge a$; the last one is obvious.
$$C\ge\int_{\T^p\times\R^p}
\left[
\frac{1}{2h}|v|^2\g(x,v)+\g(x,v)\log\g(x,v)
\right] \dr\mu(x)  \dr v\ge$$
$$-e^{-1}(2\pi h)^\frac{p}{2}+
\int_{\T^p}\dr\mu(x)\int_{ \{ \g\ge a \} }\left[
\frac{1}{2h}|v|^2\g(x,v)+\g(x,v)\log\g(x,v)
\right]   \dr v\ge$$
$$-e^{-1}(2\pi h)^\frac{p}{2}+
\log a\int_{ \{ \g\ge a \} }\g(x,v)\dr\mu(x)\dr v   .  $$
This implies immediately that $E_\mu$ is uniformly integrable.

We prove point 2), i. e. that for all $\e>0$ we can find 
$R>0$ such that
$$\int_{\T^p\times B(0,R)^c}\g(x,v)\dr\mu(x)\dr v\le\e
\qquad\forall\mu\in\Mt,\quad\forall \g\in E_\mu  .  \eqno (1.10)$$
If we show that
$$\int_{\T^p\times\R^p}
\cinh{v}\g(x,v)\dr\mu(x)\dr v\le C_5\qquad
\forall\mu\in\Mt,\quad\forall\g\in E_\mu    $$
then (1.10) follows by the Chebishev inequality. By Fubini, the last formula is equivalent to
$$\int_{\T^p}r_\g(x)\dr\mu(x)\le C_6\qquad
\forall\mu\in\Mt,\quad\forall\g\in E_\mu
\eqno (1.11)   $$
where $r_\g$ is defined by 
$$\frac{p}{2}\cdot r_\g(x)\colon=\int_{\R^p}
\cinh{v}\g(x,v)\dr v    .    $$
Since $\mu$ is a probability measure, (1.11) follows if we prove that, for some $A>0$, there is $C_7>0$ such that, 
$$\int_{
\{ x\st r(x)>A \}
}   r_\g(x)\dr\mu(x)   \le C_7\qquad
\forall\mu\in\Mt,\quad\forall\g\in E_\mu    .   \eqno (1.12)$$
We call $g(\e)$ the function on the right hand side of (1.7); the first inequality below comes from (1.9) and Fubini, the second one from (1.7).
$$C\ge\int_{\T^p}\dr\mu(x)
\int_{\R^p}A_h(\g,(x,v))\dr v\ge
\int_{\T^p}g(r_\g(x))\dr\mu(x)   .   $$
Since the logarithmic term in the definition of $g$ grows less than linearly, we easily get that there is $A>0$ such that, for $y\ge A$, we have 
$g(y)\ge\frac{y}{4}$; since $g$ is bounded from below, the last formula implies that there is $C_8>0$, independent on $\g$ and 
$\mu$, such that
$$C_8\ge\int_{
\{ x\st r_\g(x)\ge A \}
}     
\frac{r(x)}{4}\dr\mu(x)\qquad
\forall\mu\in\Mt,\quad\forall\g\in E_\mu    .   $$
But this is (1.12).

We prove point 3). Let $B\subset\T^p$ be a Borel set; the function 
$$\fun{}{\g}{
\int_{B}\dr\mu(x)\int_{\R^p}\g(x,v)\dr v
}    $$
is continuous for the weak topology of $L^1(\mu\otimes\L^p)$; moreover, if $\g\in\dcal_\mu$,
$$\int_{B}\dr\mu(x)\int_{\R^p}\g(x,v)\dr v=
\mu(B)   .   $$
As a result, if $\bar\g$ belongs to the weak closure of 
$\dcal_\mu$, then
$$\int_{
B
}   \dr\mu(x)\int_{\R^p}\bar\g(x,v)\dr v=
\mu(B)       $$
for every Borel set $B\subset\T^p$. If we set
$$R(x)=\int_{\R^p}\bar\g(x,v)\dr v    $$
the last formula implies that
$$\mu(B)=\int_{B}R(x)\dr\mu(x)  $$
for every Borel set $B\subset\T^p$. 
It is standard that this implies that $R(x)=1$ for 
$\mu$ a. e. $x\in\T^p$, i. e. that $\bar\g\in\dcal_\mu$.

\fin

\lem{1.3} Let $U\in C(\Mt)$ and let $\mu\in\Mt$; then the function
$$\fun{I}{\dcal_\mu}{\R}$$
$$\fun{I}{\g}{\int_{\T^p\times\R^p}}\left[
hL^{\2\mu}(t,x,\frac{1}{h}v)+\log\g(x,v)
\right]   \g(x,v)\dr\mu(x)\dr v+   U(\mu\ast\g)$$
is l. s. c. for the weak topology of $L^1(\mu\otimes\L^p)$. 

\proof {\bf Step 1.} We begin to show that the function 
$$\fun{}{\g}{U(\mu\ast\g)}$$
is continuous; since we are supposing that $\fun{U}{\Mt}{\R}$ is continuous, it suffices to prove that 
$\fun{}{\g}{\mu\ast\g}$ is continuous from ${\cal D}_\mu$ endowed with the weak topology of $L^1(\mu\otimes\L^p)$ to the 
weak$\ast$ topology of 
$\Mt$. Let $\g\in\dcal_\mu$ be fixed and let $f\in C(\T^p)$; it suffices to note that we can write the weak neighbourhood of $\g$
$$\left\{
\g^\prime\st
\left\vert
\int_{\T^p\times\R^p}f(x-v)\g^\prime(x,v)\dr\mu(x)\dr v-
\int_{\T^p\times\R^p}f(x-v)\g(x,v)\dr\mu(x)\dr v
\right\vert  <\e
\right\}   $$
as
$$\left\{ 
\g^\prime\st
\left\vert
\int_{\T^p}f(z)\dr(\mu\ast\g^\prime)(z)-
\int_{\T^p}f(z)\dr(\mu\ast\g)(z)   
\right\vert <\e  
\right\}$$
by the definition of $\mu\ast\g^\prime$ and $\mu\ast\g$.

\noindent{\bf Step 2.} We note that the linear function
$$\fun{I_{pot}}{\g}{
\int_{\T^p\times\R^p}P^{
\2\mu
}   (t,x)\g(x,v)\dr\mu(x)\dr v
}    $$
does not depend on $\g$ by (1.3).  

\noindent{\bf Step 3.} We prove that
$$\fun{I_{gauss}}{\g}{
\int_{\T^p\times\R^p}A_h(\g,(x,v))\dr\mu(x)\dr v
}  $$
is weakly l. s. c.. Since $I_{gauss}$ is convex, it suffices to prove that it is l. s. c. for the strong topology of 
$L^1(\mu\otimes\L^p)$. We saw after formula (1.6) that
$$\frac{1}{2h}|v|^2\g(v)+\g(v)\log\g(v)\ge
e^{
-1-\frac{1}{2h}|v|^2
}       .   $$
Since the term on the right is integrable, lower semicontinuity follows from Fatou's lemma.

\fin

\prop{1.4} Let $U\in C(\Mt,\R)$ and let $\mu\in\Mt$; then, the $\inf$ in the definition of $(G^h_tU)(\mu)$ is a minimum.

\proof {\bf Step 1.} We begin to show that $(G^h_tU)(\mu)$ is finite.

If we substitute 
$$\g(x,v)=\left(\frac{1}{2\pi h}\right)^\frac{p}{2}e^{
-\frac{1}{2h}|v|^2
}  $$
into (1.4) (or (1.5), which is the same up to a constant), we immediately get that $(G^h_tU)(\mu)<+\infty$; to prove that 
$(G^h_tU)(\mu)>-\infty$, it suffices to prove that the functional $I$ defined in the last lemma is a sum of functions, each of which is bounded from below. 

First, the function bringing $\g$ in $U(\mu\ast\g)$, i. e. 
$$\fun{}{\g}{U(\mu\ast\g)}$$
is bounded from below because $U$, a continuous function on a compact space, is bounded from below.

Second, the functional $I_{pot}$ defined in the last lemma does not depend on $\g$ by (1.3); we have an explicit bound on its value since 
$$||P^{\2\mu}||_{C^4}\le M\colon=(
||V||_{C^4}+
||W||_{C^4} 
)  .    \eqno (1.13)$$

Third, 
$$I_{gauss}(\g)=
\int_{\T^p}\dr\mu(x)
\int_{\R^p}\left[
\frac{1}{2h}|v|^2+\log\g(x,v)
\right]    \g(x,v)\dr v$$
is bounded below by (1.6).

\noindent {\bf Step 2.} Let $\{ \g_n \}$ be a sequence minimizing in (1.4); we assert that, up to subsequences, 
$\g_n\weak\g\in\dcal_\mu$.

We begin to show that 
$$\int_{\T^p\times\R^p} A_h(\g_n,(x,v))\dr\mu(x)\dr v\le C  
\eqno (1.14)$$
for some $C>0$ independent on $n$.

Since $\{  \g_n  \}$ is minimizing in (1.4), by step 1 and (1.13) there is 
$C_1>0$ such that
$$\int_{\T^p\times\R^p}  A_h(\g_n,(x,v))  \dr\mu(x)\dr v+
U(\mu\ast\g_n)\le C_1\qquad\forall n\ge 1  .  $$
Now (1.14) follows by the fact that $U$, being a continuous function on the compact space $\Mt$, is bounded. 

By (1.14), we get that $\{ \g_n \}$ satisfies points 1) and 2) of lemma 1.2; it is well-known that this implies that $\{ \g_n \}$ is weakly compact in $L^1(\mu\otimes\L^p)$. The weak limit $\g$ belongs to ${\cal D}_\mu$ by point 3) of lemma 1.2. 

\noindent{\bf End of the proof.} By step 2, any minimizing sequence $\{ \g_n \}$ has a subsequence $\{ \g_{n_k} \}$ such that $\g_{n_k}\weak\g\in\dcal_\mu$. Since the function $I$ is l. s. c. by lemma 1.3, $\g$ is a minimizer and the thesis follows.

\fin

\vskip 2pc

\centerline{\bf \S 2}
\centerline{\bf The time step: properties of the minimal}

\vskip 1pc

In this section, we prove proposition 2.3 below, which says that the modulus of continuity of $G^h_tU$ is only slightly larger than the modulus of continuity of $U$; and proposition 2.8, which says that, if $\g$ is minimal, then  the $L^\infty$ norm of $\g$ (and that of 
$\frac{1}{\g}$ on $\T^p\times B(0,2\sqrt p)$) is bounded in terms of the Lipschitz constant of $U$.

We begin with a standard fact from [14].

\lem{2.1} Let $U_1,U_2\in C(\Mt)$. Then, the following three points hold.

\noindent 1) If $U_1\le U_2$, then $G^h_t U_1\le G^h_tU_2$.

\noindent 2) For all $a\in\R$, $G^h_t(U_1+a)=G^h_tU_1+a$.

\noindent 3) $||G^h_tU_1-G^h_tU_2||_\infty\le||U_1-U_2||_\infty$.

\proof Points 1) and 2) are immediate consequences of the definition of the operator $G^h_t$, i. e. of formula (1.4); point 3) follows from 1) and 2) in a standard way.

\fin

We need a technical fact, lemma 2.2. below, and some notation; the readers of [3] will recognize the "push forward by plans".

\vskip 1pc

\noindent{\bf Definition.} Let $\mu_0,\mu_1\in\Mt$, let $\G$ be a transfer plan between $\mu_0$ and $\mu_1$ and let 
$\g_0\in{\cal D}_{\mu_0}$.

Here and in the following, we shall always reserve the variable 
$x\in\T^p$ for integration in $\mu_0$, and $y\in\T^p$ for integration in $\mu_1$. 

We disintegrate $\G$ as $\G=\G_y\otimes\mu_1$ (see [5], II.70 for the precise statement and proof of the disintegration theorem) and we set
$$\g_1(y,v)=\int_{\T^p}\g_0(x,v)\dr\G_y(x)  .  \eqno (2.1)$$
Formula (2.1) is just a generalized way of composing $\g$ with a map; indeed, if $\G$ were induced by an invertible map $g$, then we would have
$$\g_1(y,v)=\g_0(g^{-1}(y),v)   .   $$

\lem{2.2} Let $\mu_0$, $\mu_1$, $\g_0$ and $\g_1$ be as in the definition above. Then
$$\g_1\in{\cal D}_{\mu_1}   \txt{and}   \leqno 1) $$
$$\int_{\T^p\times\R^p}A_h(\g_1,(y,v))\dr\mu_1(y)\dr v\le
\int_{\T^p\times\R^p}A_h(\g_0,(x,v))\dr\mu_0(x)\dr v  .
\leqno 2)    $$
Moreover, if $\G$ is a transfer plan on which $d_1(\mu_0,\mu_1)$ is attained, we have that
$$d_1(\mu_0\ast\g_0,\mu_1\ast\g_1)\le d_1(\mu_0,\mu_1)  .  
\leqno 3)$$ 

\proof The first equality below follows from (2.1), the second one from Fubini.  
$$\int_{\R^p}\g_1(y,v)\dr v=
\int_{\R^p}\dr v\int_{\T^p}\g_0(x,v)\dr\G_y(x)=
\int_{\T^p}\dr\G_y(x)\int_{\R^p}\g_0(x,v)\dr v  .  $$
Now recall that $\g_0\in{\cal D}_{\mu_0}$, and thus
$$\int_{\R^p}\g_0(x,v)\dr v=1
\txt{for $\mu_0$ a. e. $x$.}  $$
Since a $\mu_0$-null set is a $\G_y$-null set for $\mu_1$ a. e. 
$y$, the last two formulas imply that
$$\int_{\R^p}\g_1(y,v)\dr v=1
\txt{for $\mu_1$ a. e. $y$.}  $$
This proves point 1); we turn to point 2). The first equality below is (2.1); for the inequality, we consider the strictly convex function 
$\phi(z)=z\log z$ and apply Jensen.
$$A_h(\g_1,(y,v))=
\int_{\T^p}\cinh{v}\g_0(x,v)\dr\G_y(x)+
\int_{\T^p}\g_0(x,v)\dr\G_y(x)\log\int_{\T^p}\g_0(x,v)\dr\G_y(x)\le$$
$$\int_{\T^p}[
\cinh{v}\g_0(x,v)+\g_0(x,v)\log\g_0(x,v)
]   \dr\G_y(x)=\int_{\T^p}A_h(\g_0,(x,v))\dr\G_y(x)  .  $$
Since $\phi(z)=z\log z$ is strictly convex, equality holds if there is an invertible minimal transfer map. Integrating, we get the inequality below. 
$$\int_{\T^p\times\R^p}A_h(\g_1,(y,v))\dr\mu_1(y)\dr v\le
\int_{\T^p\times\R^p}\dr\mu_1(y)\dr v\int_{\T^p}
A_h(\g_0,(x,v))\dr\G_y(x)=$$
$$\int_{\T^p\times\T^p\times\R^p}A_h(\g_0,(x,v))\dr\G(x,y)\dr v=
\int_{\T^p\times\R^p}A_h(\g_0,(x,v))\dr\mu_0(x)\dr v  .    $$
The first equality above follows because 
$\G=\G_y\otimes\mu_1$, the second one because the first marginal of $\G$ is $\mu_0$. 

We prove 3). The first equality below is (1.1), while the second one is the definition of $\mu_1\ast\g_1$ and 
$\mu_0\ast\g_0$; the third one is the definition of $\g_1$ in (2.1); the fourth one follows from the fact that $\G=\G_y\otimes\mu_1$ and the marginals of $\G$ are $\mu_0$ and $\mu_1$.
$$d_1(\mu_1\ast\g_1,\mu_0\ast\g_0)=
\sup_{f\in Lip^1(\T^p)}\left\vert
\int_{\T^p}f(y)\dr(\mu_1\ast\g_1)(y)-
\int_{\T^p}f(x)\dr(\mu_0\ast\g_0)(x)
\right\vert   =   $$
$$\sup_{f\in Lip^1(\T^p)}\left\vert
\int_{\T^p\times\R^p}f(y-v)\g_1(y,v)\dr\mu_1(y)\dr v-
\int_{\T^p\times\R^p}f(x-v)\g_0(x,v)\dr\mu_0(x)\dr v
\right\vert   =   $$
$$\sup_{f\in Lip^1(\T^p)}\left\vert
\int_{\T^p\times\R^p}f(y-v)\dr\mu_1(y)\dr v
\int_{\T^p}\g_0(x,v)\dr\G_y(x)-
\int_{\T^p\times\R^p}f(x-v)\g_0(x,v)\dr\mu_0(x)\dr v
\right\vert   =   $$
$$\sup_{f\in Lip^1(\T^p)}\left\vert
\int_{\T^p\times\T^p\times\R^p}f(y-v)\g_0(x,v)\dr\G(x,y)\dr v-
\int_{\T^p\times\T^p\times\R^p}f(x-v)\g_0(x,v)\dr\G(x,y)\dr v
\right\vert   \le   $$
$$\int_{\T^p\times\T^p\times\R^p}
|x-y|_{\T^p}\g_0(x,v)\dr\G(x,y)\dr v   .   $$
Recalling that $\g_0$ satisfies (1.2) for $\mu_0$ a. e. $x\in\T^p$, the formula above yields the inequality below; the equality comes from the fact that $\G$ is a minimal transfer plan.
$$d_1(\mu_1\ast\g_1,\mu_0\ast\g_0)\le
\int_{\T^p\times\T^p}|x-y|_{\T^p}\dr\G(x,y)=
d_1(\mu_0,\mu_1)   .      $$

\fin

\noindent{\bf Definition.} Let $U\in C(\Mt)$; we say that 
$\fun{\o}{[0,+\infty)}{[0,+\infty)}$ is a 1-modulus of continuity for 
$U$ if

\noindent 1) $\o$ is concave.

\noindent 2) $\o(0)=0$.

\noindent 3) $|U(\mu_1)-U(\mu_0)|\le\o(d_1(\mu_1,\mu_2))$ for all
$\m_1,\mu_2\in\M_t$.

\prop{2.3} There is a constant $C>0$, depending only on the potentials $V$ and $W$, such that the following holds. Let $\o$ be a 1-modulus of continuity for $U\in C(\Mt)$; then, 
$\tilde\o(z)\colon=Chz+\o(z)$ is a 1-modulus of continuity for 
$G^h_tU$.

In particular, if $U$ is $L$-Lipschitz for the 1-Wasserstein distance, then $G^h_tU$ is $(Ch+L)$-Lipschitz; if $U$ is continuous, then $G^h_tU$ is continuous.

\proof We assert that it suffices to show the following: if 
$\mu_0,\mu_1\in\Mt$ and $\g_0$ minimizes
$$\fun{}{\g}{S(U,\mu_0,\g)}$$
($\g_0$ exists by proposition 1.4), then we can find 
$\g_1\in\dcal_{\mu_1}$ such that
$$S(U,\mu_1,\g_1)
+\int_{\T^p}P^{\2\mu_1}(t,x)\dr\mu_1(x)\le
S(U,\mu_0,\g_0)+
\int_{\T^p}P^{\2\mu_0}(t,x)\dr\mu_0(x)+
\tilde\o(d_1(\mu_1,\mu_0))  .  
\eqno (2.2)$$
Indeed, this implies by (1.4) that
$$(G^h_tU)(\mu_1)\le(G^h_tU)(\mu_0)+
\tilde\o(d_1(\mu_1,\mu_0))   .  $$
Exchanging the r\^oles of $\mu_1$ and $\mu_0$, we get that 
$\tilde\o$ is a modulus of continuity for $G^h_tU$, and the assertion follows.

To prove (2.2), we let $\G$ be a minimal transfer plan between 
$\mu_0$ and $\mu_1$ and define $\g_1$ by (2.1). Now, (1.3) implies the equality below.
$$\int_{\T^p\times\R^p}
[hL^{\2\mu_1}(t,y,\frac{1}{h}v)+\log\g_1(y,v)]
\g_1(y,v)\dr\mu_1(y)\dr v-$$
$$\int_{\T^p\times\R^p}
[hL^{\2\mu_0}(t,x,\frac{1}{h}v)+\log\g_0(x,v)]
\g_0(x,v)\dr\mu_0(x)\dr v 
=$$
$$\int_{\T^p\times\R^p}A_h(\g_1,(y,v))\dr\mu_1(y)\dr v-
\int_{\T^p\times\R^p}A_h(\g_0,(x,v))\dr\mu_0(x)\dr v
   +   \eqno (2.3)_a$$
$$h\int_{\T^p}V(t,x)\dr\mu_0(x)-
h\int_{\T^p}V(t,y)\dr\mu_1(y)
+    \eqno (2.3)_b  $$
$$h\int_{\T^p}W^{\2\mu_0}(x)\dr\mu_0(x)-
h\int_{\T^p}W^{\2\mu_1}(y)
\dr\mu_1(y)  .   \eqno (2.3)_c  $$
Let us tackle the terms $(2.3)_a$, $(2.3)_b$ and $(2.3)_c$; first of all, point 2) of lemma 2.2 implies that
$$(2.3)_a\le 0  .  \eqno (2.4)$$
As for the term $(2.3)_b$, we have that
$$(2.3)_b=
h\int_{\T^p\times\T^p}[
V(t,x)-V(t,y)
]  \dr\G(x,y)  \le
h\int_{\T^p\times\T^p}
C_1|x-y|_{\T^p}\dr\G(x,y)=C_1hd_1(\mu_0,\mu_1)  .  
\eqno (2.5)$$
The first equality above follows because the marginals of $\G$ are $\mu_0$ and $\mu_1$; the inequality follows because $V$ is 
$C_1$-Lipschitz. The last equality follows from the fact that $\G$ is minimal in the definition of $d_1(\mu_0,\mu_1)$. 

Analogously, we get that
$$(2.3)_c=h\int_{\T^p\times\T^p}[
W^{\2\mu_0}(x)-W^{\2\mu_1}(y)
]  \dr\G(x,y)\le  $$
$$h\int_{\T^p\times\T^p}
|
W^{\2\mu_0}(x)-W^{\2\mu_0}(y)
|  \dr\G(x,y)+
h\int_{\T^p\times\T^p}
|
W^{\2\mu_1}(y)-W^{\2\mu_0}(y)
|  \dr\G(x,y)   .  $$
We can see as in (1.13) that the Lipschitz constant of 
$W^{\2\mu_0}$ is bounded by one half of the Lipschitz constant 
$C_2$ of $W$; this implies as in (2.5) that
$$\int_{\T^p\times\T^p}
|
W^{\2\mu_0}(x)-W^{\2\mu_0}(y)
|  \dr\G(x,y)   \le   \2  C_2 d_1(\mu_0,\mu_1)   .   $$
On the other hand, since $W$ is $C_2$-Lipschitz, (1.1) yields the inequality below; the equality comes from the definition of 
$W^{\2\mu_i}$.
$$|
W^{\2\mu_1}(y)-W^{\2\mu_0}(y)
|   =$$
$$\left\vert
\2\int_{\T^p}W(x-y)\dr\mu_1(x)-
\2\int_{\T^p}W(x-y)\dr\mu_0(x)
\right\vert    \le  
\2C_2 d_1(\mu_1,\mu_0)   .  $$
By the last three formulas, we get that
$$(2.3)_c\le C_2h d_1(\mu_0,\mu_1)   .   $$
From (2.3), (2.4), (2.5) and the last formula, we get that
$$\int_{\T^p\times\R^p}\left[
hL_c^{\2\mu_1}(t,y,\frac{1}{h}v)\g_1(y,v)
+\g_1(y,v)\log\g_1(y,v)
\right]   \dr\mu_1(y)\dr v-$$
$$\int_{\T^p\times\R^p}  \left[
hL_c^{\2\mu_0}(t,x,\frac{1}{h}v)\g_0(x,v)+
\g_0(x,v)\log\g_0(x,v)
\right]
\dr\mu_0(x)\dr v   \le  
C_3h d_1(\mu_0,\mu_1)   .   $$
Since $\o$ is the modulus of continuity of $U$, point 3) of lemma 2.2 implies that
$$|
U(\mu_1\ast\g_1)-U(\mu_0\ast\g_0)
|     \le   \o(d_1(\mu_0,\mu_0))    .    $$
Setting $C=C_3$, (2.2) is implied by the last two formulas.

\fin

We begin the estimate on $||\g||_\infty$ with a technical lemma.

\lem{2.4} Let $\mu\in\Mt$, and let $\g_0,\g_1\in{\cal D}_\mu$. Let us suppose that the functions $\g_0(x,\cdot)$ and $\g_1(x,\cdot)$ coincide whenever $x$ does not belong to a Borel set 
$E\subset\T^p$. Then,
$$d_1(\mu\ast\g_0,\mu\ast\g_1)\le
\sqr p 
\int_{E\times\R^p}|
\g_0(x,v)-\g_1(x,v)
|     \dr\mu(x)\dr v   .   $$

\proof We use the dual formulation (1.1) for the first equality below; the second one is the definition of $\mu\ast\g_i$; the third one follows because $\g_0$ and $\g_1$ coincide on 
$E^c\times\R^p$; the inequality comes from the fact that, since 
$f\in Lip^1(\T^p)$ and 
$\T^p$ has diameter $\sqr p$, we can as well suppose that 
$||f||_\infty\le\sqr p$.
$$d_1(\mu\ast\g_0,\mu\ast\g_1)=
\sup_{f\in Lip^1(\T^p)}\left[
\int_{\T^p}f(z)\dr(\mu\ast\g_0)(z)-\int_{\T^p}f(z)\dr(\mu\ast\g_1)(z)
\right]   =   $$
$$\sup_{f\in Lip^1(\T^p)}\left[
\int_{\T^p\times\R^p}
f(x-v)[
\g_0(x,v)-\g_1(x,v)
]\dr\mu(x)\dr v
\right]      =$$
$$\sup_{f\in Lip^1(\T^p)}\left[
\int_{E\times\R^p}
f(x-v)[
\g_0(x,v)-\g_1(x,v)
]\dr\mu(x)\dr v
\right]         \le
\sqr p \int_{E\times\R^p}  |
\g_0(x,v)-\g_1(x,v)
|    \dr\mu(x)\dr v    .   $$

\fin

\lem{2.5} There is a constant 
$C_1(L,h)$, depending only on $L,h>0$, for which the following happens. Let $U$ be $L$-Lipschitz for the 1-Wasserstein distance $d_1$, let $\mu\in\Mt$ and let 
$\g$ minimize in the definition of $(G^h_tU)(\mu)$. Then,
$$||\g||_{L^\infty(\T^p\times\R^p,\mu\otimes\L^p)}     \le C_1(L,h)     .   $$

\proof  We are going to use the fact that the superlinear entropy term becomes huge when $\g$ is large; thus, if $|| \g ||_\infty$ is too large, we can take some mass from the region where $\g$ is big, smear it where $\g$ is small and obtain a function $\tilde\g$ such that 
$$S(U,\mu,\tilde\g)< S(U,\mu,\g)  ,   $$ 
contradicting the minimality of $\g$.  

\noindent{\bf Step 1.} We define the set where $\g$ is large.

Let us suppose that 
$||\g||_{L^\infty(\mu\otimes\L^p)}\ge 2A$; then, there is a Borel set 
$D_A\subset\T^p\times\R^p$ such that
$$0<(\mu\otimes\L^p)(D_A)   \txt{and}
\g(x,v)\ge A\qquad   \forall (x,v)\in D_A  .   \eqno (2.6)$$
We denote by $D_A(x)$ its sections:
$$D_A(x)\colon=\{
v\in\R^p\st (x,v)\in D_A
\}   .  $$
Since $\g\in{\cal D}_\mu$, Chebishev's inequality implies that 
$\L^p(D_A(x))\le\frac{1}{A}$ for $\mu$ a. e. $x$.

We set $2a(p)=\L^p(B(0,1))$ and we define
$$B_A(x)=\int_{D_A(x)}\left(
\g(x,v)-a(p)
\right)     \dr v      .    $$
We shall suppose that $A>\max(a(p),a(p)^{-1})$ (otherwise there is nothing to prove); as a consequence, $B_A(x)\ge 0$. 
Since $\g\in{\cal D}_\mu$, we have that $B_A(x)\in[0,1]$ for 
$\mu$ a. e. $x$.

\noindent {\bf Step 2.} We show that set where $\g$ is small has room enough to accommodate some mass from $D_A$.

We let
$$Z=\{
(x,v)\in\T^p\times B(0,1)\st\g(x,v)\le a(p)^{-1}
\}    .    $$
As above, we call $Z(x)$ its sections. Since 
$\int_{\R^p}\g(x,v)\dr v=1$ for $\mu$ a. e. $x$, we get by the Chebishev inequality that $\L^p(B(0,1)\setminus Z(x))\le a(p)$ for $\mu$ a. e. $x$. Since $\L^p(B(0,1))=2a(p)$, this implies that 
$\L^p(Z(x))\ge a(p)$ for $\mu$ a. e. $x$.

A standard consequence of this is that we can find a Borel set 
$\tilde Z\subset Z$ such that $\mu$ a. e. section $\tilde Z(x)$ satisfies 
$\L^p(\tilde Z(x))=a(p)$. 

\noindent{\bf Step 3.} We build $\tilde\g$.

Since we have chosen $A> a(p)^{-1}$, we have that $\tilde Z(x)$ and $D_A(x)$ are disjoint; we set $M(x)=\R^p\setminus(\tilde Z(x)\cup B_A(x))$ and 
$$\tilde\g(x,v)=\left\{
\eqalign{
\g(x,v)   &\txt{if} v\in M(x)\cr
\g(x,v)+ B_A(x)a(p)^{-1}
&\txt{if}  v\in \tilde Z(x)\cr
a(p)  &\txt{if} v\in D_A(x)   .  
}          \right.    $$
The first equality below comes from the fact that $\g$ and 
$\tilde\g$ coincide on $M(x)$, the second one from the fact that 
$\L^p(\tilde Z(x))=a(p)$ and the third one from the definition of $B_A(x)$.
$$\int_{\T^p\times\R^p}|
\g(x,v)-\tilde\g(x,v)
|   \dr\mu(x)\dr v   =  $$
$$\int_{\T^p}\dr\mu(x)\int_{\tilde Z(x)}a(p)^{-1} B_A(x)\dr v +
\int_{\T^p}\dr\mu(x)\int_{D_A(x)}[
\g(x,v)-a(p)
]   \dr v=     $$
$$\int_{\T^p}B_A(x)\dr\mu(x)+
\int_{\T^p}\dr\mu(x)\int_{D_A(x)}[
\g(x,v)-a(p)
]   \dr v  = 
2 \int_{\T^p}B_A(x)\dr\mu(x)  .   \eqno (2.7)$$
The same argument without the modulus shows the first equality below; the second one follows since $\g\in{\cal D}_\mu$.
$$\int_{\R^p}\tilde\g(x,v)\dr v=
\int_{\R^p}\g(x,v)\dr v=1
\txt{for} \mu \txt{a. e.} x   .   $$
In other words, $\tilde\g\in{\cal D}_\mu$. 
In order to compare the actions of $\g$ and $\tilde\g$ we note that, for all $x$,
$$\int_{\R^p}
\left[
\cinh{v}+\log\tilde\g(x,v)
\right]    \tilde\g(x,v)\dr v=
\int_{M(x)}
\left[
\cinh{v}+\log\g(x,v)
\right]    \g(x,v)\dr v+$$
$$\int_{\tilde Z(x)}
\left[
\cinh{v}+\log\tilde\g(x,v)
\right]    \tilde\g(x,v)\dr v+
\int_{D_A(x)}
\left[
\cinh{v}+\log\tilde\g(x,v)
\right]    \tilde\g(x,v)\dr v    \eqno (2.8)$$
because $\g$ and $\tilde\g$ coincide on $M(x)$.

\noindent{\bf Step 4.} We compare the actions of $\g$ and 
$\tilde\g$ on $\tilde Z(x)$.

If $v\in\tilde Z(x)$, then $\g(x,v)\le a(p)^{-1}$; this yields the first  inequality below; for the second one, we recall that, by step 1,  $B_A(x)\le 1$.
$$\tilde\g(x,v)=\g(x,v)+B_A(x) a(p)^{-1}\le
a(p)^{-1}(1+B_A(x))\le 2a(p)^{-1} 
\qquad\forall v\in\tilde Z(x)   .  \eqno (2.9)$$
The inequality below follows because 
$\tilde\g(x,v)\ge\g(x,v)$ and $\cinh{v}\le\frac{1}{2h}$ on 
$\tilde Z(x)\subset B(0,1)$; moreover, we have used the Lagrange inequality and the fact, which follows from (2.9), that 
$[\g(x,v),\tilde\g(x,v)]\subset [0,2 a(p)^{-1}]$. The second equality follows by the definition of $\tilde\g$ on $\tilde Z(x)$, and the fact that 
$\L^p(\tilde Z(x))=a(p)$.
$$\int_{\tilde Z(x)}
\left[
\cinh{v}+\log\tilde\g(x,v)
\right]    \tilde\g(x,v)\dr v=
\int_{\tilde Z(x)}
\left[
\cinh{v}+\log\g(x,v)
\right]    \g(x,v)\dr v+ $$
$$\int_{\tilde Z(x)}
\cinh{v}[
\tilde\g(x,v)-\g(x,v)
]    \dr v+
\int_{\tilde Z(x)}[
\tilde\g(x,v)\log\tilde\g(x,v)-\g(x,v)\log\g(x,v)
]    \dr v\le$$
$$\le
\int_{\tilde Z(x)}
\left[
\cinh{v}+\log\g(x,v)
\right]    \g(x,v)\dr v+ $$
$$\int_{\tilde Z(x)}\frac{1}{2h}[
\tilde\g(x,v)-\g(x,v)
]     \dr v+
\max_{t\in[0,2a(p)^{-1}]}\frac{\dr}{\dr t}(t\log t)\cdot
\int_{\tilde Z(x)}[
\tilde\g(x,v)-\g(x,v)
]      \dr v=$$
$$\int_{\tilde Z(x)}
\left[
\cinh{v}+\log\g(x,v)
\right]    \g(x,v)\dr v+
\frac{1}{2h} B_A(x)+[
\log(2a(p)^{-1})+1
]B_A(x)   .   
\eqno (2.10)$$

\noindent {\bf Step 5.} We compare the actions of $\g$ and 
$\tilde\g$ on $D_A(x)$.

Since $\fun{}{t}{t\log t}$ is superlinear, there is $M(A)$, tending to 
$+\infty$ as $A\tends+\infty$, such that 
$$\g(x,v)\log\g(x,v)-\tilde\g(x,v)\log\tilde\g(x,v)\ge
M(A)[
\g(x,v)-a(p)
]    \qquad\forall v\in D_A(x)       .  $$
This implies the inequality below, while the last equality comes from the definition of $B_A(x)$. 
$$\int_{D_A(x)}\tilde\g(x,v)\log\tilde\g(x,v)\dr v=$$
$$\int_{D_A(x)}\g(x,v)\log\g(x,v)\dr v+
\int_{D_A(x)}[\tilde\g(x,v)\log\tilde\g(x,v)-
\g(x,v)\log\g(x,v)]\dr v\le  $$
$$\int_{D_A(x)}\g(x,v)\log\g(x,v)\dr v  -
M(A) \int_{D_A(x)}[\g(x,v)-a(p)]\dr v=
\int_{D_A(x)}\g(x,v)\log\g(x,v)\dr v  -
M(A) B_A(x)   .  $$
The first inequality below follows from the fact that $\g\ge\tilde\g$ on $D_A(x)$; the second one, from the last formula.
$$\int_{D_A(x)}
\left[
\cinh{v}+\log\tilde\g(x,v)
\right]    \tilde\g(x,v)\dr v
\le 
\int_{D_A(x)}\cinh{v}\g(x,v)\dr v+
\int_{D_A(x)}\tilde\g(x,v)\log\tilde\g(x,v)\dr v\le$$
$$\int_{D_A(x)}\cinh{v}\g(x,v)\dr v+
\int_{D_A(x)}\g(x,v)\log\g(x,v)\dr v-   
M(A)B_A(x)    .   $$

\noindent{\bf End of the proof.} From the last formula, (2.8) and (2.10) we get that
$$\int_{\R^p}
\left[
\cinh{v}+\log\tilde\g(x,v)
\right]    \tilde\g(x,v)\dr v\le
\int_{\R^p}\left[
\cinh{v}+\log\g(x,v)
\right]    \g(x,v)\dr v+$$
$$\left(
\frac{1}{2h}+1+\log(2a(p)^{-1}) -M(A)
\right)      
B_A(x)  .  \eqno (2.11)$$
The first inequality below follows by lemma 2.4; the second one, from (2.7).
$$d_1(\mu\ast\g,\mu\ast\tilde\g)\le
\sqrt p\int_{\T^p\times\R^p}|
\g(x,v)-\tilde\g(x,v)
|\dr\mu(x)\dr v\le
2\sqrt p\int_{\T^p} B_A(x)\dr\mu(x)  .  $$
Since $U$ is $L$-Lipschitz, this implies
$$|
U(\mu\ast\g)-U(\mu\ast\tilde\g)
|    \le 2\sqrt p L\int_{\T^p}B_A(x)\dr\mu(x)  .  $$
By the last formula and (2.11) we get that
$$\int_{\T^p\times\R^p}A_h(\tilde\g,(x,v)) \dr\mu(x)\dr v+
U(\mu\ast\tilde\g)\le
\int_{\T^p\times\R^p}A_h(\g,(x,v))\dr\mu(x)\dr v+
U(\mu\ast\g)+$$
$$\left(
\frac{1}{2h}+1+\log 2a(p)^{-1}+
2\sqrt p L   -M(A)
\right)    
\int_{\T^p}B_A(x)\dr\mu(x)    .  $$
We have seen that $M(A)\tends+\infty$ as $A\tends+\infty$; thus, if $A$ is large enough, we have contradicted the minimality of 
$\g$.

\fin

Before proving the estimate from below on $\g$, we need a result on the tightness on the set of all minimal $\g$'s.

\lem{2.6} Let $U$ be $L$-Lipschitz. Then, for all $\e>0$ there is $R>0$ such that the following happens. If $\mu\in\Mt$ and $\g$ minimizes the single particle functional $S(U,\mu,\g)$, then for 
$\mu$ a. e. $x\in\T^p$ we have that
$$\int_{B(0,R)^c}\g(x,v)\dr v\le\e  .  $$

\proof Let us suppose by contradiction that the thesis does not hold; then there is $\e>0$ such that for infinitely many $l\in\N$ we can find 

\noindent 1) a measure $\mu_l\in\Mt$, 

\noindent 2) a minimal $\g_l$ and 

\noindent 3) a set $E_l\subset\T^p$ with 
$\mu_l(E_l)>0$ such that 
$$\int_{B(0,l)^c}\g(x,v)\dr v>\e\txt{for $\mu$-a. e. }x\in\T^p  .  $$ 

Point 3) implies the second inequality below. 
$$\inf_{x\in E_l}
\int_{\R^p}\cinh{v}\g_l(x,v)\dr v\ge   \inf_{x\in E_l}
\int_{B(0,l)^c}\cinh{v}\g_l(x,v)\dr v\ge
\frac{\e}{2h}|l|^2      \tends+\infty   .  $$
Then, (1.7) implies that
$$\inf_{x\in E_l}
\int_{\R^p}A_h(\g_l,(x,v))\dr v
\ge M_l     \eqno (2.12)$$
with $M_l\tends+\infty$ as $l\tends+\infty$.

We set
$$\tilde\g_l(x,v)=
\left\{
\eqalign{
\g_l(x,v)&\txt{if}x\not\in E_l\cr
\left(
\frac{1}{2\pi h}
\right)^\frac{p}{2}
e^{
-\frac{1}{2h}|v|^2
}   &\txt{if}x\in E_l    .
}
\right.   $$
Since $\g_l(x,\cdot)$ and the Gaussian have integral one over 
$\R^p$, we have that $\tilde\g_l\in{\cal D}_\mu$; by (2.12), we get that
$$\int_{\R^p}A(\g_l,(x,v))\dr v-
\int_{\R^p}A(\tilde\g_l,(x,v))\dr v\ge M_l^\prime \txt{for}   
x\in E_l$$
with $M_l^\prime\tends+\infty$ as $l\tends+\infty$. Integrating over $\T^p$ and recalling the definition of $\tilde\g_l$, we get that 
$$\int_{\T^p\times\R^p}A(\tilde\g_l,(x,v))\dr\mu_l(x)\dr v\le
\int_{\T^p\times\R^p}A(\g_l,(x,v))\dr\mu_l(x)\dr v-
M_l^\prime  \mu_l(E_l)   .  $$
Lemma 2.4 yields the first inequality below; since 
$\g,\tilde\g\in{\cal D}_{\mu_l}$, the second one follows.
$$d_1(\mu_l\ast\g_l,\mu_l\ast\tilde\g_l)\le
\sqrt{p}\int_{E_l}\dr\mu_l(x)
\int_{\R^p}|
\tilde\g_l(x,v)-\g_l(x,v)
|    \dr v\le
2\sqrt{p}\mu_l(E_l)    .    $$
From the last two formulas and the fact that $U$ is $L$-Lipschitz, we get that
$$\int_{\T^p\times\R^p}A_h(\tilde\g_l,(x,v))\dr\mu(x)\dr v+
U(\mu_l\ast\tilde\g_l)\le$$
$$\int_{\T^p\times\R^p}A_h(\g_l,(x,v))\dr\mu(x)\dr v+
U(\mu_l\ast\g_l)   +(2L\sqrt{p}-M^\prime_l)\mu_l(E_l)   .   $$
Since $M^\prime_l\tends+\infty$ and $\mu_l(E_l)>0$, if we take 
$l$ large enough we contradict  the minimality of $\g_l$.

\fin

\lem{2.7} There is a constant $C_2(L,h)$, depending only on $L,h>0$, for which the following happens. Let $U$ be 
$L$-Lipschitz, let $\mu\in\Mt$ and let 
$\g$ minimize in the definition of $(G^h_tU)(\mu)$. Then,
$$||\frac{1}{\g}||_{L^\infty(\T^p\times B(0,2\sqrt p),\mu\otimes\L^p)}     \le 
C_2(L,h)     .   $$

\proof Reversing the procedure of the lemma 2.5, we add some mass to the region $D_A$ where $\g$ is small, taking it from the region $Z_\d$ where it is larger; some work (step 2 below) is necessary to check that $(\mu\otimes\L^p)(Z_\d)$ is large enough. 

\noindent{\bf Step 1.} We settle the notation. 

We begin with the set where $\g$ is small. Let us suppose that 
$||\frac{1}{\g}||_{L^\infty(\T^p\times B(0,2\sqrt{p}))}\ge 
\frac{2}{A}$; we define
$$D_A=\{
(x,v)\st v\in B(0,2\sqrt{p})\txt{and} \g(x,v)\le A
\}   .  $$
Clearly, $(\mu\otimes\L^p)(D_A)>0$. 
As in lemma 2.5, we set $2a(p)=\L^p(B(0,1))$; we define
$$B_A(x)=\int_{D_A(x)}\g(x,v)    \dr v   .  $$
Let us set $P=(2\sqrt{p})^p$. 
Since $D_A(x)\subset B(0,2\sqrt{p})$, we have that 
$B_A(x)\le 2a(p)A P$ for $\mu$ a. e. $x$. 

Now we define a set of points which are not too far from the origin and where $\g$ is not too small. For $\d>0$, we define 
$$Z_\d=\{
(x,v)\in\T^p\times  B(0,\frac{1}{(8a(p)\d)^\frac{1}{p}})\st
\d\le\g(x,v)
\}  $$
and we call $Z_\d(x)$ its sections.

\noindent{\bf Step 2.} By lemma 2.6, we can find $l\in\N$ such that 
$$\int_{B(0,l)}\g(x,v)\dr v\ge\frac{1}{2}
\txt{for $\mu$ a. e. } x\in\T^p  .  \eqno (2.13)$$

We want to exclude the possibility that the mass of (2.13) is concentrated on a set of very small measure.

More precisely, we assert that, for all $\d>0$ small enough and independent on $\mu$, any minimal $\g$ satisfies 
$\L^p(Z_\d(x))>\frac{\d^p}{8a(p)}$ for $\mu$ a. e. $x$. 

Indeed, if this were not the case, for all $k\ge 1$ we could find 
$\mu_k\in\Mt$, a minimal $\g_k$  and a set $E_k\subset\T^p$ with $\mu_k(E_k)>0$ such that, for all $x\in E_k$ we have
$\L^p(Z_\frac{1}{k}(x))\le\frac{1}{8a(p)k^p}$.

Formula (2.13) implies that, for $k$ large enough, the first inequality below holds; the first equality is the definition of 
$Z_\frac{1}{k}(x)$, while the last one comes from the fact  that $2a(p)$ is the measure of the unit ball. 
$$\int_{Z_\frac{1}{k}(x)}\g_k(x,v)\dr v=$$
$$\int_{
\{ v\in B(0,\left(\frac{k}{8a(p)}\right)^\frac{1}{p})\st 
\g_k(x,v)\ge\frac{1}{k} \}
}    
\g_k(x,v)\dr v\ge
\2-
\int_{
\{ v\in B(0,\left(\frac{k}{8a(p)}\right)^\frac{1}{p})\st 
\g_k(x,v)<\frac{1}{k} \}
}    
\g_k(x,v)\dr v\ge$$
$$\2-\int_{B(0,\left(\frac{k}{8a(p)}\right)^\frac{1}{p})}
\frac{1}{k}\dr v=\frac{1}{4}\qquad
\forall x\in E_k   .  $$
In other words, the integral of $\g_k(x,\cdot)$ over 
$Z_\frac{1}{k}(x)$ is larger than $\frac{1}{4}$, while we are supposing that 
$\L^p(Z_\frac{1}{k}(x))\le\frac{1}{8a(p)k^p}$; this implies that 
$\{ \g_k(x_k,\cdot) \}_k$ is not uniformly integrable, however we choose the sequence $x_k\in E_k$. Using this and arguing as in point 1) of lemma 1.2, we get that
$$\int_{
\R^p
}       
A_h(\g_k,(x,v))\dr v
\ge M_k  \qquad
\forall x\in E_k    $$
with $M_k\tends+\infty$ as $k\tends+\infty$. If we define
$$\hat\g_k(x,v)=
\left\{
\eqalign{
\g_k(x,v)&\txt{if}x\not\in E_k\cr
\left(
\frac{1}{2\pi h}
\right)^\frac{p}{2}
e^{
-\frac{1}{2h}|v|^2
}   &\txt{if}x\in E_k    
}
\right.   $$
it is easy to see, using the last formula and arguing as in lemma 2.6, that, for $k$ large, 
$$S(U,\mu,\hat\g_k)<S(U,\mu,\g_k)$$
contradicting the minimality of $\g$.

\noindent{\bf Step 3.} Here we build the function $\tilde\g$; we shall show in the next two steps that, if $A$ is small enough, the action of $\tilde\g$ is lower than the optimal $\g$; this contradiction will end the proof.

Let us fix $\d>0$ such that step 2 holds; we shall suppose that 
$A<\d$ (otherwise there is nothing to prove); with this choice, 
$D_A(x)$ and $Z_\d(x)$ are disjoint. The first inequality below follows from the fact, shown in step 1, that $B_A(x)\le 2a(P)AP$; the second one, from the fact that $A<\d$; we choose $\e$ so small that the third one also holds.
$$\d-\e B_A(x)\frac{8a(p)}{\d^p}\ge
\d-\frac{16\e a(p)^2AP}{\d^p}\ge
\d-\frac{16\e a(p)^2\d P}{\d^p}\ge 
\frac{\d}{2}   .    \eqno (2.14)$$

By step 2, we can find $\tilde Z\subset Z$ such that, for $\mu$ a. e. $x$, $\L^p(\tilde Z(x))=\frac{\d^p}{8a(p)}$. We define 
$M(x)=\R^p\setminus(\tilde Z(x)\cup D_A(x))$ and we set
$$\tilde\g(x,v)=\left\{
\eqalign{
\g(x,v)&\txt{if}v\in M(x)\cr
\g(x,v)-\e B_A(x)\frac{8a(p)}{\d^p}&\txt{if} v\in\tilde Z(x)\cr
(1+\e)\g(x,v)&\txt{if}v\in D_A(x)   .  
}
\right.     $$
We have to prove that $\tilde\g(x,v)\in{\cal D}_{\mu}$; we begin to show that, if $\e$ is small enough, $\tilde\g(x,v)\ge 0$; by the definition above, it suffices to prove that $\tilde\g(x,v)\ge 0$ when $v\in\tilde Z(x)$. Since $v\in\tilde Z(x)$, we get the first inequality below; the second one is (2.14).
$$\g(x,v)- \e B_A(x)\frac{8a(p)}{\d^p}\ge
\d- \e B_A(x)\frac{8a(p)}{\d^p}\ge\frac{\d}{2}   .    \eqno (2.15)   $$
This shows that $\tilde\g(x,v)\ge 0$. We prove that 
$$\int_{\R^p}\tilde\g(x,v)\dr v=1\txt{for $\mu$ a. e.} x\in\T^p  .  $$
Since $\g$ and $\tilde\g$ coincide on $M(x)$, we have the first equality below; the second one comes from the definition of $\tilde\g$ and the fact that $\L^p(\tilde Z(x))=\frac{\d^p}{8a(p)}$; the last one, from the definition of $B_A(x)$.
$$\int_{\R^p}[
\tilde\g(x,v)-\g(x,v)
]    \dr v=$$
$$\int_{\tilde Z(x)}[
\tilde\g(x,v)-\g(x,v)
]    \dr v    +
\int_{D_A(x)}[
\tilde\g(x,v)-\g(x,v)
]    \dr v     =$$
$$-\e B_A(x)\frac{8a(p)}{\d^p}\cdot\frac{\d^p}{8a(p)}+
\e\int_{D_A(x)}\g(x,v)\dr v   =0\txt{for $\mu$ a. e. $x$.}    $$
Since $\g\in{\cal D}_\mu$, this ends the proof that 
$\tilde\g\in{\tilde D}_\mu$. With the same argument,
$$\int_{\T^p\times\R^p}|
\g(x,v)-\tilde\g(x,v)
|\dr\mu(x)\dr v\le
2\e\int_{\T^p}B_A(x)\dr\mu(x)  .  $$

\noindent{\bf Step 4.} We compare $U(\mu\ast\g)$ with 
$U(\mu\ast\tilde\g)$; actually, by lemma 2.4 and the last formula, we get that
$$|
U(\mu\ast\g)-U(\mu\ast\tilde\g)
|    \le
2\sqrt{p}L\e\int_{\T^p}B_A(x)\dr\mu(x)   .   \eqno (2.16)$$

\noindent{\bf Step 5.} We compare Lagrangian actions on 
$\tilde Z(x)$.

We recall that, if $v\in\tilde Z(x)$, then

\noindent 1) $\tilde\g(x,v)\le\g(x,v)$ and

\noindent 2) the derivative of $t\log t$ on $[\tilde\g(x,v),\g(x,v)]$ is greater than $(1-\log\frac{2}{\d})$; indeed, by (2.15),  
$[\tilde\g(x,v),\g(x,v)]\subset [\frac{\d}{2},+\infty)$.

Point 1) yields the first inequality below, point 2) the second one; the last equality follows from the fact that 
$\L^p(\tilde Z(x))=\frac{\d^p}{8a(p)}$. 
$$\int_{\tilde Z(x)}\left[
\cinh{v}+\log\tilde\g(x,v)
\right]     \tilde\g(x,v)\dr v=
\int_{\tilde Z(x)}\left[
\cinh{v}+\log \g(x,v)
\right]    \g(x,v)\dr v+$$
$$\int_{\tilde Z(x)}
\cinh{v}[\tilde\g(x,v)-\g(x,v)]\dr v+
\int_{\tilde Z(x)}[
\tilde\g(x,v)\log\tilde\g(x,v)-\g(x,v)\log\g(x,v)
]    \dr v   \le$$
$$\int_{\tilde Z(x)}\left[
\cinh{v}+\log\g(x,v)
\right]     \g(x,v)\dr v-
\inf_{t\in[\tilde\g(x,v),\g(x,v)]}\frac{\dr}{\dr t}(t\log t)
\int_{\tilde Z(x)}[\g(x,v)-\tilde\g(x,v)]\dr v\le$$
$$\int_{\tilde Z(x)}\left[
\cinh{v}+\log\g(x,v)
\right]     \g(x,v)\dr v+
(
-1+\log\frac{2}{\d}
)    \int_{\tilde Z(x)}  \e B_A(x)\frac{8a(p)}{\d^p}\dr v=       $$
$$\int_{\tilde Z(x)}\left[
\cinh{v}+\log\g(x,v)
\right]     \g(x,v)\dr v+
(
-1+\log\frac{2}{\d}
)    \e B_A(x)    .   \eqno (2.17)$$

\noindent{\bf Step 6.} We compare Lagrangian actions on 
$D_A(x)$.

With the same calculations as in step 5, and using the fact that the derivative of $A\log A$ tends to $-\infty$ as $A\searrow 0$,
$$\int_{D_A(x)}\left[
\cinh{v}+\log\tilde\g(x,v)
\right]     \tilde\g(x,v)\dr v=$$
$$\int_{D_A(x)}\left[
\cinh{v}+\log \g(x,v)
\right]    \g(x,v)\dr v+
\int_{D_A(x)}
\cinh{v}[\tilde\g(x,v)-\g(x,v)]\dr v+$$
$$\int_{D_A(x)}[
\tilde\g(x,v)\log\tilde\g(x,v)-\g(x,v)\log\g(x,v)
]    \dr v   \le
\int_{D_A(x)}\left[
\cinh{v}+\log\g(x,v)
\right]     \g(x,v)\dr v   +$$
$$\int_{D_A(x)}
\cinh{v}[\tilde\g(x,v)-\g(x,v)]\dr v
- M(A)\int_{D_A(x)}[\tilde\g(x,v)-\g(x,v)]\dr v$$
for a constant $M(A)\tends+\infty$ as $A\searrow 0$. Since 
$D_A(x)\subset B(0,2\sqrt p)$, we get that 
$\frac{1}{2h}|v|^2\le\frac{1}{2h}4p$ on $D_A(x)$; this and the last formula yield the first inequality below, while the equality comes by the definition of $\tilde\g$ and $B_A(x)$.
$$\int_{D_A(x)}\left[
\cin{v}+\log\tilde\g(x,v)
\right]   \tilde\g(x,v)\dr v\le$$
$$\int_{D_A(x)}\left[
\cinh{v}+\log\g(x,v)
\right]  \g(x,v)\dr v+
\int_{D_A(x)}\frac{1}{2h}4p\e\g(x,v)\dr v-
M(A)\int_{D_A(x)}\e\g(x,v)\dr v=$$
$$\int_{D_A(x)}\left[
\cin{v}+\log\g(x,v)
\right]   \g(x,v)\dr v+
\frac{1}{2h}4p\e B_A(x)-\e M(A)B_A(x)  .  \eqno (2.18)$$

\noindent{\bf End of the proof.} By (2.16), (2.17), (2.18) and the fact that Lagrangian actions on $M(x)$ coincide, we get that
$$\int_{\T^p\times\R^p}A_h(\tilde\g,(x,v))
\dr\mu(x)\dr v+U(\mu\ast\tilde\g)\le
\int_{\T^p\times\R^p}A_h(\g,(x,v))
\g\dr\mu(x)\dr v+U(\mu\ast\g)+$$
$$\left[
\left( 
-1+\log\frac{2}{\d}
\right)    +
2\sqrt p L 
+\frac{1}{2h}4p  -M(A)
\right]
\e\int_{\T^p}B_A(x)\dr\mu(x)    .  $$
Recall that $M(A)\tends+\infty$ as $A\searrow 0$; thus, if $A$ is small enough, the last formula contradicts the minimality of $\g$.

\fin

\prop{2.8} There is a constant 
$C(L,h)$, depending only on $L,h>0$, for which the following happens. Let $U$ be $L$-Lipschitz, let $\mu\in\Mt$ and let 
$\g$ minimize in the definition of $(G^h_tU)(\mu)$. Then,

\noindent 1) the function $\g$ satisfies
$$\max\left(
||\g||_{L^\infty(\T^p\times\R^p,\mu\otimes\L^p)},
||\frac{1}{\g}||_{L^\infty(\T^p\times B(0,2\sqrt p),\mu\otimes\L^p)}
\right)       \le C(L,h)     .   \eqno (2.19)$$

\noindent 2) Let us denote by $\r_{\mu\ast\g}$ the density of 
$\mu\ast\g$; then
$$||\frac{1}{\r_{\mu\ast\g}}||_{L^\infty{(\T^p,\L^p)}}    \le C(L,h)   .  \eqno (2.20)$$

\noindent 3)  The set 
$$\{ 
\r_{\mu\ast\g} \st \mu\in\Mt\txt{and}\g\in{\cal D}_\mu
\txt{is minimal}
\}$$ 
is uniformly integrable in $L^1(\T^p,\L^p)$.

\proof Point 1) is just the statement of lemmas 2.5 and 2.7. We prove point 2).

Let ${\cal F}$ denote the class of all continuous probability densities on $\T^p$. The first equality below is standard; the second one comes from the fact that $\r_{\mu\ast\g}$ is the density of $\mu\ast\g$, and the third one from the definition of 
$\mu\ast\g$; for the fourth one, we have set $Q\colon=[-\2,\2)^p$ and used the fact that $f$ is periodic. The first inequality below comes from the fact that $f$ and $\g$ are non negative; for the second one, we have pushed the measure $\mu$, which lives on $\T^p$, to a measure on $Q$, which we denote by the same letter; for the last one, we use (2.19) and the fact that, if $x,w\in Q$, then 
$x-w\in B(0,2\sqrt p)$. 
$${\rm ess}\inf \r_{\mu\ast\g}=
\inf_{f\in{\cal F}}\int_{\T^p}f(z)\r_{\mu\ast\g}(z)\dr z=
\inf_{f\in{\cal F}}\int_{\T^p}f(z)\dr(\mu\ast\g)(z)=$$
$$\inf_{f\in{\cal F}}\int_{\T^p\times\R^p}f(x-v)\g(x,v)\dr\mu(x)\dr v=
\inf_{f\in{\cal F}}\int_{\T^p}\dr \mu(x)\sum_{k\in\Z^p}\int_{Q}
f(w)\g(x,x-w-k)\dr w\ge$$
$$\inf_{f\in{\cal F}}\int_{\T^p}\dr\mu(x)\int_Qf(w)\g(x,x-w)\dr w\ge
\int_{Q}[{\rm ess}\inf_{w\in Q}\g(x,x-w)]\dr\mu(x)  \ge  
\frac{1}{C(L,h)}   .  $$
But this is (2.20). 

We prove point 3). Let $\mu_k\in\Mt$ and let 
$\g_k\in{\cal D}_{\mu_k}$ be minimal. Let $\r_k$ be the density of 
$\mu_k\ast\r_k$; we want to prove that, up to subsequences, 
$\r_k\weak\r$ in $L^1(\T^p,\L^p)$. Thus, let 
$g\in L^\infty(\T^p,\L^p)$; the first equality below is the definition of $\r_k$, the second one is the definition of $\mu_k\ast\g_k$, the last one is the change of variables in $\R^p$ $\fun{}{v}{w=x-v}$.
$$\int_{\T^p}g(z)\r_k(z)\dr z=
\int_{\T^p}
g(z)\dr(\mu_k\ast\g_k)(z)=
\int_{\T^p}\dr\mu_k(x)
\int_{\R^p}g(x-v)\g_k(x,v)\dr v  =  $$
$$\int_{\R^p}g(w)\dr w\int_{\T^p}\g_k(x,x-w)\dr\mu_k(x)  .  $$
Thus, point 3) holds if we prove that
$$a_k(w)\colon=
\int_{\T^p}\g_k(x,w-x)\dr\mu_k(x)$$
has a subsequence converging weakly in $L^1(\R^p)$. To prove this, we recall from point 1) that $\g_k(x,v)\le C(L,h)$ for 
$x\in E_k$, with $\mu_k(E_k^c)=0$; this implies that 
$||a_k||_\infty\le C(L,h)$ and thus $a_k$ is uniformly integrable on $\R^p$. 

Next, we have to show that the measures $a_k\L^p$ are tight. For the first equality below, we lift $\mu$ to a measure on 
$Q$ and use Fubini; the first inequality comes from the fact that, if $x\in Q$ and $w\ge R$, then $|w-x|\ge R-\sqrt p$; the second inequality follows by lemma 2.6 if we take $R$ large enough.
$$\int_{B(0,R)^c}a_k(w)\dr w=
\int_{Q}\dr\mu_k(x)\int_{B(0,R)^c}\g_k(x,w-x)\dr w\le$$
$$\int_{Q}\dr\mu_k(x)
\int_{B(0,R-\sqrt p)^c}\g_k(x,v)\dr v\le
\int_Q\e\dr\mu_k(x)=\e   .   $$
This proves tightness.

\fin

\vskip 2pc
\centerline{\bf \S 3}
\centerline{\bf Discrete characteristics and value functions}
\vskip 1pc

In this section, we define the characteristics and value function for the problem with discrete time, and we show that the discrete value function is bounded as the time-step tends to zero. From now on, the parameter $h$ in the definition of $G^h_tU$ (formula (1.4)) will be set to 
$h=\frac{1}{n}$, with $n\in\N$. 

\vskip 1pc

\noindent{\bf Definitions.} $\bullet$) Let $m,n\in\N$ and let $U\in C(\Mt)$ be $L$-Lipschitz for the 1-Wasserstein distance; we can define inductively the following sequence of functions.
$$\matrix{
\hat U_n(0,\mu)=U(\mu)\cr
\hat U_n(-\frac{1}{n},\mu)=
\left[
G^\frac{1}{n}_{-\frac{1}{n}}\hat U_n(0,\cdot)
\right]   (\mu)-
\log\left(
\frac{n}{2\pi}
\right)^\frac{p}{2}  \cr
\dots\cr
\hat U_n(-\frac{mn}{n},\mu)=\left[
G^\frac{1}{n}_{\frac{-mn}{n}}\hat U_n(-\frac{mn-1}{n},\cdot)
\right]   (\mu)-
\log\left(
\frac{n}{2\pi}
\right)^\frac{p}{2}     .   
}   $$
Applying iteratively proposition 2.3, we see that 
$\hat U_n(\frac{j}{n},\cdot)$ is $(L+Cm)$-Lipschitz if 
$j\in(-mn,\dots,0)$.

\noindent $\bullet$) Let $s\in(1,\dots,mn)$; we say that 
$\{ \mu^\frac{1}{n}_{\frac{j}{n}} \}_{j=-s}^{0}$ is a 
$\{ \g^\frac{1}{n}_{\frac{j}{n}} \}_{j=-s}^{-1}$-sequence starting at 
$(-\frac{s}{n},\mu)$ if the following three points hold.

\noindent 1) $\mu^\frac{1}{n}_{-\frac{s}{n}}=\mu$.

\noindent 2) $\g^\frac{1}{n}_{\frac{j}{n}}\in\dcal_{\mu^\frac{1}{n}_{\frac{j}{n}}}$ for $j=-s,\dots,-1$.

\noindent 3) $\mu^\frac{1}{n}_{\frac{j+1}{n}}=
\mu^\frac{1}{n}_{\frac{j}{n}}\ast\g^\frac{1}{n}_{\frac{j}{n}}$ for
$j=-s,\dots,-1$.

\noindent $\bullet$) If 
$\{ \mu^\frac{1}{n}_{\frac{j}{n}} \}_{j=-s}^{0}$ is a 
$\{ \g^\frac{1}{n}_{\frac{j}{n}} \}_{j=-s}^{-1}$-sequence starting at 
$(-\frac{s}{n},\mu)$, and 
$\g^\frac{1}{n}_{\frac{j}{n}}$ is minimal in the definition of
$$\left[
G^\frac{1}{n}_\frac{j}{n}\hat U(\frac{j}{n},\cdot)
\right]   \left(
\mu^\frac{1}{n}_{\frac{j}{n}}
\right)   $$
for $j=-s,\dots,-1$, then we say that 
$\{ \mu^\frac{1}{n}_{\frac{j}{n}} \}_{j=-s}^{0}$ is a minimal
$\{ \g^\frac{1}{n}_{\frac{j}{n}} \}_{j=-s}^{-1}$-sequence starting at 
$(-\frac{s}{n},\mu)$.

\noindent $\bullet$) If 
$\{ \mu^\frac{1}{n}_{\frac{j}{n}} \}_{j=-s}^{0}$ is a 
$\{ \g^\frac{1}{n}_{\frac{j}{n}} \}_{j=-s}^{-1}$-sequence and 
$t\in[\frac{-s}{n},0]$, say $t\in[\frac{j}{n},\frac{j+1}{n}]$ for some 
$j\in(-s,\dots,-1)$, we let $\mu^\frac{1}{n}_t$ be the geodesic for the 2-Wasserstein distance which connects 
$\mu^\frac{1}{n}_{\frac{j}{n}}$ at time $\frac{j}{n}$ with 
$\mu^\frac{1}{n}_{\frac{j+1}{n}}$ at time $\frac{j+1}{n}$.

\noindent $\bullet$ For $j\in(-s+1,\dots,0)$ and 
$t\in(\frac{j-1}{n},\frac{j}{n}]$ we define
$$\hat U_n(t,\mu)=\hat U_n(\frac{j}{n},\mu)  .  $$

\noindent $\bullet$) For $t\in [-m,0]$, we let
$$\hat U(t,\mu)=\liminf_{n\tends+\infty}\hat U_n(t,\mu)  .  $$

When there is no ambiguity, we shall drop the $\frac{1}{n}$, and denote $\g^\frac{1}{n}_{\frac{j}{n}}$, 
$\mu^\frac{1}{n}_{\frac{j}{n}}$ and $\mu^\frac{1}{n}_{t}$ by 
$\g_{\frac{j}{n}}$, $\mu_{\frac{j}{n}}$ and $\mu_{t}$ respectively. 

\vskip 1pc

The definitions above raise at least two questions: the first one is the convergence of a $\{ \g^\frac{1}{n}_{\frac{j}{n}} \}_j$-minimal sequence $\{ \m^\frac{1}{n}_{\frac{j}{n}} \}_j$ to a minimal characteristic as $n\tends+\infty$; this will have to wait until section 6 for an answer. The second one is whether 
$\hat U(t,\mu)$ is finite; this is the content of proposition 3.2 below. Before proving it, we need a definition and a lemma. 

\noindent $\bullet$) Let 
$\{ \mu^\frac{1}{n}_{\frac{j}{n}} \}_{j=-s}^{0}$ be a 
$\{ \g^\frac{1}{n}_{\frac{j}{n}} \}_{j=-s}^{-1}$-sequence starting at 
$(-\frac{s}{n},\mu)$; we define the functional
$$I(\mu,\g^\frac{1}{n}_\frac{-s}{n},\dots,\g^\frac{1}{n}_\frac{-1}{n})=
\sum_{j=-s}^{-1}\int_{
\T^p\times\R^p
}\left[
\frac{1}{n}L^{
\2\mu_\frac{j}{n}
}  (t,x,nv)+\log\g_\frac{j}{n}(x,v)
\right]    \g_\frac{j}{n}(v)\dr\mu_\frac{j}{n}(x,v)\dr v  .  $$

We omit the proof of the next lemma, since the fact that the value function defines the Hopf-Lax semigroup is standard.

\lem{3.1} Let $\{ \bar\mu_\frac{j}{n} \}_j$ be a minimal 
$\{ \bar\g_\frac{j}{n} \}_j$-sequence starting at $(-\frac{s}{n},\mu)$; then, $\{ \bar\g_\frac{j}{n} \}_j$ minimizes the functional which brings $(\g_{\frac{-s}{n}},\g_\frac{-s+1}{n},\dots,\g_\frac{-1}{n})$ to
$$I(\mu,\g^\frac{1}{n}_\frac{-s}{n},\dots,\g^\frac{1}{n}_\frac{-1}{n})+
U(\mu_0) - 
s\log\left(
\frac{n}{2\pi}
\right)^\frac{p}{2}    .   \eqno (3.1)     $$
Moreover, the value of the minimum is equal to 
$\hat U(-\frac{s}{n},\mu)$.

In other words, if for $s>j$ we define
$$T_{-\frac{s}{n},-\frac{j}{n}}U(\mu)\colon=
\min_{(\g_{\frac{-s}{n}},\g_\frac{-s+1}{n},\dots,\g_\frac{-j}{n})}
I(\mu,\g^\frac{1}{n}_\frac{-s}{n},\dots,\g^\frac{1}{n}_\frac{-j}{n})+
U(\mu_0) - 
(s-j)\log\left(
\frac{n}{2\pi}
\right)^\frac{p}{2}  $$
then $T_{-\frac{s}{n},-\frac{j}{n}}$ is a semigroup in the past: for $t>j>s$, 
$T_{-\frac{t}{n},-\frac{s}{n}}\circ T_{-\frac{s}{n},-\frac{j}{n}}=
T_{-\frac{t}{n},-\frac{j}{n}}$.

\fin

\prop{3.2} There is $C>0$, only depending on $m\in\N$, such that 
for all $\mu\in\Mt$ and $t\in[-m,0]$, we have 
$$|\hat U(t,\mu)|\le C  .  $$ 

\proof It suffices to show that there is $C>0$ such that
$$\left\vert
\hat U_n(\frac{-s}{n},\mu)
\right\vert   \le C\quad
\forall s\in(0,1,\dots,mn),\quad\forall n\ge 1,\quad
\forall\mu\in\Mt   .   \eqno (3.2)$$
For $j\in(-s,\dots,-1)$ let us set
$$\tilde\g_{\frac{j}{n}}(v)=\left(
\frac{n}{2\pi}
\right)^\frac{p}{2} e^{
-\frac{n}{2}|v|^2
}      $$
and let $\{ \tilde\mu_j \}$ be a $\{ \tilde\g_j \}$-sequence starting at $(-\frac{s}{n},\mu)$. By lemma 3.1 we get the first inequality below; by (1.13) and the fact that $U$ is bounded, the second one follows.
$$\hat U_n(-\frac{s}{n},\mu)\le$$
$$\sum_{j=-s}^{-1}\int_{
\T^p\times\R^p
}\left[
\frac{1}{n}L^{
\2\tilde\mu_\frac{j}{n}
}  (t,x,nv)+\log\tilde\g_\frac{j}{n}(x,v)
\right]    \tilde\g_\frac{j}{n}(x,v)\dr\tilde\mu_\frac{j}{n}(x)\dr v+
U(\tilde\mu_0) - 
s\log\left(
\frac{n}{2\pi}
\right)^\frac{p}{2}\le$$
$$\sum_{j=-s}^{-1}\int_{\T^p\times\R^p}\left[
\frac{n}{2}|v|^2+\log\tilde\g_\frac{j}{n}(x,v)
\right]  \tilde\g_\frac{j}{n}(x,v)\dr\tilde\mu_\frac{j}{n}(x)\dr v +
\frac{s}{n}(
||V||_\infty+||W||_\infty
)        + ||U||_{\sup}
-s\log\left(\frac{n}{2\pi}\right)^\frac{p}{2}  .  $$
Since
$$\int_{\T^p\times\R^p}
\left[
\cinn{v}+\log\tilde\g_\frac{j}{n}(x,v)
\right]     \tilde\g_\frac{j}{n}(x,v)\dr\tilde\mu_\frac{j}{n}(x)\dr v=
\log\left(
\frac{n}{2\pi }
\right)^\frac{p}{2}     $$
and $s\le nm$, we get that
$$\hat U_n(\frac{-s}{n},\mu)\le 
m(||V||_\infty+||W||_\infty)+||U||_{\sup}   .  \eqno (3.3)$$
To prove the opposite inequality, we let 
$\{ \mu^\frac{1}{n}_{\frac{j}{n}} \}_j$ be a minimal 
$\{ \g^\frac{1}{n}_{\frac{j}{n}} \}_j$-sequence starting at 
$(-\frac{s}{n},\mu)$, with $s\in(1,2,\dots,mn)$.  By lemma 1.1, the function $\tilde\g_\frac{j}{n}$ defined above minimizes the integral of $A_\frac{1}{n}$; as a consequence,
$$\int_{\T^p\times\R^p}\left[
\frac{n}{2}|v|^2+\log\g_\frac{j}{n}(x,v)
\right]  \g_\frac{j}{n}(x,v)\dr\mu_\frac{j}{n}(x)\dr v\ge$$
$$\int_{\T^p\times\R^p}\left[
\frac{n}{2}|v|^2+\log\tilde\g_\frac{j}{n}(x,v)
\right]  \tilde\g_\frac{j}{n}(x,v)\dr\mu_\frac{j}{n}(x)\dr v=
\log\left(
\frac{n}{2\pi }
\right)^\frac{p}{2}   .  $$
Lemma 3.1 implies the equality below; the first inequality comes from (1.13) and the fact that $U$ is bounded; the second one comes from the formula above.
$$\hat U_n\left(
\frac{-s}{n},\mu
\right) =$$   
$$\sum_{j=-s}^{-1}\int_{
\T^p\times\R^p
}\left[
\frac{1}{n}L^{
\2\mu_\frac{j}{n}
}  (t,x,nv)+\log\g_\frac{j}{n}(x,v)
\right]    \g_\frac{j}{n}(x,v)\dr\mu_\frac{j}{n}(x)\dr v+
U(\mu_0)  -
s
\log\left(
\frac{n}{2\pi}
\right)^\frac{p}{2}   \ge$$
$$\sum_{j=-s}^{-1}
\int_{\T^p\times\R^p}\left[
\frac{n}{2}|v|^2+\log\g_\frac{j}{n}(x,v)
\right]  \g_\frac{j}{n}(x,v)\dr\mu_\frac{j}{n}(x)\dr v -C_7-
s
\log\left(
\frac{n}{2\pi}
\right)^\frac{p}{2} \ge -C_7   .  $$
This inequality and (3.3) imply (3.2).

\fin

\vskip 2pc
\noindent\centerline{\bf \S 4}
\noindent\centerline{\bf Differentiability of $U$}
\vskip 1pc

In this section, we want to show that the minimal of 
$\fun{}{\psi}{S(U,\mu,\psi)}$ also minimizes a problem for a linear final condition. Proposition 4.1 below deals with a single time step, while proposition 4.6 deals with the whole history. 

\prop{4.1} Let $U$ be Lipschitz and differentiable on densities (see below for a definition). Let 
$\mu\in\Mt$, and let $(G^\frac{1}{n}_tU)(\mu)$ be achieved on 
$\g$. Let 
$\r_{\g}$ be the density of $\mu\ast\g$, and let $f$ be the differential of $U$ at $\mu\ast\g=\r_{\g}\L^p$. Then, there is a bounded Borel function $a(x)$ such that
$$\g(x,v)=e^{-\cinn{v}-f(x-v)+a(x)}   .  \eqno (4.1)$$
Moreover, there is a constant $M>0$, independent of $\mu$ and 
$h$, such that
$$||f||_\infty\le M   .   \eqno (4.2)$$

\rm

The following definition will come handy.

\vskip 1pc

\noindent{\bf Definition.} If $\fun{f}{\T^p}{\R}$ is a bounded Borel function, we define
$$\fun{U_f}{\Mt}{\R},\qquad U_f(\mu)\colon=\int_{\T^p}f\dr\mu  .  
$$

\vskip 1pc

We sketch the proof of proposition 4.1: we are going to  concentrate on the particles which lie in a small ball centered in 
$x\in\T^p$, i. e. we shall consider the probability measure 
$\mu_{in}\colon=\frac{1}{\mu(B(x,r))}\mu|_{B(x,r)}$. We shall see that the optimal strategy for $\mu_{in}$ approximates, as 
$r\tends 0$, the optimal strategy for a single particle problem; the final condition is the linear $U_f$, where $f$ is the derivative of $U$ at $\mu\ast\g$. When the final condition is linear, the minimizer can be written explicitly by [14], and (4.1) will follow. As for the potential, in the case of the single time step it won't even enter the picture; for more time steps, i. e. in lemma 4.5, we shall see that it tends to the mean field generated by all the particles. 

In general, the gradient of a function on the space of probability measures is defined as a vector field ([2] and [16]) rather delicate to find; however, in our case $\mu\ast\g$ is the convolution of a probability measure with a $L^1$ function, and thus has a $L^1$ density. As a consequence, we can use the standard definition of derivative in $L^1(\T^p)$. 

\vskip 1pc

\noindent{\bf Definition.} We say that $\fun{U}{\Mt}{\R}$ is differentiable at densities if the following two points hold.

\noindent $i$) There is a function 
$\fun{\tilde U}{L^1(\T^p)}{\R}$ such that 
$U(\phi\L^p)=\tilde U(\phi)$ when $\phi\L^p$ is a probability measure. We ask that $\tilde U$ be differentiable at every probability density $\phi$; in other words, there is a function 
$h\in L^\infty(\T^p)$ such that
$$\left\vert
\tilde U((\phi+\psi)\L^p)-\tilde U(\phi\L^p)-
\int_{\T^p}h\cdot\psi\dx
\right\vert   =   o(||\psi||_{L^1(\T^p)})    .   $$
We set $U^\prime(\phi\L^p)\colon=h$. Actually, we shall need the formula above only on the affine space of probability densities, i. e. when $h$ has zero mean. 

\noindent $ii$) We also ask that there is $M>0$ such that, if 
$h=U^\prime(\phi\L^p)$ for a probability density $\phi$, then
$$||h||_{L^\infty(\T^p)}\le M   .  $$
Note that, if $U$ is Lipschitz for $d_1$, then point $ii$) holds automatically; indeed, in this case $U$ is Lipschitz also for the total variation distance, which easily implies point $ii$). 

\vskip 1pc

The typical example of a function $U$ differentiable on densities is the usual one: we take $k$ bounded Borel functions 
$\fun{f_1,\dots,f_k}{\T^p}{\R}$ and we set
$$U(\mu)=\left(
\int_{\T^p}f_1\dr\mu
\right)\cdot
\dots
\cdot\left(
\int_{\T^p}f_k\dr\mu
\right)   .  $$

As we saw above, $\mu\ast\g$ has a density; this leads us to the first of the following definitions.

\vskip 1pc

\noindent {\bf Definitions.} $\bullet$) Let $U\in C(\Mt)$; we define
$$||U||_{den}\colon=\sup\{
|U(\r\L^p)|  \st \r\in L^1(\T^p),\quad  \r\ge 0,\quad
\int_{\T^p}\r(x)\dx=1
\}    .  $$
If $f\in L^\infty(\T^p)$, then it is easy to see that 
$$||U_f||_{den}= ||f||_{L^\infty(\T^p)}  .  $$
Since in $L^\infty(\T^p)$ we disregard null sets with respect to the Lebesgue measure,  the $\sup$ of $|U_f|$ on $\Mt$ could be larger than $||f||_{L^\infty(\T^p)}$.

\noindent$\bullet$) We want to isolate the particles in $B(x,r)$; thus, for $\mu\in\Mt$ and $x\in\T^p$, we define
$$\mu_{ext}=\mu\vert B(x,r)^c,\qquad
\mu_{in}=\frac{1}{\mu(B(x,r))}\mu\vert B(x,r)  .  $$

\noindent$\bullet$) Let $U$, $f$, $\mu$ and $\g$ be as in the hypotheses of proposition 4.1; for $\psi\in{\cal D}_\mu$, we define the function 
$U^r_\psi$ as 
$$U^r_{\psi}(\l)=\frac{1}{\mu(B(x_0,r))}\cdot
\left\{
U[
\mu_{ext}\ast\psi+\mu(B(x,r))\l
]   -U(\mu\ast\g)  +\mu(B(x,r))\cdot U_f(\mu_{in}\ast\g)
\right\}     =   $$
$$\frac{1}{\mu(B(x,r))}\cdot\big\{
U[
\mu\ast\psi+\mu(B(x,r))\cdot(\l-\mu_{in}\ast\psi)
]   -U(\mu\ast\g)   +\mu(B(x,r))\cdot U_f(\mu_{in}\ast\g)
\big\}     .  \eqno (4.3)$$
As for the second equality above, it comes from the fact, easy to  check, that the operator $\ast$ is linear in $\mu$: 
$(\mu_1+\mu_2)\ast\g=\mu_1\ast\g+\mu_2\ast\g$.

\vskip 1pc

Lemma 4.2 below shows that $U^r_\g$ is the final condition that is seen by the particles in $B(x,r)$; by lemma 4.3 below, $U^r_\g$ is very close to the derivative of $U$ at $\mu\ast\g$.

\lem{4.2} Let $U$, $\mu$ and $\g$ be as in proposition 4.1, let 
$x\in\T^p$ and let $U^r_\g$ be defined as in (4.3); then $\g|B(x,r)$ minimizes $\fun{}{\psi}{S(U^r_\g,\mu_{in},\psi)}$.

\proof We set 
$$\hat U(\nu)=U(\nu)-U(\mu\ast\g)+
\mu(B(x,r))\cdot U_f(\mu_{in}\ast\g)  $$
and we note that the minima of 
$S(U,\mu,\cdot)$ coincide with the minima of 
$S(\hat U,\mu,\cdot)$. Indeed, adding a constant to the final condition does not change the set of minima.

Next, we isolate the particles in $B(x,r)$: taking 
$\psi\in{\cal D}_\mu$, the first equality below is the definition of $S$, while the second one comes from the definition of $\mu_{in}$ and $\mu_{ext}$. 
$$S(\hat U,\mu,\psi)=
\int_{\T^p\times\R^p}\left[
\cinn{v}+\log\psi
\right]     \psi\dr\mu(x)\dr v+
\hat U(\mu\ast\psi)=$$
$$\int_{\T^p\times\R^p}\left[
\cinn{v}+\log\psi(x,v)
\right]     \psi(x,v)\dr\mu_{ext}(x)\dr v+$$
$$\mu(B(x,r))\cdot
\int_{\T^p\times\R^p}\left[
\cinn{v}+\log\psi(x,v)
\right]     \psi(x,v)\dr\mu_{in}(x)\dr v   +
\hat U(\mu_{ext}\ast\psi+\mu(B(x,r))\mu_{in}\ast\psi) .  $$
If $U^r_{\psi}$ is defined as in (4.3), then we can write the formula above as 
$$S(\hat U,\mu,\psi)=$$
$$\int_{\T^p\times\R^p}\left[
\cinn{v}+\log\psi(x,v)
\right]     \psi(x,v)\dr\mu_{ext}(x)\dr v+   \eqno (4.4)_a$$
$$\mu(B(x,r))\cdot S(U^r_{\psi},\mu_{in},\psi)   .
\eqno (4.4)_b$$
Since $\g$ minimizes in (4.4) and  $(4.4)_a$ does not depend on 
$\g |_{B(x,r)\times\R^p}$, we have that 
$\g|_{B(x,r)\times\R^p}$ must minimize $(4.4)_b$, i. e.  
$\fun{}{\psi}{S(U^r_\psi,\mu_{in},\psi)}$. This is almost the thesis, save for the fact that we have $U^r_\psi$ in stead of $U^r_\g$. But we can restrict to the functions $\psi$ such that 
$\psi |_{B(x,r)^c\times\R^p}=\g |_{B(x,r)^c\times\R^p}$; in this way we have that $\mu_{ext}\ast\psi=\mu_{ext}\ast\g$; since $\psi$ enters the definition of $U^r_\psi$ only through 
$\mu_{ext}\ast\psi$ (this is the first equality of (4.3)), we get that 
$U^r_\psi=U^r_\g$, and the thesis follows.

\fin

\lem{4.3} Let $U$, $f$, $\mu$ and $\g$ be as in the hypotheses of proposition 4.1; let us suppose that $x$ is not an atom of $\mu$. Then,
$$\lim_{r\tends 0}||U^r_\g-U_f||_{den}=0 .  \eqno (4.5)$$

\proof Let $\eta$ be a probability density on $\T^p$. The first quality below is the definition of $U^r_\g$, the second one comes from the fact that $U_f$ is the differential of $U$ at $\mu\ast\g$; in the "small oh" we have denoted by $\r_{in}$ the density of 
$\mu_{in}\ast\g$. 
$$|
U^r_\g(\eta\L^p)-U_f(\eta\L^p)
|   =   $$
$$\Bigg\vert
\frac{1}{\mu(B(x,r))}  \cdot
\{
U[
\mu\ast\g+\mu(B(x,r))\cdot(\eta\L^p-\mu_{in}\ast\g)
]-U(\mu\ast\g)+$$
$$\mu(B(x,r))\cdot U_f(\mu_{in}\ast\g)
\}  -U_f(\eta\L^p)  
\Bigg\vert    =$$
$$  \Bigg\vert
\frac{1}{\mu(B(x,r))}
\{
U_f[
\mu(B(x,r))\cdot(\eta\L^p-\mu_{in}\ast\g)
]+\mu(B(x,r))\cdot U_f(\mu_{in}\ast\g)+$$
$$o[\mu(B(x,r))\cdot ||\eta-\r_{in}||_{L^1(\T^p)}]
\}
-U_f(\eta\L^p)
\Bigg\vert    .  $$
Now we note that, for all probability density $\eta$, 
$$\mu(B(x,r))\cdot ||\eta-\r_{in}||_{L^1(\T^p)}\le 2\mu(B(x,r))    $$
and thus
$$|
U^r_\g(\eta\L^p)-U_f(\eta\L^p)
|   =
\frac{1}{\mu(B(x,r))} o(\mu(B(x,r)))     $$
where the "small oh" does not depend on $\eta$. Since $x$ is not an atom, $\mu(B(x,r))\tends 0$ and the thesis follows.

\fin

\lem{4.4} Let $U$ and $\mu$ be as in proposition 4.1 and let $\g$ minimize $S(U,\mu,\cdot)$. Then, for $\mu$ a. e. $x\in\T^p$ which is not an atom,
$$\liminf_{r\tends 0}S(U^r_\g,\mu_{in},\g)\ge
S(U_f,\d_x,\g)      \eqno (4.6)$$
and
$$\limsup_{r\tends 0} S(U^r_\g,\mu_{in},\g)\le
\min_\psi S(U_f,\d_x,\psi)   .   \eqno (4.7)$$

\proof  {\bf Step 1.} We begin to show that, for $\mu$ a. e. $x$, atom or not, the three fact below hold. First, that 
$$\lim_{r\tends 0}
\frac{1}{\mu(B(x,r))}
\int_{B(x,r)\times\R^p}
|
\g(x,v)-\g(y,v)
|    \dr\mu(y)\dr v=0   .  \eqno (4.8)$$
Second, if $f=U^\prime(\mu\ast\g)$, then
$$\lim_{r\tends 0}
\frac{1}{\mu(B(x,r))}
\int_{B(x,r)\times\R^p}  |
f(y-v)\g(y,v)-f(x-v)\g(x,v)
|   \dr\mu(y)\dr v  =0  .   \eqno (4.9)$$
Third, if $\hat\g(v)=\left(\frac{n}{2\pi p}\right)^\frac{p}{2}
e^{-\cinn{v}}$, then 
$$\lim_{r\tends 0}\frac{1}{{\mu(B(x,r))}}
\int_{B(x,r)\times\R^p}|f(y-v)-f(x-v)|\hat\g(v)\dr\mu_{in}(y)\dr v 
=0  .  \eqno (4.10)$$
We begin with the standard proof of (4.8): we let 
$\{ \g_m \}_{m\ge 1}$ be a dense sequence in $L^1(\R^p)$ and consider the Borel measures on $\T^p$
$$\mu_m(A)\colon=\int_{A\times\R^p}
|\g(y,v)-\g_m(v)|\dr\mu(y)\dr v  .  $$
By the Lebesgue differentiation theorem, for all $x\in E$ with 
$\mu(E^c)=0$ we have that, for all $m$,
$$\lim_{r\tends 0} \frac{1}{\mu(B(x,r))}
\int_{B(x,r)\times\R^p}|\g(y,v)-\g_m(v)|\dr\mu(y)\dr v=
\int_{\R^p}|\g(x,v)-\g_m(v)|\dr v   .  \eqno (4.11)$$
For $x\in E$ and $\e>0$, we choose $\g_m$ such that
$$\int_{\R^p}|\g(x,v)-\g_m(v)|\dr v\le\e   .   \eqno (4.12)$$
The first inequality below is obvious, while the equality follows by (4.11) and the last inequality by (4.12).
$$\limsup_{r\tends 0}\frac{1}{\mu(B(x,r))}
\int_{\mu(B(x,r))\times\R^p}
|\g(y,v)-\g(x,v)|\dr\mu(y)\dr v   \le  $$
$$\lim_{r\tends 0}\frac{1}{\mu(B(x,r))}
\int_{\mu(B(x,r))\times\R^p}
|\g(y,v)-\g_m(v)|
\dr\mu(y)\dr v    +$$
$$\lim_{r\tends 0}\frac{1}{\mu(B(x,r))}
\int_{\mu(B(x,r))\times\R^p}
|\g_m(v)-\g(x,v)|
\dr\mu(y)\dr v = $$
$$2\int_{\R^p}|\g_m(v)-\g(x,v)|\dr v\le 2\e  .  $$
Since $\e$ is arbitrary, (4.8) follows. Formulas (4.9) and (4.10) follow by the same argument, but applied to $f(y-v)\g(y,v)$ and 
$f(y-v)\hat\g(v)$ respectively.  

For the next steps, we suppose that $x\in E$ and $x$ is not an atom of $\mu$.

\noindent{\bf Step 2.} Here we prove (4.6). For $\e>0$ and 
$x\in E$ let us set 
$$F\colon=\{
y\in B(x,r)\st ||\g(y,\cdot)-\g(x,\cdot)||_{L^1(\R^p)}<\e
\}   .  $$
By (4.8) and the Chebishev inequality we have that
$$\mu_{in}(F)\tends 1 \txt{and}
\mu_{in}(B(x,r)\setminus F)\tends 0  
\txt{as}r\tends 0  .   \eqno (4.13)$$
Now,
$$\int_{B(x,r)}
\dr\mu_{in}(y)\int_{\R^p}A_\frac{1}{n}(\g,(y,v))\dr v=$$
$$\int_{F}
\dr\mu_{in}(y)\int_{\R^p}A_\frac{1}{n}(\g,(y,v))\dr v+
\int_{B(x,r)\setminus F}
\dr\mu_{in}(y)\int_{\R^p}A_\frac{1}{n}(\g,(y,v))\dr v  . 
\eqno (4.14) $$
We saw in lemma 1.3 that the map
$$\fun{}{\psi}{\int_{\R^p}A_\frac{1}{n}(\psi,v)}\dr v$$
is l. s. c. for the $L^1$ topology; thus, there is $\d(\e)\tends 0$ as 
$\e\tends 0$ such that 
$$\int_{\R^p}A_\frac{1}{n}(\g(y,\cdot),v)\dr v\ge
\int_{\R^p}A_\frac{1}{n}(\g(x,\cdot),v)\dr v-\d(\e)\quad
\forall y\in F. $$
This implies the first equality below, while the second one comes from (4.13).
$$\liminf_{r\tends 0}\int_F\dr\mu_{in}(y)
\int_{\R^p}A_\frac{1}{n}(\g,(y,v))\dr v\ge
\liminf_{r\tends 0}\mu_{in}(F)\left[
\int_{\R^p}A_\frac{1}{n}(\g(x,\cdot),v)\dr v 
-\d(\e) 
\right]     =  $$
$$\int_{\R^p}A_\frac{1}{n}(\g(x,\cdot),v)\dr v -\d(\e) .  $$
By (4.13) and (1.6) we get that
$$\liminf_{r\tends 0}\int_{B(x,r)\setminus F}\dr\mu_{in}(y)
\int_{\R^p}A_\frac{1}{n}(\g,(y,v))\dr v\ge 0  .  $$
Thus, by the last two formulas and (4.14),
$$\liminf_{r\tends 0} \int_{B(x,r)}\dr\mu_{in}(y)
\int_{\R^p}A_\frac{1}{n}(\g,(x,v))\dr v\ge
\int_{\R^p}A_\frac{1}{n}(\g(x,\cdot),v)\dr v  -\d(\e)  .  $$
Since $\e$ is arbitrary and $\d(\e)$ tends to zero as $\e\tends 0$, we have that 
$$\liminf_{r\tends 0} \int_{B(x,r)}\dr\mu_{in}(y)
\int_{\R^p}A_\frac{1}{n}(\g,(x,v))\dr v\ge
\int_{\R^p}A_\frac{1}{n}(\g(x,\cdot),v)\dr v    .  $$
The equality below comes from the definition of $U_f$ and of 
$\mu\ast\g$; the limit comes from (4.9). 
$$|
U_f(\mu_{in}\ast\g)-U_f(\d_x\ast\g)
|  =$$
$$\left\vert
\int_{\T^p\times\R^p}f(y-v)\g(y,v)\dr\mu_{in}(y)\dr v-
\int_{\R^p}f(x-v)\g(x,v)\dr v
\right\vert \tends 0   .  $$
The first inequality below comes from the definition of $S$; the first equality follows from (4.5); the last inequality comes from the last two formulas above.
$$\liminf_{r\tends 0} S(U^r_\g,\mu_{in},\g)\ge
\liminf_{r\tends 0}[
S(U_f,\mu_{in},\g)-  ||U^r_\g-U_f||_{den}
]   =$$
$$\liminf_{r\tends 0} S(U_f,\mu_{in},\g)=
\liminf_{r\tends 0}\left[
\int_{B(x,r)}\dr\mu_{in}(y)
\int_{\R^p} A_\frac{1}{n}(\g,(y,v))\dr v+
U_f(\mu_{in}\ast\g)
\right]   \ge$$
$$\int_{\R^p}A_\frac{1}{n}(\g(x,\cdot),v)\dr v+
U_f(\d_x\ast\g)=S(U_f,\d_x,\g)   .  $$
This proves (4.6).

\noindent{\bf Step 3.} We prove (4.7). 
We know from [14] that the function 
$\fun{}{\psi}{S(U_f,\d_x,\psi)}$ has a unique minimum, given by 
$$\tilde\g(v)=e^{
-\cinn{v}-f(x-v)+a(x)
}   $$
with $a(x)$ such that $\tilde\g$ is a probability density. In order to compare the actions, we are going to plug $\tilde\g$ into 
$S(U_\g^{r},\mu_{in},\cdot)$. The first inequality below holds because we have seen in lemma 4.2 that $\g$ minimizes 
$\fun{}{\psi}{S(U^r_\g,\mu_{in},\psi)}$; the second one holds by the definition of $S$; the first equality is the definition of $S$ while the last one comes from the fact that 
$\tilde\g$ does not depend on $y$ and $\mu_{in}$ is a probability measure. 
$$S(U^r_\g,\mu_{in},\g)\le
S(U^{r}_\g,\mu_{in},\tilde\g)\le
S(U_f,\mu_{in},\tilde\g)+||U_\g^{r}-U_f||_{den}=  $$
$$\int_{B(x,r)}\dr\mu_{in}(y)
\int_{\R^p}A_\frac{1}{n}(\tilde\g,v)\dr v+U_f(\mu_{in}\ast\tilde\g)+
||U_\g^{r}-U_f||_{den}=$$
$$\int_{\R^p}A_\frac{1}{n}(\tilde\g,v)\dr v+
U_f(\mu_{in}\ast\tilde\g)+||U_\g^{r}-U_f||_{den}  .  $$
By (4.5) we have that 
$$||U_\g^{r}-U_f||_{den}\tends 0
\txt{as} r\tends 0   .  $$
The first inequality below comes from the definitions of 
$\mu_{in}\ast\g$ and $\d_x\ast\g$; the second one, from the definition of $\tilde\g$, 
$\hat\g$ (the definitions are just above and in step 1 respectively) and the fact that $||f||_\infty\le M$; the limit comes from (4.10).
$$|U_f(\mu_{in}\ast\tilde\g)-U_f(\d_x\ast\tilde\g)|\le
\int_{\T^p\times\R^p}|f(y-v)-f(x-v)|\tilde\g(v)
\dr\mu_{in}(y)\dr v\le$$
$$e^{M}\int_{\R^p}|f(y-v)-f(x-v)|\hat\g(v)\dr\mu_{in}(y)\dr v
\tends 0  . $$
From the last three formulas, we get that
$$\limsup_{r\tends 0}S(U^{r}_\g,\mu_{in},\g)\le
S(U_f,\d_x,\tilde\g)      $$
as we wanted.

\fin 

\noindent{\bf End of the proof of proposition 4.1.} By (4.6) and (4.7) we see that, for $\mu$ a. e. $x\in\T^p$ which is not an atom, 
$$S(U_f,\d_x,\g(x,\cdot)) \le
\min_\psi S(U_f,\d_x,\psi)   .  $$
In other words, $\g(x,\cdot)$ coincides with a minimum of 
$S(U_f,\d_x,\cdot)$; by [14], this functional has just one minimum, which is the right hand side of formula (4.1). This proves (4.1), while (4.2) follows from point $ii$) of the definition of differentiability on densities. 

It remains to prove (4.1) when $x$ is an atom of $\mu$. To show this, we have to enlarge our set of controls. Namely, let us suppose for simplicity that $\mu$ has just one atom, say $x_0$ with $\mu(\{ x_0 \})=\l$; let us write 
$$\mu=\tilde\mu+\l\d_{x_0}   .   $$
Then we assign to each $x\not=x_0$ a strategy $\g(x,\cdot)$ as before, but we assign to $x_0$ an enlarged set of controls, say 
$\g_w(x_0,\cdot)$ with $w\in [0,\l]$: in a sense, we are supposing that in $x_0$ sits a continuum of particles, each parametrized by $w$ and each with a strategy $\g_w(x_0,\cdot)$. We define
$$K(U,\mu,\g)=
\int_{\T^p\times\R^p}A_\frac{1}{n}(\g,(x,v))\dr\tilde\mu(x)\dr v+
\int_0^\l\dr w\int_{\R^p}A_\frac{1}{n}(\g_w(x_0,\cdot),v)\dr v+$$
$$U\left(
\tilde\mu\ast\g+\int_0^\l (\d_{x_0}\ast\g_w(x_0,\cdot))\dr w
\right)    .  $$

Two things are clear:

\noindent 1) first, that the minimum of $K(U,\mu,\cdot)$ is lower than the minimum of $S$, simply because we have a larger set of controls.

\noindent 2) Second, that if we find a minimum $\g$ of $K$ such that $\g_w(x_0,\cdot)$ does not depend on $w$, then it is also a minimum of $S$: indeed, in this case $K(U,\mu,\g)=S(U,\mu,\g)$ and point 1) implies the assertion. 

Thus, (4.1) follows if we show that any minimum 
$(\g_w(x_0,\cdot),\g(x,\cdot))$ of 
$K(U,\mu,\cdot)$ is given by (4.1) for $\L^1$ a. e. 
$w\in[0,\l]$ and $\tilde\mu$ a. e. $x\in\T^p$. This is done exactly as in lemma 4.4: indeed, instead of the torus we are considering 
$(\T^p\setminus \{ x_0 \})\sqcup [0,\l]$ with the measure 
$\tilde\m$ on $\T^p\setminus \{ x_0 \}$ and $\L^1$ on $[0,l]$; to this space and measure the proof of lemma 4.4 applies. We avoid repeating the details: if $w_0\in(0,l)$, as in lemma 4.4 we isolate that particles $w$ with $|w-w_0|<r$ and, letting $r\tends 0$, we show that 
$$S(U_f,\d_{x_0},\g_{w_0}(x_0,\cdot))\le
\min_\psi S(U_f,\d_{x_0},\psi)  .  $$
Since the unique minimal of the expression on the right is the linear one given by [14], i. e. formula (4.1), we are done.

\fin

Proposition 4.1 gives an explicit expression for the minimal 
$\g_\frac{-1}{n}$; now we want to extend this result to more than one time-step, i. e. the situation of section 3. However, if we want to find an explicit expression for the minimizer 
$(\g_\frac{-s}{n},\dots,\g_\frac{-1}{n})$ of 
$$I(\mu,\g_\frac{-s}{n},\dots,\g_\frac{-1}{n})+U(\mu_0),$$
we need a slightly different proof. The reason for this is that, even if we isolate the particles in $B(z_0,r)$ at the initial time 
$\frac{-s}{n}$, they are going to spread over all $\T^p$ at time 
$\frac{-s+1}{n}$, and after this time their trajectories coincide with the rest of the pack. In other words, after the first step, there is no way to control some particles separately from the other ones. To tackle this problem, we are going to minimize over a larger set of controls which keeps track of the initial position. 

\vskip 1pc

\noindent{\bf Definition.} We consider the functions 
$\psi^\frac{1}{n}_{\frac{-s}{n},z},\psi^\frac{1}{n}_{\frac{-s+1}{n},z},\dots, \psi^\frac{1}{n}_{\frac{-1}{n},z}$ depending measurably on 
$z\in\T^p$; we let
$$\mu^\psi_{\frac{-s}{n},z}=\d_z,\quad
\mu^\psi_{\frac{-s+1}{n},z}=\mu^\psi_{\frac{-s}{n},z}\ast
\psi^\frac{1}{n}_{\frac{-s}{n},z},\quad\dots,\quad
\mu^\psi_{0,z}=\mu^\psi_{\frac{-1}{n},z}\ast
\psi^\frac{1}{n}_{\frac{-1}{n},z}   $$
be a $(\psi^\frac{1}{n}_{\frac{-s}{n},z},\psi^\frac{1}{n}_{\frac{-s+1}{n},z},\dots, \psi^\frac{1}{n}_{\frac{-1}{n},z})$-sequence starting at 
$(-\frac{s}{n},\d_z)$. In other words, $\mu^\psi_{\frac{j}{n},z}$ is the distribution at time $\frac{j}{n}$ of the particle which started at $z$. 

If $\mu_\frac{-s}{n}$ is the initial distribution of the particles at time 
$\frac{-s}{n}$, we define the total distribution of all the particles at time   $\frac{j}{n}\ge\frac{-s}{n}$ by
$$\mu^\psi_\frac{j}{n}=
\int_{\T^p}\mu_{\frac{j}{n},z}\dr\mu^\psi_\frac{-s}{n}(z)  .  $$
We say that the mean field generated by all the particles at time 
$\frac{j}{n}$ is
$$W^{\mu^\psi_\frac{j}{n}}(x)  .  $$
We define the cost for particle $z$ analogously as the functional 
$I$ of lemma 3.1; it considers the history of a particle subject to the mean field generated by the whole community. 
$$I_\2(\d_z,\psi^\frac{1}{n}_{\frac{-s}{n},z},\dots,
\psi^\frac{1}{n}_{\frac{-1}{n},z})=$$
$$\sum_{j=-s}^{-1}\int_{\T^p\times\R^p}
[
\frac{1}{n} L^{\2\mu_\frac{j}{n}^\psi}(\frac{j}{n},x,nv)
+\log\psi_{\frac{j}{n},z}(x,v)
]\psi_{\frac{j}{n},z}(x,v)
\dr\mu^\psi_{\frac{j}{n},z}(x)\dr v=$$
$$\sum_{j=-s}^{-1}\int_{\T^p\times\R^p}\{
A_\frac{1}{n}(\psi_\frac{j}{n}^\frac{1}{n},(x,v))-
\frac{1}{n}[
V(\frac{j}{n},x)+ W^{\2\mu^\psi_\frac{j}{n}}(x)]
\}
\dr\mu_{\frac{j}{n},z}^\psi(x)\dr v  .  $$
We have called it $I_\2$ because of the coefficient $\2$ in 
$W^{\2\mu^\psi_\frac{j}{n}}$; we shall call $I_1$ its counterpart with $W^{\mu^\psi_\frac{j}{n}}$.

Integrating over the initial distribution $\mu_\frac{-s}{n}$ and adding the final condition, we define the cost for all particles.
$$J_\2[U,\quad\mu_\frac{-s}{n},\quad(\psi^\frac{1}{n}_{\frac{-s}{n},z},\psi^\frac{1}{n}_{\frac{-s+1}{n},z},\dots, 
\psi^\frac{1}{n}_{\frac{-1}{n},z})]  
\colon=$$
$$\int_{\T^p}
I_\2(\d_z,\psi^\frac{1}{n}_{\frac{-s}{n},z},\psi^\frac{1}{n}_{\frac{-s+1}{n},z},\dots \psi^\frac{1}{n}_{\frac{-1}{n},z})
\dr\mu_\frac{-s}{n}(z)   +
U\left(
\mu_0^\psi
\right)   .  $$

We omit the proof that $J_\2(U,\mu_\frac{-s}{n},\cdot)$ has a minimum, since it is identical to proposition 1.4. 

\lem{4.5} Let 
$(\g^\frac{1}{n}_{\frac{-s}{n},z},\dots,\g^\frac{1}{n}_{\frac{-1}{n},z})$ minimize the functional
$$\colon (\psi^\frac{1}{n}_{\frac{-s}{n},z},\psi^\frac{1}{n}_{\frac{-s+1}{n},z},\dots \psi^\frac{1}{n}_{\frac{-1}{n},z})\tends
J_\2[U,\quad\mu_\frac{-s}{n},\quad(\psi^\frac{1}{n}_{\frac{-s}{n},z},\psi^\frac{1}{n}_{\frac{-s+1}{n},z},\dots \psi^\frac{1}{n}_{\frac{-1}{n},z})]  .  $$
Then, for $\mu_\frac{-s}{n}$ a. e. $z\in\T^p$ which is not an atom, 
$(\g^\frac{1}{n}_{\frac{-s}{n},z},\dots,\g^\frac{1}{n}_{0,z})$ does not depend on $z$ and has the following expression. Let $f_0$ be the derivative of $U$ at  the measure $\mu_0^\g$ defined above. For 
$j\in(-s+1,\dots,0)$ we define by backward induction
$$f_\frac{j-1}{n}(x)=-\log\int_{\R^p}
e^{-\cinn{v}+
\frac{1}{n}V(\frac{j}{n},x)+\frac{1}{n}W^{\mu_\frac{j}{n}}(x)-
f_\frac{j}{n}(x-v)
}   \dr v   .   $$
Then, we have that
$$\g_\frac{j}{n}(x,v)=
e^{
-\cinn{v}-f_\frac{j}{n}(x-v)+a_\frac{j}{n}(x)
}   \eqno (4.15) $$
with $a_\frac{j}{n}(x)$ chosen in such a way that 
$\g_\frac{j}{n}(x,\cdot)$ is a probability density for all $x\in\T^p$. 

\proof In the first two steps below, which correspond to lemma 4.2, we isolate the particles in $B(z_0,r)$; in the third one we let 
$r\tends 0$ as in lemma 4.3. It is in this step that we need that 
$\mu_\frac{-s}{n}(\{ z_0 \})=0$.

\noindent {\bf Step 1.} In this step, we set some notation and add a constant to $U$, as we did at the beginning of lemma 4.1. 

For $j\ge -s$ we set
$$\mu^\psi_{\frac{j}{n},int}=
\frac{1}{\mu_\frac{-s}{n}(B(z_0,r))}
\int_{B(z_0,r)}\mu^\psi_{\frac{j}{n},z}\dr\mu_\frac{-s}{n}(z)  .  $$
This is the distribution at time $\frac{j}{n}$ of the particles which started in $B(z_0,r)$ at time $\frac{-s}{n}$. 

Let $(\g^\frac{1}{n}_{\frac{-s}{n},z},
\g^\frac{1}{n}_{\frac{-s+1}{n},z},\dots, \g^\frac{1}{n}_{\frac{-1}{n},z})$ be as in the hypotheses; as at the beginnig of this section, we define 
$$\tilde U(\mu)\colon=U(\mu)-
U\left(
\mu_0^\g
\right)    +
\mu_\frac{-s}{n}(B(z_0,r))U_f(\mu^\g_{0,in}) $$
and
$$U^r_\psi(\l)=
\frac{1}{\mu(B(z_0,r))}
\tilde U\left[
\int_{\T^p}\mu_{0,z}^\psi\dr\mu_\frac{-s}{n}(z)+
\mu_\frac{-s}{n}(B(z_0,r))\cdot
\left(
\l-\mu_{0,int}^\psi
\right)  
\right]   .  $$
As in lemma 4.2, we shall see that $U^r_\g$ is the final condition seen by the particles in $B(z_0,r)$. Note that, as in lemma 4.2, 
$$\tilde U(\mu^\psi_0)=\mu(B(z_0,r))U^r_\psi(\mu^\psi_{0,int})  .   \eqno (4.16)$$

Since the addition of a constant to $U$ does not change the set of minima neither of $J_\2(\mu_\frac{-s}{n},\cdot)$ nor of 
$I_\2(\mu_\frac{-s}{n},\cdot)$, we have that 
$(\g^\frac{1}{n}_{\frac{-s}{n},z},\dots,\g^\frac{1}{n}_{\frac{-1}{n},z})$ is a minimum of the functional
$$\fun{}{
(\psi^\frac{1}{n}_{\frac{-s}{n},z},\dots,\psi^\frac{1}{n}_{\frac{-1}{n},z})
}{
J_\2[\tilde U,\quad\mu_\frac{-s}{n},\quad 
(\psi^\frac{1}{n}_{\frac{-s}{n},z},\dots,\psi^\frac{1}{n}_{\frac{-1}{n},z})]
}    .   $$

\noindent{\bf Step 2.} In this step, we deal with the mutual interaction; this is the main difference with lemma 4.2, where there was none of it. We define $W_{\mu_\frac{j}{n},in}$ and
$W_{\mu_\frac{j}{n},ext}$ as the potentials generated by the particles starting in $B(z_0,r)$ and $B(z_0,r)^c$ respectively, i. e. 
$$W_{\mu_\frac{j}{n},in}(x)\colon=
\int_{B(z_0,r)}\dr\mu_\frac{-s}{n}(z)
\int_{\T^p}W(x-y)\dr\mu^\psi_{\frac{j}{n},z}(y)\dr y  ,  \eqno (4.17)_a$$
$$W_{\mu_\frac{j}{n},ext}(x)\colon=
\int_{B(z_0,r)^c}\dr\mu_\frac{-s}{n}(z)
\int_{\T^p}W(x-y)\dr\mu^\psi_{\frac{j}{n},z}(y)\dr y    .   \eqno (4.17)_b$$
Now our particles interact among themselves only through the potential $W$. Note that $W$ appears in 
$I_\2(\d_z,\psi^\frac{1}{n}_{\frac{-s}{n},z},\dots,
\psi^\frac{1}{n}_{\frac{-1}{n},z})$ in terms of the form
$$\int_{\T^p}W^{\2\mu^\psi_\frac{j}{n}}(x)\dr\mu_{\frac{j}{n},z}(x) . 
$$
Integrating in $\mu_{\frac{-s}{n}}$, we get that $W$ appears in
$J_\2[\tilde U^r_\psi,\quad\mu_\frac{-s}{n},\quad 
(\psi^\frac{1}{n}_{\frac{-s}{n},z},\dots,\psi^\frac{1}{n}_{\frac{-1}{n},z})]$ in terms which have the form of the left hand side of the equation below; the first equality is the definition of 
$W^{\2\mu^\psi_\frac{j}{n}}$, while the second one is the definition of $\mu^\psi_\frac{j}{n}$. 
$$\int_{\T^p}\dr\mu_\frac{-s}{n}(z)
\int_{\T^p}W^{\2\mu^\psi_\frac{j}{n}}(x)\dr\mu^\psi_{\frac{j}{n},z}(x)=
\2\int_{\T^p}\dr\mu_\frac{-s}{n}(z)
\int_{\T^p}\dr\mu^\psi_{\frac{j}{n},z}(x)
\int_{\T^p}W(x-y)\dr\mu^\psi_{\frac{j}{n}}(y)=$$
$$\2\int_{\T^p}\dr\mu_\frac{-s}{n}(z)
\int_{\T^p}\dr\mu_\frac{-s}{n}(w)
\int_{\T^p\times\T^p}W(x-y)\dr\mu_{\frac{j}{n},z}^\psi(x)
\dr\mu_{\frac{j}{n},w}^\psi(y)  .  $$
The term on the right in the formula above is the sum of the three terms below: the first one is the interaction of $B(z_0,r)^c$ with itself, the last one is the interaction of $B(z_0,r)$ with itself, while the middle one is the interaction of $B(z_0,r)$ with $B(z_0,r)^c$; note that here a factor $\2$ is missing due to the symmetry of the potential. 
$$\int_{\T^p}\dr\mu_\frac{-s}{n}(z)
\int_{\T^p}
W^{\2\mu^\psi_\frac{j}{n}}(x)\dr\mu_{\frac{j}{n},z}^\psi(x)=$$
$$\2\int_{B(z_0,r)^c}\dr\mu_\frac{-s}{n}(z)
\int_{B(z_0,r)^c}\dr\mu_\frac{-s}{n}(w)
\int_{\T^p\times\T^p}W(x-y)\dr\mu_{\frac{j}{n},z}^\psi(x)
\dr\mu_{\frac{j}{n},w}^\psi(y)  +  $$
$$\int_{B(z_0,r)}\dr\mu_\frac{-s}{n}(z)
\int_{B(z_0,r)^c}\dr\mu_\frac{-s}{n}(w)
\int_{\T^p\times\T^p}W(x-y)\dr\mu_{\frac{j}{n},z}^\psi(x)
\dr\mu_{\frac{j}{n},w}^\psi(y)    +$$
$$\2\int_{B(z_0,r)}\dr\mu_\frac{-s}{n}(z)
\int_{B(z_0,r)}\dr\mu_\frac{-s}{n}(w)
\int_{\T^p\times\T^p}W(x-y)\dr\mu_{\frac{j}{n},z}^\psi(x)
\dr\mu_{\frac{j}{n},w}^\psi(y)  .  $$
Using this and $(4.17)_{a-b}$ above, we can write 
$$\int_{\T^p}\dr\mu_\frac{-s}{n}(z)
\int_{\T^p}W^{\2\mu^\psi_\frac{j}{n}}(x)
\dr\mu_{\frac{j}{n},z}^\psi(x)=
\2\int_{B(z_0,r)^c}\dr\mu_\frac{-s}{n}(z)
\int_{\T^p}
W_{\mu_{\frac{j}{n},ext}}(x)\dr\mu^\psi_{\frac{j}{n},z}(x)+$$
$$\int_{B(z_0,r)}\dr\mu_\frac{-s}{n}(z)
\int_{\T^p}
W_{\mu_{\frac{j}{n},ext}}(x)\dr\mu^\psi_{\frac{-s}{n},z}(x)+
\2\int_{B(z_0,r)}\dr\mu_\frac{-s}{n}(z)
\int_{\T^p}
W_{\mu_{\frac{j}{n},int}}(x)\dr\mu^\psi_{\frac{j}{n},z}(x)  .  $$
The first term above crops up in $(4.18)_a$ below, the second one crops up in $(4.18)_b$ and the third one in $(4.18)_c$; we have used (4.16) to get $(4.18)_b$.  
$$J_\2[\tilde U,\quad\mu_\frac{-s}{n},\quad(\psi^\frac{1}{n}_{\frac{-s}{n},z},\psi^\frac{1}{n}_{\frac{-s+1}{n},z},\dots, \g^\frac{1}{n}_{\frac{-1}{n},z})]=  $$
$$\int_{B(z_0,r)^c}\dr\mu_\frac{-s}{n}(z)
\sum_{j=-s}^{-1}
\int_{\T^p\times\R^p}\left[
A_\frac{1}{n}(\g_{\frac{j}{n},z},(x,v))-
\frac{1}{n}V(\frac{j}{n},x)-\frac{1}{2n} W_{\mu_\frac{j}{n},ext}(x)
\right]
\dr\mu_{\frac{j}{n},z}(x)\dr v  +\eqno (4.18)_a$$
$$\mu_\frac{-s}{n}(B(z_0,r))\Bigg[
\int_{B(z_0,r)}
\dr\mu_{\frac{-s}{n},in}(z)
\sum_{j=-s}^{-1}
\int_{\T^p\times\R^p}
[A_\frac{1}{n}(\g_{\frac{j}{n},z},(x,v))-\frac{1}{n}V(\frac{j}{n},x)-
\frac{1}{n}W_{\mu_\frac{j}{n},ext}(x)]\dr\mu_{\frac{j}{n},z}(x)\dr v+$$
$$\tilde U^r_\psi\left(
\mu_0^\psi
\right)
\Bigg]    -   \eqno (4.18)_b$$
$$\int_{B(z_0,r)}\dr\mu_\frac{-s}{n}(z)
\sum_{j=-s}^{-1}\int_{\T^p}
\frac{1}{2n} W_{\mu_\frac{j}{n},int}(x)\dr\mu_{\frac{j}{n},z}(x)   .   
\eqno (4.18)_c$$
We shall call $\hat J$ the term in the square parentheses in  $(4.18)_b$; it is almost equal to the functional $J_1$, the only difference being that the potential is 
$V(\frac{j}{n},x)+W_{\mu^\psi_{\frac{j}{n},ext}}$ instead of 
$V(\frac{j}{n},x)+W^{\mu^\psi_{\frac{j}{n}}}$. Note that we have lost  the constant $\2$ before the potential $W$. 

Now $(4.18)_a$ is not affected by 
$(\g^\frac{1}{n}_{\frac{-s}{n},z},\dots,
\g^\frac{1}{n}_{\frac{-1}{n},z})$ when $z\in B(z_0,r)$; this prompts us to restrict, as in lemma 4.2, to functions 
$(\psi^\frac{1}{n}_{\frac{-s}{n},z},\dots, \psi^\frac{1}{n}_{\frac{-1}{n},z})$ which coincide with 
$(\g^\frac{1}{n}_{\frac{-s}{n},z},\dots, \g^\frac{1}{n}_{\frac{-1}{n},z})$ for $z\not\in B(z_0,r)$; since 
$(\g^\frac{1}{n}_{\frac{-s}{n},z},\dots, \g^\frac{1}{n}_{\frac{-1}{n},z})$ is minimal, we see that 
$(\g^\frac{1}{n}_{\frac{-s}{n},z},\dots, \g^\frac{1}{n}_{\frac{-1}{n},z})|_{z\in B(z_0,r)}$ must minimize $(4.18)_{b-c}$. Note that, by our choice of 
$(\psi^\frac{1}{n}_{\frac{-s}{n},z},\dots, \psi^\frac{1}{n}_{\frac{-1}{n},z})$, $\tilde U^r_\psi=\tilde U^r_\g$; in other words, 
$(\g^\frac{1}{n}_{\frac{-s}{n},z},\dots, \g^\frac{1}{n}_{\frac{-1}{n},z})|_{z\in B(z_0,r)}$ minimizes
$$\mu_\frac{-s}{n}(B(z_0,r))\cdot
\hat J[\tilde U^r_\g,\quad\mu_{in},\quad(\g^\frac{1}{n}_{\frac{-s}{n},z},\dots,
\g^\frac{1}{n}_{\frac{-1}{n},z})]    +   \eqno (4.19)_a$$
$$\int_{B(z_0,r)}\dr\mu_\frac{-s}{n}(z)
\sum_{j=-s}^{-1}\int_{\T^p}
\2 W_{\mu_\frac{j}{n},int}(x)\dr\mu_{\frac{j}{n},z}(x)   .   
\eqno (4.19)_b$$

\noindent {\bf Step 3.} We want to use the fact that 
$(\g^\frac{1}{n}_{\frac{-s}{n},z},\dots, \g^\frac{1}{n}_{\frac{-1}{n},z})|_{z\in B(z_0,r)}$ minimizes $(4.19)_{a-b}$ to get (4.15). First of all, we fix $z_0$, a Lebesgue point of 
$$\fun{}{z}{(\g^\frac{1}{n}_{\frac{-s}{n},z},\g^\frac{1}{n}_{\frac{-s+1}{n},z},\dots, \g^\frac{1}{n}_{\frac{-1}{n},z})}  \eqno (4.20)$$
for the measure $\mu_\frac{-s}{n}$. 

Since we are supposing that $\{ z_0 \}$ is not an atom of 
$\mu_\frac{-s}{n}$, by $(4.17)_a$ we have that 
$$(4.19)_b=o(\mu(B(z_0,r)))  .  \eqno (4.21)$$
Since this term is negligible with respect to $(4.19)_a$, with a proof similar to that of lemma 4.3 we get that 
$$\limsup_{r\tends 0}
J_1[\tilde U^r_\g,\quad\mu_{in},\quad(\g^\frac{1}{n}_{\frac{-s}{n},z},\dots,
\g^\frac{1}{n}_{\frac{-1}{n},z})]\le
\min_{\psi_{\frac{-s}{n},z_0},\dots,\psi_{\frac{-1}{n},z_0}}
J_1[\tilde U_{f_0},\quad\d_{z_0},\quad(\psi_{\frac{-s}{n},z_0},\dots,
\psi_{\frac{-1}{n},z_0})]   .  $$
Note that here we are dealing with $J_1$: the coefficient $\2$ in 
$W^{\2\mu^\psi_\frac{j}{n}}$ was shed already in (4.18). 

Moreover, we can see as in formula (4.6) of lemma 4.3 that 
$$\liminf_{r\tends 0}
J_1[\tilde U^r_\g,\quad\mu_{in},\quad(\g^\frac{1}{n}_{\frac{-s}{n},z},\dots,
\g^\frac{1}{n}_{\frac{-1}{n},z})]\ge
J_1[\tilde U_{f_0},\quad\d_{z_0},\quad(\g_{\frac{-s}{n},z_0},\dots,
\g_{\frac{-1}{n},z_0})]   .  $$
As in the proof of lemma 4.2, the last two formulas imply that 
$(\g_{\frac{-s}{n},z_0},\dots,\g_{\frac{-1}{n},z_0})$ minimizes the term on the right in the formula above; now [14] prescribes that 
$(\g_{\frac{-s}{n},z_0},\dots,\g_{\frac{-1}{n},z_0})$ satisfies (4.15).

\fin

\prop{4.6} Let $U$ be $L$-Lipschitz and differentiable on densities, let $s\in(1,\dots,mn)$ and  let $\{ \mu_\frac{j}{n} \}_j$ be a minimal 
$\{ \g_\frac{j}{n} \}_j$-sequence starting at $\mu_{-\frac{s}{n}}$. Let $f_0$ be the derivative of $U$ at $\mu_0$. For 
$j\in(-s+1,\dots,0)$ we define by backward induction
$$f_\frac{j-1}{n}(x)=-\log\int_{\R^p}
e^{-\cinn{v}+
\frac{1}{n}V(\frac{j}{n},x)+\frac{1}{n}W^{\mu_\frac{j}{n}}(x)-
f_\frac{j}{n}(x+v)
}   \dr v   .   $$
Then, we have that
$$\g_\frac{j}{n}(x,v)=
e^{
-\cinn{v}-f_\frac{j}{n}(x+v)+a_\frac{j}{n}(x)
}    \eqno (4.22)$$
with $a_\frac{j}{n}(x)$ chosen in such a way that 
$\g_\frac{j}{n}(x,\cdot)$ is a probability density for all $x\in\T^p$. 

\proof We shall prove the assertion when $\mu_\frac{-s}{n}$ has no atoms; the argument for the atoms of $\mu_\frac{-s}{n}$ is identical to the one in the proof of proposition 4.1, and we skip it.

For the functional $I$ we defined before lemma 3.1, let us set 
$$\hat I_\2[U,\quad\mu_\frac{-s}{n},\quad(\psi_\frac{-s}{n},\dots,\psi_\frac{-1}{n})]=
I(\mu_\frac{-s}{n},\psi_\frac{-s}{n},\dots,\psi_\frac{-1}{n})
+U(\mu_0)  .  $$
By lemma 3.1 it suffices to show that, if 
$(\g_\frac{-s}{n},\dots,\g_\frac{-1}{n})$ minimizes $\hat I_\2$, then it satisfies (4.22). We prove this. If we compare $\hat I_\2$ with the function $J_\2$ of the last lemma, we see two things:

\noindent 1) the minimum of $J_\2(U,\mu_\frac{-s}{n},\cdot)$ is smaller than the minimum of $\hat I_\2(U,\mu_\frac{-s}{n},\cdot)$, simply because the dependence on $z\in\T^p$ gives us a larger set of strategies.

\noindent 2) If $(\g^\frac{1}{n}_{\frac{-s}{n},z},\g^\frac{1}{n}_{\frac{-s+1}{n},z},\dots, \g^\frac{1}{n}_{\frac{-1}{n},z})$ minimizes 
$J_\2(U,\mu_\frac{-s}{n},\cdot)$ and does not depend on $z$, then 
$$\hat I_\2[U,\quad\mu_\frac{-s}{n},\quad(\g^\frac{1}{n}_{\frac{-s}{n},z},\g^\frac{1}{n}_{\frac{-s+1}{n},z},\dots, \g^\frac{1}{n}_{\frac{-1}{n},z})]=
J_\2[U,\quad\mu_\frac{-s}{n},\quad(\g^\frac{1}{n}_{\frac{-s}{n},z},\g^\frac{1}{n}_{\frac{-s+1}{n},z},\dots, \g^\frac{1}{n}_{\frac{-1}{n},z})]  .  $$

A consequence is the following: suppose we can find a minimizer of $J_\2(U,\mu_\frac{-s}{n},\cdot)$ which does not depend on $z$, then it is also a minimal of $\hat I_\2(U,\mu_\frac{-s}{n},\cdot)$; thus, the value of the minimum for the two functional is the same and any minimizer of $\hat I_\2(U,\mu_\frac{-s}{n},\cdot)$ is a minimizer of $J_\2(U,\mu_\frac{-s}{n},\cdot)$ too. In other words, the proposition follows if we prove that any minimizer $(\g^\frac{1}{n}_{\frac{-s}{n},z},\g^\frac{1}{n}_{\frac{-s+1}{n},z},\dots, \g^\frac{1}{n}_{\frac{-1}{n},z})$ of $J_\2(U,\mu_\frac{-s}{n},\cdot)$ has the form (4.15); but that is the content of lemma 4.5.

\fin

\vskip 2pc
\centerline{\bf \S 5}
\centerline{\bf Regularity of the linearized action}
\vskip 1pc

Thanks to proposition 4.6, we can express the minimals 
$\g^\frac{1}{n}_\frac{j}{n}$ in terms of the functions 
$f^\frac{1}{n}_\frac{j}{n}$; in this section, we shall suitably normalize these functions and show, in proposition 5.2 below, that they are regular; by Ascoli-Arzel\`a\ this will imply (lemma 5.3 below) that, up to subsequences, they converge to a function $u$. We shall use proposition 5.2 in the next section, when we prove that $u$ solves Hamilton-Jacobi and that the discretized characteristics converge.  

\noindent{\bf Definitions.} $\bullet$ Let $Q$ be a symmetric, positive-definite matrix and let $\a\in\R^n$; we denote by $N(\a,Q)$ the Gaussian of mean $\a$ and variance $Q$, i. e.
$$N(\a,Q)(v)=
\frac{1}{\sqrt{(2\pi)^p{\rm det}Q}}
e^{
-\2\inn{Q^{-1}(v-\a)}{v-\a}
}   .  $$

\vskip 1pc

\noindent$\bullet$ Let $m\in\N$, $s\in(0,\dots,mn)$; 
let $\{ \mu^{\frac{1}{n}}_\frac{j}{n} \}_j$ be  a minimal 
$\{ \g^{\frac{1}{n}}_\frac{j}{n} \}_j$-sequence starting at 
$\mu_{\frac{-s}{n}}$ and let $f_0$ be the derivative of 
$U$ at $\mu_0^\frac{1}{n}$. As in proposition 4.6, we define by backward induction 
$$e^{  
-f_\frac{j-1}{n}(x)
}         
\colon =   
\int_{\R^p}e^{
-\cinn{v}-f_\frac{j}{n}(x-v)+
\frac{1}{n}V(\frac{j}{n},x)+\frac{1}{n}W^{\mu_\frac{j}{n}}(x)
}       \dr v =$$
$$e^{
\frac{1}{n}P_\frac{j}{n}(x)
}
\int_{\R^p}e^{-\cinn{v}}
e^{
-f_\frac{j}{n}(x-v)
}   \dr v   , 
\qquad j\in(-s+1,\dots,0)   \eqno (5.1)$$
where we have set
$$P_\frac{j}{n}(x)=V(\frac{j}{n},x)+W^{\mu_\frac{j}{n}}(x)  .
\eqno (5.2)$$
Once more by proposition 4.6, we have for the minimal 
$\g^\frac{1}{n}_\frac{j}{n}$ the expression 
$$\g^{\frac{1}{n}}_\frac{j}{n}(x,v)=
e^{
-\cinn{v}-f_\frac{j}{n}(x-v)+a_\frac{j}{n}(x)
}     \eqno (5.3)$$
where $a_\frac{j}{n}(x)$ is such that 
$\g^{\frac{1}{n}}_\frac{j}{n}(x,\cdot)$ is a probability density for all $x$. 

\noindent $\bullet$ We normalize the functions $f_\frac{j}{n}$, setting
$$\bar f_\frac{j}{n}(x)\colon=f_\frac{j}{n}(x)-
|j|\log\left(\frac{n}{2\pi}\right)^\frac{p}{2}    $$
and we see that (5.1) becomes 
$$e^{  
-\bar f_\frac{j-1}{n}(x)
}         
=   
e^{
\frac{1}{n}P_\frac{j}{n}(x)
}
\int_{\R^p}N(0,\frac{1}{n}Id)(v)e^{
-\bar f_\frac{j}{n}(x-v)
}   \dr v     \eqno (5.4)  $$
or, equivalently,
$$\bar f_\frac{j-1}{n}(x)=
-\frac{1}{n}P_\frac{j}{n}(x)-
\log\left[
\int_{\R^p}
N(0,\frac{1}{n}Id)(v)e^{-\bar f_\frac{j}{n}(x-v)}\dr v     
\right]    . \eqno (5.5)$$

\noindent $\bullet$ We set 
$$b_\frac{j}{n}(x)=-|j-1|\log\left(
\frac{n}{2\pi}
\right)^\frac{p}{2}        +a_\frac{j}{n}(x)$$
and (5.3) becomes
$$\g^\frac{1}{n}_\frac{j}{n}(x,v)=\left(
\frac{n}{2\pi}
\right)^\frac{p}{2}
e^{
-\cinn{v}-\bar f_\frac{j}{n}(x-v)-
|j-1|\log\left(\frac{n}{2\pi}\right)^\frac{p}{2}+
a_\frac{j}{n}(x)
}    =N(0,\frac{1}{n}Id)(v)e^{
-\bar f_\frac{j}{n}(x-v)+b_\frac{j}{n}(x)
}     .  \eqno (5.6)$$
In the following, we shall drop the bar from $\bar f_\frac{j}{n}$ and call it simply $f_\frac{j}{n}$. 

\noindent $\bullet$ We shall say that 
$\{ f_\frac{j}{n}^\frac{1}{n} \}_{j=-s}^0$ is the linearized cost for the minimal characteristic starting at 
$\left( \frac{-s}{n},\mu \right)$. 

\noindent $\bullet$ We gather here two other bits of notation: if $P_\frac{j}{n}$ is as in (5.2), we set 
$${\cal P}(x_\frac{j+1}{n},x_\frac{j+2}{n},\dots,x_\frac{-1}{n})=
\exp\left\{
\frac{1}{n}[P_\frac{j+1}{n}(x_\frac{j}{n})+
P_\frac{j+2}{n}(x_\frac{j+1}{n})+\dots+
P_0(x_\frac{-1}{n}) ]  
\right\}  .   \eqno (5.7)$$

\noindent $\bullet$ We also give a name to the linear path which at time $t=\frac{j}{n}<0$ is in $0$ and at time $t=0$ is in $y$:
$$a_{y}(t)=
\frac{n}{|j|}\left( t-\frac{j}{n} \right) y   .  \eqno (5.8)$$

\vskip 1pc

In the next lemma we shall see that (5.4) is simply a version of the Feynman-Kac formula. This is by no means surprising: indeed, the Hopf-Cole transform $\fun{}{f}{e^{-f}}$ brings Hamilton-Jacobi into the Schr\"odinger equation, for which Feynman-Kac provides a solution. We refer the reader to [14] for a discussion of this.

\lem{5.1} Let $U$ be Lipschitz and differentiable on densities. 
Let $\{ \mu^\frac{1}{n}_\frac{j}{n} \}_j$ be a minimal 
$\{ \g^\frac{1}{n}_\frac{j}{n} \}_j$-sequence starting at 
$\left( \frac{-s}{n},\mu \right)$; we saw above that 
$\g^\frac{1}{n}_\frac{j}{n}$ has the form (5.6) for a function 
$f_\frac{j}{n}$ defined by (5.4). 

Let $E_{0,0}$ denote the expectaction of the Brownian bridge 
$\tilde w$ which is in $0$ at $t=\frac{j}{n}$ and at $t=0$ (see [11] for a definition). Let ${\cal P}$ and $a_y$ be as in (5.7) and (5.8) respectively. 

Then,
$$e^{-f_\frac{j}{n}(x)}=$$
$$\int_{\R^p}N\left(
0,\frac{|j|}{n}Id
\right)   (x-z)  e^{-f_0(z)}
E_{0,0}\left[
{\cal P}(
x-a_{x-z}(\frac{j}{n})-\tilde w(\frac{j}{n}) , \dots,
x-a_{x-z}(-\frac{1}{n})-\tilde w(-\frac{1}{n})
)
\right]  
\dr z   .  \eqno (5.9)$$

\proof Let $j\in(s,s+1,\dots,-1)$; for 
$v_\frac{j+1}{n},\dots,v_0\in\R^p$, we set
$$\tilde v_\frac{j+1}{n}=v_\frac{j+1}{n}
\txt{and, if}\frac{l}{n}>\frac{j+1}{n},\quad
\tilde v_\frac{l}{n}=v_\frac{j+1}{n}+
\dots+v_\frac{l}{n}  .  $$
Heuristically, our particle will be in $x$ at time $\frac{j}{n}$, in 
$x+\tilde v_\frac{j+1}{n}$ at time $\frac{j+1}{n}$, ending up in 
$x+\tilde v_0$ at time $0$; the increment at each step is 
$v_\frac{j}{n}$.

Given $f_0$, which is the derivative of $U$ at $\mu_0$, we can use (5.4) to get $f_\frac{-1}{n}$ and then, iterating backwards, 
$f_\frac{-2}{n}$, $f_\frac{-3}{n}$, etc...; in this way, we get the first equality below, while the second one comes from the fact that the map 
$\fun{}{(v_\frac{j+1}{n},\dots,v_0)}{(\tilde v_\frac{j+1}{n},\dots,\tilde v_0)}$ has determinant one. 
$$e^{-f_\frac{j}{n}(x)}=$$
$$e^{
\frac{1}{n}P_\frac{j+1}{n}(x)
}
\int_{\R^p}N(0,\frac{1}{n}Id)(v_\frac{j+1}{n})
e^{\frac{1}{n}
P_\frac{j+2}{n}(x-\tilde v_\frac{j+1}{n})}\dr v_\frac{j+1}{n}
\int_{\R^p}N(0,\frac{1}{n}Id)(v_\frac{j+2}{n})
e^{
\frac{1}{n}P_\frac{j+3}{n}(x-\tilde v_\frac{j+2}{n})
}   \dr v_\frac{j+2}{n}
\dots$$
$$\dots\int_{\R^p}N(0,\frac{1}{n}Id)(v_\frac{-1}{n})
e^{\frac{1}{n}P_0(x-\tilde v_\frac{-1}{n})}
\dr v_\frac{-1}{n}
\int_{\R^p}N(0,\frac{1}{n}Id)(v_0)
e^{-f_0(x-\tilde v_0)}
\dr v_0    =$$
$$e^{
\frac{1}{n}P_\frac{j+1}{n}(x)
}
\int_{\R^p}N(0,\frac{1}{n}Id)(\tilde v_\frac{j+1}{n})
e^{\frac{1}{n}
P_\frac{j+2}{n}(x-\tilde v_\frac{j+1}{n})}\dr\tilde v_\frac{j+1}{n}
\cdot$$
$$\int_{\R^p}N(0,\frac{1}{n}Id)(\tilde v_\frac{j+2}{n}-
\tilde v_\frac{j+1}{n})
e^{
\frac{1}{n}P_\frac{j+3}{n}(x-\tilde v_\frac{j+2}{n})
}   \dr\tilde v_\frac{j+2}{n}
\dots$$
$$\dots\int_{\R^p}N(0,\frac{1}{n}Id)(\tilde v_\frac{-1}{n}-
\tilde v_\frac{-2}{n})
e^{\frac{1}{n}P_0(x-\tilde v_\frac{-1}{n})}
\dr\tilde v_\frac{-1}{n}
\int_{\R^p}N(0,\frac{1}{n}Id)(\tilde v_0-\tilde v_\frac{-1}{n})
e^{-f_0(x-\tilde v_0)}
\dr\tilde v_0    \eqno (5.10)$$
This equality looks complicated only because we have written in full the Wiener measure on cylinders; indeed, let $w$ be the Brownian motion with $w(\frac{j}{n})=0$ and let us denote by $E_w$ the expectation with respect to the Wiener measure; by the definition of the latter, (5.10) becomes
$$e^{-f_\frac{j}{n}(x)}=
E_w\left[
{\cal P}\left(
x-w\left(\frac{j}{n}\right) ,x-w\left(\frac{j+1}{n}\right),\dots, 
x-w(-\frac{1}{n})
\right)       
e^{-f_0(x-w(0))}
\right]       $$
where ${\cal P}$ has been defined in (5.7). 

We denote by $E_{0,y}$ the expectation of the Brownian bridge which is in $0$ at $t=\frac{j}{n}$ and in $y$ at $t=0$. By the properties of the Brownian bridge (see for instance [11]), the formula above becomes
$$e^{-f_\frac{j}{n}(x)}=
\int_{\R^p}N\left( 0,\frac{|j|}{n}Id \right) (y)
e^{-f_0(x-y)}
E_{0,y}\left[
{\cal P}\left(
x-w\left(\frac{j}{n}\right) ,x-w\left(\frac{j+1}{n}\right),\dots, 
x-w\left( -\frac{1}{n}\right)
\right)       
\right]   \dr y  .  \eqno (5.11)      $$
If $a_y$ is as in (5.8) and $\tilde w$ is a Brownian bridge which is in $0$ at $t=\frac{j}{n}$ and at $t=0$, then we have that
$$w(t)=a_{y}(t)+\tilde w(t)$$
is a Brownian bridge which is in $0$ at $t=\frac{j}{n}$ and in  $y$ at $t=0$. Thus, (5.11) becomes
$$e^{-f_\frac{j}{n}(x)}=$$
$$\int_{\R^p}N(0,\frac{|j|}{n}Id)(y)
e^{-f_0(x-y)}
E_{0,0}\left[
{\cal P}\left(
x-a_{y}\left(\frac{j}{n}\right)-\tilde w\left(\frac{j}{n}\right) , \dots,
x-a_{y}\left(-\frac{1}{n}\right)-\tilde w\left(-\frac{1}{n}\right)
\right)
\right]  
\dr y   .  $$
By the change of variables $z=x-y$ we get (5.9).

\fin

We fix $m\in\N$, which basically will be the time of formula (1) of theorem 1; in the following proofs, $D_i$ will always denote an increasing function from $[-m,0)$ to $(0,+\infty)$, independent of $n$ and of the starting point $\left( \frac{-s}{n},\mu \right)$ of the minimal characteristic, provided that $\frac{-s}{n}\in[-m,0)$. 

\prop{5.2} There is an increasing function 
$\fun{D_1}{[-m,0)}{(0,+\infty)}$ such that the following happens. If $(\frac{-s}{n},\mu)\in[-m,0)\times\Mt$, if 
$\{ \mu^\frac{1}{n}_\frac{j}{n} \}_j$ is a minimal 
$\{ \g^\frac{1}{n}_\frac{j}{n} \}_j$-sequence starting at 
$(\frac{-s}{n},\mu)$, if $f_\frac{j}{n}$ is as in (5.4), then we have that
$$||f_\frac{j}{n}||_{C^4(\T^p)}\le D_1(\frac{j}{n})   
\txt{for} -s\le j\le -1    .  \eqno (5.12) $$

\proof  If we set 
$$g_\frac{j}{n}(x,z,\tilde w) =$$
$$\exp\left\{
\frac{1}{n}\left[
P_\frac{j+1}{n}\left( x-a_{x-z}\left(\frac{j}{n}\right)-
\tilde w\left(\frac{j}{n}\right) \right)
+\dots+
P_0\left( x-a_{x-z}\left(\frac{-1}{n}\right)-
\tilde w\left(\frac{-1}{n}\right) \right)
\right]
\right\}$$
then (5.9) becomes 
$$e^{-f_\frac{j}{n}(x)}=
\int_{\R^p}N\left( 0,\frac{|j|}{n}Id \right)(x-z)
e^{-f_0(z)}E_{0,0}(g(x,z,\tilde w))\dr z  .  $$
If we differentiate under the the integral sign in the last formula, we see that 
$\partial^l_x e^{-f_\frac{j}{n}(x)}$ is the sum of terms of the form
$$a_k(x)\colon=\int_{\R^p}
\partial_x^{l-k}N\left( 0,\frac{|j|}{n}Id \right)
(x-z)e^{-f_0(z)}
\cdot  
E_{0,0}\left[
\partial_x^kg_\frac{j}{n}(x,z,\tilde w)
\right]    \dr z   \eqno (5.13)$$
with $0\le k\le l$. We are going to estimate each of the terms in the integral above. 

From (5.8) and (1.13) we see that, for $j\le l\le 0$, 
$$\left\vert
\partial^r_xP_\frac{l}{n}\left( x-a_{x-z}\left(\frac{l}{n}\right)-
\tilde w\left(\frac{l}{n}\right) \right)
\right\vert   \le C ,
\qquad 0\le r\le 4  $$
for a constant $C$ independent of $\tilde w$, $x$ and $z$. If we sum up in the definition of $g_\frac{j}{n}$, we get that there is 
$D_2>0$ such that 
$$|\partial_x^rg_\frac{j}{n}(x,z,\tilde w)|\le 
D_2
\qquad\forall x\in\T^p, \qquad 0\le r\le 4 .  \eqno (5.14)$$

On the other side, a simple calculation on the Gaussian shows that there is an increasing function 
$D_3$ on $[-m,0)$ such that
$$\int_{\R^p}
\left\vert
\partial^r_x N(0,\frac{|j|}{n}Id)(x-z)
\right\vert  \dr z \le D_3\left(
\frac{j}{n}
\right)     \txt{for} 0\le r\le 4.     \eqno (5.15)$$
By point $ii$) of the definition of differentiability on densities, we have that 
$$||e^{-f_0}||_{L^\infty(\T^p)}\le e^M  .   \eqno (5.16)$$
The first inequality below follows from (5.13) and H\"older, the second one from (5.14), (5.15) and (5.16).
$$||a_k||_\infty\le
||\partial^{l-k}_x N\left( 0,\frac{|j|}{n}Id \right) (\cdot)||_{L^1(\R^p)}
\cdot ||e^{-f_0}||_{L^\infty(\R^p)}\cdot
||\partial^k_x g_\frac{j}{n}(\cdot,z,\tilde w)||_{L^\infty(\R^p)}\le
D_4\left( \frac{j}{n} \right)   .  $$
Summing over $k\in(0,\dots,4)$, the thesis follows. 

\fin

\noindent{\bf Definition.} Let $f^\frac{1}{n}_\frac{j}{n}$ be the linearized cost for a discrete minimal characteristic starting at 
$\left( \frac{-s}{n},\mu \right)$. For $t\ge\frac{-s}{n}$, we define the function $f^\frac{1}{n}(t,x)$ by 
$$f^\frac{1}{n}(t,x)=f^\frac{1}{n}_\frac{j}{n}(x)
\txt{if} t\in[\frac{j}{n},\frac{j+1}{n})   .  $$

\lem{5.3} Let $\mu\in\Mt$, let $T\in[-m,0]$, let $s$ be the largest integer such that $T\le\frac{-s}{n}$ and let 
$\{ f^\frac{1}{n}_\frac{j}{n} \}_{j\ge s}$ be the linearized cost for a discrete minimal characteristic starting at 
$\left( \frac{-s}{n},\mu \right)$. Then there is 
$u\in Lip_{loc}([-m,0),C^2(\T^p))$ such that, up to subsequences, for all $\e\in(0,\frac{s}{n})$ we have that 
$$\sup_{(t,x)\in[\frac{-s}{n},-\e]\times\T^p}
|\partial^l_xf^\frac{1}{n}(t,x)-\partial^l_xu(t,x)| \tends 0
\txt{as}n\tends+\infty\txt{for}l=0,1,2.    $$

\proof Let us consider the maps
$$\fun{
F^\frac{1}{n}
}{
\left( \frac{-s}{n},\frac{-s+1}{n},\dots,\frac{-[n\e]}{n} \right)
}{
C^2(\T^p)
}  ,  \qquad
\fun{F^\frac{1}{n}}{\frac{j}{n}}{f^\frac{1}{n}(\frac{j}{n},\cdot)}  $$
where $[\cdot]$ denotes the integer part. By Ascoli-Arzel\`a\ the lemma follows if we prove that

\noindent 1) $F^\frac{1}{n}$ arrives in the same compact subset of $C^2(\T^p)$ for all $n$ and

\noindent 2) the functions $F^\frac{1}{n}$ are equilipschitz. 

Point 1) follows by (5.12); as for point 2), we shall show that there is an increasing function $\fun{D_2}{[-m,0)}{(0,+\infty)}$ such that
$$||
\partial^l_xf^\frac{1}{n}_\frac{j}{n}-
\partial^l_xf^\frac{1}{n}_\frac{j+1}{n}
||_{C^0(\T^p)}   \le\frac{1}{n} D_2\left(\frac{j+1}{n}\right)
\txt{for} l=0,1,2   .  $$
We begin to show the estimate above when $l=0$. 

By (5.5) we have that 
$$f_\frac{j-1}{n}(x)-f_\frac{j}{n}(x)=$$
$$\frac{-1}{n}P_\frac{j}{n}(x)-
\log\left[
\int_{\R^p}N\left( 0,\frac{1}{n}Id \right) (v)
e^{-f_\frac{j}{n}(x-v)}
\right]     -f_\frac{j}{n}(x)  .  $$
Thus, by (1.13), it suffices to show that
$$\left\vert
\log\left[
\int_{\R^p}N\left( 0,\frac{1}{n}Id \right) (v)
e^{-f_\frac{j}{n}(x-v)}
\right]   +f_\frac{j}{n}(x)
\right\vert   \le\frac{1}{n}D_3\left( \frac{j}{n} \right)\quad
\forall x\in\T^p  .  $$
We can take exponentials and get that the formula above is equivalent to 
$$\exp\left\{
\frac{-1}{n}D_3\left( \frac{j}{n} \right)
\right\}   \le 
\int_{\R^p}N\left( 0,\frac{1}{n}Id \right) (v)
e^{-f_\frac{j}{n}(x-v)+f_\frac{j}{n}(x)}      \dr v
\le\exp\left\{
\frac{1}{n}D_3\left( \frac{j}{n} \right)
\right\}          \quad
\forall x\in\T^p  .  \eqno (5.17)$$
We shall prove the estimate from above, since the one from below is analogous. Let us consider  
$[-f_\frac{j}{n}(x-v)+f_\frac{j}{n}(x)]$ when 
$v\in B(0,n^\frac{-1}{3})$; by (5.12) we can develop this function in Taylor series and get that  
$$\int_{\R^p}N\left( 0,\frac{1}{n}Id \right)
e^{-f_\frac{j}{n}(x-v)+f_\frac{j}{n}(x)}   \dr v   \le$$
$$\int_{B(0,n^\frac{-1}{3})}
N\left( 0,\frac{1}{n}Id \right) (v)
e^{f^\prime_\frac{j}{n}(x)\cdot v-\2 f^\pprime_\frac{j}{n}(x)(v,v)+
r(x,v)}\dr v+$$
$$\int_{B(0,n^\frac{-1}{3})^c}
N\left( 0,\frac{1}{n}Id \right) (v)
e^{-f_\frac{j}{n}(x-v)+f_\frac{j}{n}(x)}
\dr v    \eqno (5.18)$$
where
$$|r(x,v)|\le D_5\left( \frac{j}{n} \right)  |v|^3  .  $$
As for the first exponential in (5.18), we develop it in Taylor series, getting
$$\int_{B(0,n^\frac{-1}{3})}
N\left( 0,\frac{1}{n}Id \right) (v)
e^{f^\prime_\frac{j}{n}(x)\cdot v-\2 f^\pprime_\frac{j}{n}(x)(v,v)+r(x,v)}\dr v=$$
$$\int_{B(0,n^\frac{-1}{3})}
N\left( 0,\frac{1}{n}Id \right) (v)\left[
1+f^\prime_\frac{j}{n}(x)\cdot v-\2 f^\pprime_\frac{j}{n}(x)(v,v)+
Bil^1_\frac{j}{n}(x)(v,v)+r^\prime(x,v)
\right]   \dr v$$
where $Bil^1_\frac{j}{n}$ is a positive bilinear form bounded by $D_7\left( \frac{j}{n} \right) Id$ and
$$|r^\prime(x,v)|\le D_5^\prime\left( \frac{j}{n} \right)  |v|^3  .  $$
By the last two formulas and by standard properties of the Gaussian we get that
$$\int_{B(0,n^\frac{-1}{3})}
N\left( 0,\frac{1}{n}Id \right) (v)
e^{f^\prime_\frac{j}{n}(x)\cdot v-\2 f^\pprime_\frac{j}{n}(x)(v,v)+r(x,v)}\dr v\le
e^{
\frac{1}{n}D_6\left( \frac{j}{n} \right)
}    .   $$
On the other side, (5.12) and standard properties of the Gaussian imply that
$$\int_{B(0,n^\frac{-1}{3})}
N\left( 0,\frac{1}{n}Id \right)  (v)
e^{
-f_\frac{j}{n}(x-v)+f_\frac{j}{n}(x)
}    \dr v\le
e^{
D_7\left( \frac{j}{n} \right)  n^\frac{1}{6}
}   .  $$
If we apply the last two formulas to (5.18), (5.17) follows. 

Note that we need to know that $f_\frac{j}{n}$ is $C^2$ to get an estimate on 
$||f_\frac{j-1}{n}-f_\frac{j}{n}||_{C^0(\T^p)}$; the method for the estimate on $f^\prime_\frac{j}{n}$ and $f^\pprime_\frac{j}{n}$ is exactly analogous; the reason for the $C^4$ estimate on 
$f_\frac{j}{n}$ is that, to estimate the norm of the second derivative of $f_\frac{j-1}{n}-f_\frac{j}{n}$, we need two derivatives more on $f_\frac{j}{n}$.

\fin

\vskip 2pc
\centerline{\bf \S 6}
\centerline{\bf Fokker-Planck and Hamilton-Jacobi}
\vskip 1pc

By the last section, $f^\frac{1}{n}_\frac{j}{n}$ is a regular function; we shall use this information in proposition 6.2 below to prove that the discrete minimal characteristic 
$\{ \mu^\frac{1}{n}_\frac{j}{n} \}$ converges, as $n\tends+\infty$, to a weak solution of Fokker-Planck. Moreover we shall prove, in proposition 6.3 below, that the function $u$ we defined in lemma 5.3 is a solution of the Hamilton-Jacobi equation. 

We begin to show that $\g_\frac{j}{n}^\frac{1}{n}(x,\cdot)$ is a good approximation of a Gaussian.

\lem{6.1} Let $\{ f_\frac{j}{n} \}$ be as in (5.4); we set
$$Q_\frac{j}{n}(x)=
\left[
Id+\frac{1}{n}f^\pprime_\frac{j}{n}(x)
\right]^{-1}  \txt{and}
\b_\frac{j}{n}(x)=
\frac{1}{n}Q_\frac{j}{n}(x)f^\prime_\frac{j}{n}(x)    \eqno (6.1)$$
Then, there are increasing functions 
$\fun{D_{3},D_{4},D_{5},D_6}{[-m,0)}{[0,+\infty)}$ such that the following holds. 

For $a>0$, let $L(a)=[-\2 a,\2 a)^p$ and let 
$\g^\frac{1}{n}_\frac{j}{n}$ be as in the last lemma; we have that 
$$e^{-D_{3}\left(\frac{j}{n}\right)\frac{1}{n}}\cdot
N(\b_\frac{j}{n}(x),\frac{1}{n}Q_\frac{j}{n}(x))(v)
e^{-d_\frac{j}{n}(x,v)}\le  $$
$$\g^\frac{1}{n}_\frac{j}{n}(x,v)    \le
e^{D_{3}\left(\frac{j}{n}\right)\frac{1}{n}}\cdot
N(\b_\frac{j}{n}(x),\frac{1}{n}Q_\frac{j}{n}(x))(v)
e^{-d_\frac{j}{n}(x,v)}     \eqno (6.2)$$
where $d_\frac{j}{n}$ is a function such that 
$$|d_\frac{j}{n}(x,v)|\le
\frac{D_4(\frac{j}{n})}{n} \txt{if}v\in L(n^\frac{-1}{3}) .
          \eqno (6.3)   $$
Moreover,
$$\int_{\R^p\setminus L(n^\frac{-1}{3})}\g^\frac{1}{n}_\frac{j}{n}(x,v)    \le
e^{D_{5}\left(\frac{s_2}{n}\right)}e^{-n^\frac{1}{6}}     
\eqno (6.4)$$
and
$$\sup_{\R^p\setminus L(n^\frac{1}{3})}\g^\frac{-1}{n}_\frac{j}{n}(x,v)    \le
e^{D_{6}\left(\frac{s_2}{n}\right)}e^{-n^\frac{1}{6}}     .
\eqno (6.5)$$

\proof {\bf Step 1.} We are going to use Taylor's formula to get an equivalent expression for $\g^\frac{1}{n}_\frac{j}{n}$. We begin to  note that, by (5.12), there is an increasing function 
$\fun{D_7}{[-m,0)}{(0,+\infty)}$, not depending on $n$, such that 
$$|| \b_\frac{j}{n} ||_\infty  +
||Q_\frac{j}{n}(x)-Id||_\infty\le 
D_7\left(\frac{j}{n}\right)\frac{1}{n} . \eqno (6.6)$$
The first equality below is (5.6), the second one is the definition of the function $d_\frac{j}{n}(x,v)$.
$$\g^\frac{1}{n}_\frac{j}{n}(x,v)=
e^{b_\frac{j}{n}(x)}
\left(\frac{n}{2\pi}\right)^\frac{p}{2}
e^{
-\cinn{v}-f_\frac{j}{n}(x-v)
}=    \eqno (6.7)_a$$
$$\left(\frac{n}{2\pi}\right)^\frac{p}{2}
\exp\left\{
b_\frac{j}{n}(x)-f_\frac{j}{n}(x)
+\2\inn{nQ_\frac{j}{n}(x)^{-1}\b_\frac{j}{n}(x)}{\b_\frac{j}{n}(x)}   
\right\}  \cdot \eqno (6.7)_b$$
$$\exp\left[
-\2\inn{nQ_\frac{j}{n}^{-1}(x)(v-\b_\frac{j}{n}(x))}{v-\b_\frac{j}{n}(x)} -
d_\frac{j}{n}(x,v)
\right]   .   \eqno (6.7)_c$$
Since by Taylor's formula,
$$f_\frac{j}{n}(x-v)=f_\frac{j}{n}(x)-
\inn{f^\prime_\frac{j}{n}(x)}{v}+
\2\inn{f^\pprime_\frac{j}{n}(x)v}{v}+\tilde d_\frac{j}{n}(x,v)$$
and easy but lengthy computation implies that 
$\tilde d_\frac{j}{n}(x,v)=d_\frac{j}{n}(x,v)$; together with (5.12) this implies that there is a function $D_8$, bounded on $[-m,-\e]$ for all $\e>0$, such that 
$$|d_\frac{j}{n}(x,v)|\le
D_8(\frac{j}{n})|v|^3   \qquad
\forall (x,v)\in\T^p\times\R^p.   $$
By the formula above, $d_\frac{j}{n}$ satisfies (6.3).

\noindent{\bf Step 2.} We want to show that the rather complicated expression in $(6.7)_b$ and $(6.7)_c$ is not too far from a Gaussian of suitable mean and variance. Note that 
$(6.7)_c$ has the form $e^{-\2\inn{A(v-b)}{v-b}}$; thus, we have to show that $(6.7)_b$ is the "right" normalization coefficient. This is the content of this step.

Since $\g^\frac{1}{n}_\frac{j}{n}(x,\cdot)$ is a probability density for all $x$, we get that $(6.7)_b$ is the reciprocal of the integral of 
$(6.7)_c$ in the variable $v$. We calculate this integral.

First of all, since $f_\frac{j}{n}$ satisfies (5.12), we get the second  inequality below, while the third one comes from standard properties of the Gaussian; the constant $C$ does not depend on anything. 
$$0\le\int_{\R^p\setminus L(n^\frac{-1}{3})}
e^{
-\cinn{v}-f_\frac{j}{n}(x-v)
}     \dr v   \le 
e^{
D_1\left( \frac{j}{n} \right) 
}
\int_{\R^p\setminus L(n^\frac{-1}{3})}
e^{
\cinn{v}
}   \dr v\le
e^{
D_9\left(\frac{j}{n}\right)
}    e^{
-Cn^\frac{1}{6}
}    .  $$
Now, 
$$||
f_\frac{j}{n}(x)+
\2\inn{nQ^{-1}_\frac{j}{n}(x)\b_\frac{j}{n}(x)}{\b_\frac{j}{n}(x)}
||_\infty\le D_{10}\left(\frac{j}{n}\right)   \eqno (6.8)$$
by (6.6) and (5.12). By $(6.7)$ we have that
$$\int_{\R^p\setminus L(n^\frac{-1}{3})}
e^{-\cinn{v}-f_\frac{j}{n}(x-v)} \dr v=
e^{
-f_\frac{j}{n}(x)+\2\inn{nQ_\frac{j}{n}(x)^{-1}\b_\frac{j}{n}(x)}{\b_\frac{j}{n}(x)}
}    \cdot  $$
$$\int_{\R^p\setminus L(n^\frac{-1}{3})}
\exp\left[
-\2\inn{nQ_\frac{j}{n}^{-1}(x)(v-\b_\frac{j}{n}(x))}{v-\b_\frac{j}{n}(x)} -
d_\frac{j}{n}(x,v)
\right]    \dr v   .   $$
The last three formulas imply that
$$0\le \int_{\R^p\setminus L(n^\frac{-1}{3})}
\exp\left[
-\2\inn{nQ_\frac{j}{n}^{-1}(x)(v-\b_\frac{j}{n}(x))}{v-\b_\frac{j}{n}(x)} -
d_\frac{j}{n}(x,v)
\right]    \dr v\le  
e^{
D_{11}\left(\frac{s_2}{n}\right)
}    e^{
-Cn^\frac{1}{6}
}    .  \eqno (6.9)$$
Formula (6.3) implies the two inequalities below. 
$$e^{
-D_4\left(\frac{j}{n}\right)\frac{1}{n}
}\cdot
\int_{L(n^\frac{-1}{3})}
e^{
-\2\inn{nQ^{-1}_\frac{j}{n}(x)(v-\b_\frac{j}{n}(x))}{v-\b_\frac{j}{n}(x)}
}     \dr v\le  $$
$$\int_{L(n^\frac{-1}{3})}
\exp\left[
-\2\inn{nQ^{-1}_\frac{j}{n}(x)(v-\b_\frac{j}{n}(x))}{v-\b_\frac{j}{n}(x)}
-d_\frac{j}{n}(x,v)
\right]     \dr v\le$$
$$e^{
D_4\left(\frac{j}{n}\right)\frac{1}{n}
}\cdot
\int_{L(n^\frac{-1}{3})}
e^{
-\2\inn{nQ^{-1}_\frac{j}{n}(x)(v-\b_\frac{j}{n}(x))}{v-\b_\frac{j}{n}(x)}
}     \dr v  .  $$
The integrals on the left and on the right in the last formula are easy to evaluate: indeed, the Gaussian is centered in 
$\b_\frac{j}{n}(x)$, which satisfies (6.6); thus, almost all its mass (save for an exponentially small rest) lies in $L(n^\frac{-1}{3})$. In formulas,
$$\left(
\frac{(2\pi)^p\det Q_\frac{j}{n}(x)}{n^p}
\right)^\frac{1}{2}
e^{-D_{12}\left(\frac{j}{n}\right)\frac{1}{n}}\le
\int_{L(n^\frac{-1}{3})}
\exp\left[
-\2\inn{nQ_\frac{j}{n}^{-1}(x)(v-\b_\frac{j}{n}(x))}{v-\b_\frac{j}{n}(x)} -
d_\frac{j}{n}(x,v)
\right]    \dr v\le$$
$$\left(
\frac{(2\pi)^p\det Q_\frac{j}{n}(x)}{n^p}
\right)^\frac{1}{2}
e^{D_{12}\left(\frac{j}{n}\right)\frac{1}{n}}   .   $$
By the last formula and (6.9), we get that 
$$\left(
\frac{(2\pi)^p\det Q_\frac{j}{n}(x)}{n^p}
\right)^\frac{1}{2}
e^{-D_{13}\left(\frac{j}{n}\right)\frac{1}{n}}\le
\int_{\R^p}
\exp\left[
-\2\inn{nQ_\frac{j}{n}^{-1}(x)(v-\b_\frac{j}{n}(x))}{v-\b_\frac{j}{n}(x)} -
d_\frac{j}{n}(x,v)
\right]    \dr v\le$$
$$\left(
\frac{(2\pi)^p\det Q_\frac{j}{n}(x)}{n^p}
\right)^\frac{1}{2}
e^{D_{13}\left(\frac{j}{n}\right)\frac{1}{n}}   .   $$
We saw above that $(6.7)_b$ is the inverse of the integral above; thus, 
$$\left(
\frac{n^p}{(2\pi)^p\det Q_\frac{j}{n}(x)}
\right)^\frac{1}{2}
e^{-D_{13}\left(\frac{j}{n}\right)\frac{1}{n}}\le
\left(\frac{n}{2\pi}\right)^\frac{p}{2}
\exp\left\{
b_\frac{j}{n}(x)-f_\frac{j}{n}(x)
+\2\inn{nQ_\frac{j}{n}(x)^{-1}\b_\frac{j}{n}(x)}{\b_\frac{j}{n}(x)}   
\right\}    \le$$
$$\left(
\frac{n^p}{(2\pi)^p\det Q_\frac{j}{n}(x)}
\right)^\frac{1}{2}
e^{D_{13}\left(\frac{j}{n}\right)\frac{1}{n}}   .   \eqno (6.10)$$

\noindent{\bf End of the proof.} We saw at the end of step 1 that 
$d_\frac{j}{n}$ satisfies (6.3). Formula (6.2) follows by (6.7) and (6.10). We prove (6.4) and (6.5). 

We begin to write the normalization coefficient $b_\frac{j}{n}$ in the following complicated way.
$$e^{b_\frac{j}{n}(x)}\left( \frac{n}{2\pi} \right)^\frac{p}{2}=$$
$$\left( \frac{n}{2\pi} \right)^\frac{p}{2}
\exp\left\{
b_\frac{j}{n}(x)-f_\frac{j}{n}(x)+
\2\inn{nQ_\frac{j}{n}(x)^{-1}\b_\frac{j}{n}(x)}{\b_\frac{j}{n}(x)}
\right\}   \cdot$$
$$\exp\left\{
f_\frac{j}{n}(x)-
\2\inn{nQ_\frac{j}{n}(x)^{-1}\b_\frac{j}{n}(x)}{\b_\frac{j}{n}(x)}
\right\}   .   $$
Formulas (6.10) and (6.8) give an estimate on the first and second term respectively in the product above; thus,  
$$e^{-D_{14}\left( \frac{j}{n} \right)}
\left(
\frac{n^p}{(2\pi)^p\det Q_\frac{j}{n}(x)}
\right)   \le
e^{b_\frac{j}{n}(x)}
\left( \frac{n}{2\pi} \right)^\frac{p}{2}\le
e^{D_{14}\left( \frac{j}{n} \right)}
\left(
\frac{n^p}{(2\pi)^p\det Q_\frac{j}{n}(x)}
\right)   .   $$
Together with $(6.7)_a$, this implies that
$$e^{-D_{14}\left( \frac{s_2}{n} \right)}
\left(
\frac{n^p}{(2\pi)^p\det Q_\frac{j}{n}(x)}
\right)
e^{
-\cinn{v}-f_\frac{j}{n}(x-v)
}      \le$$
$$\g^\frac{1}{n}_\frac{j}{n}(x,v)\le
e^{D_{14}\left( \frac{s_2}{n} \right)}
\left(
\frac{n^p}{(2\pi)^p\det Q_\frac{j}{n}(x)}
\right)
e^{
-\cinn{v}-f_\frac{j}{n}(x-v)
}    .   $$
Now (6.4) and (6.5) follow from the last formula, (5.12) and well-known properties of the Gaussian.

\fin

\prop{6.2} Let $T\in[-m,0]$, let $s$ be the maximal integer such that $T\le\frac{-s}{n}$ and let $\mu\in\Mt$. Let 
$\{ \mu^\frac{1}{n}_\frac{j}{n} \}_{j=-s}^0$ be a minimal 
$\{ \g^\frac{1}{n}_\frac{j}{n} \}_{j=-s}^{-1}$-sequence starting at 
$(\frac{-s}{n},\mu)$ and let the interpolating curve 
$\{ \mu^\frac{1}{n}_t \}$ be defined as in section 3. Then, up to subsequences $\{ \mu^\frac{1}{n}_t \}$ converges weakly, uniformly on each interval $[T,-\e]$, to a curve $\mu_t$ which is a weak solution of $(FP)_{-T,-\partial_xu,\mu}$, where $u$ is the limit of lemma 5.3. 

\proof Throughout the proof, we shall deal with the sequence 
$\{ n_k \}$ of lemma 5.3, but we shall drop the subscript $k$ to lighten the notation.

It is standard that $(FP)_{T,-\partial_xu,\mu}$ has a weak solution $\mu_{T,t}$; we have to prove that, if $g\in C(\T^p)$ and 
$\e\in (0,T)$, then
$$\sup_{t\in [T,-\e]}\inn{\mu_t^\frac{1}{n}-\mu_{T,t}}{g}\tends 0
\txt{as}n\tends+\infty   \eqno (6.11)$$
where $\inn{\cdot}{\cdot}$ denotes the duality coupling between 
${\cal M}(\T^p)$, the space of signed measures on $\T^p$, and 
$C(\T^p)$. 

Let $(\g^\frac{1}{n}_\frac{-s}{n},\dots,\g^\frac{1}{n}_\frac{-1}{n})$ be as in the hypotheses; for $-s\le j\le -1$ we define
$$\fun{
S^\ast_{\frac{j}{n},\frac{j+1}{n}}
}{
{\cal M}(\T^p)
}{
{\cal M}(\T^p)
},\qquad
\fun{
S^\ast_{\frac{j}{n},\frac{j+1}{n}}
}{
\mu
}{
\mu\ast\g^\frac{1}{n}_\frac{j}{n}
}  .  $$
If $-s\le l\le j\le 0$, we set 
$$\fun{
S^\ast_{\frac{l}{n},\frac{j}{n}}
}{
{\cal M}(\T^p)
}{
{\cal M}(\T^p)
},\qquad
S^\ast_{\frac{l}{n},\frac{j}{n}}(\mu)=
S^\ast_{\frac{j-1}{n},\frac{j}{n}}\circ\dots\circ
S^\ast_{\frac{l}{n},\frac{l+1}{n}}(\mu)  .  $$
Clearly, with this definition $S^\ast_{\frac{l}{n},\frac{j}{n}}$ has the co-cycle property
$$S^\ast_{\frac{j}{n},\frac{i}{n}}\circ 
S^\ast_{\frac{l}{n},\frac{j}{n}}=S^\ast_{\frac{l}{n},\frac{i}{n}}  \txt{for}
-s\le l\le j\le i\le 0  $$
and  
$$\mu^\frac{1}{n}_\frac{j}{n}=S^\ast_{\frac{-s}{n},\frac{j}{n}}\mu  $$
i. e. $(\mu,S^\ast_{\frac{-s}{n},\frac{-s+1}{n}}(\mu),\dots,
S^\ast_{\frac{-s}{n},0}(\mu))$ is a 
$(\g^\frac{1}{n}_\frac{-s}{n},\dots,\g^\frac{1}{n}_\frac{-1}{n})$-sequence.

Let us also introduce the operator
$$\fun{F^\ast_{T,t}}{\mu}{\mu_{T,t}}    $$
where $\mu_{T,t}$ is the solution, at time $t\ge T$, of the Fokker-Planck equation $(FP)_{T,-\partial_x u,\mu}$. 

By the last two formulas, we can write (6.11) as 
$$\sup_{t\in[T,-\e]}\inn{(S_{\frac{-s}{n},\frac{[nt]}{n}}^\ast)\mu-
F^\ast_{T,t}\mu}{g}
\tends 0\txt{as} n\tends+\infty   .   \eqno (6.12)$$
Since 
$S^\ast_{\frac{j}{n},\frac{j+1}{n}}\mu=
\mu\ast\g^\frac{1}{n}_\frac{j}{n}$, the definition of 
$\mu\ast\g^\frac{1}{n}_\frac{j}{n}$ immediately yields that 
$S^\ast_{\frac{j}{n},\frac{j+1}{n}}$ is the adjoint of the operator
$$\fun{S_{\frac{j}{n},\frac{j+1}{n}}
}{C(\T^p)}{C(\T^p)}  \qquad
\fun{S_{\frac{j}{n},\frac{j+1}{n}}
}{g}{
\int_{\R^p}g(x-v)\g^\frac{1}{n}_\frac{j}{n}(x,v)\dr v
}   .   $$
Note that $S_{\frac{j}{n},\frac{j+1}{n}}$ arrives in $C(\T^p)$ by the results of section 5: indeed, in that section we have proven that 
$\g^\frac{1}{n}_\frac{j}{n}$ is a continuous function. Actually, proposition 5.2 implies that $S_{\frac{j}{n},\frac{j+1}{n}}$ is a bounded operator; we can associate to it a co-cycle  as we did with $S^\ast_{\frac{j}{n},\frac{j+1}{n}}$, with the only difference that $S_{\frac{j}{n},\frac{j+1}{n}}$ is going back in time: we are bringing a final condition $\phi$ to $\phi_\frac{-1}{n}$, $\phi_\frac{-2}{n}$, etc... 

Also $F^\ast_{T,t}$ has an adjoint, which is a co-cycle going back in time; namely, for $-m\le T\le t\le 0$ and $u$ as in lemma 5.3, we can define $F_{t,T}g$ to be $\psi_T$, the solution at time $T$ of 
$$\left\{
\eqalign{
\partial_t\psi&=-(\D\psi-\partial_x u\cdot\partial_x\psi)\cr
\psi_t(x)&=g(x)    
}
\right.   $$
Thus, (6.12) is equivalent to  
$$\sup_{t\in[T,-\e]}
\inn{\mu}{S_{\frac{-s}{n},\frac{[nt]}{n}}g-F_{T,t}g}
\tends 0\txt{as}n\tends+\infty   $$
which in turn is implied by
$$\sup_{t\in[T,-\e]}||
S_{\frac{-s}{n},\frac{[nt]}{n}}g-F_{T,t}g
||_\infty  \tends 0\txt{as}n\tends+\infty .  \eqno (6.13)$$
Let $B_\tau$ be the operator
$$\fun{B_\tau}{\psi}{-(\D\psi-\partial_xu(\tau,\cdot)\cdot\partial_x\psi)} . $$
Theorem 6.5 of section 1 of [7] (which holds in the autonomuos case, but is easy to adapt to our situation) says that (6.13) holds if we have that (keeping track that the time is inverted)
$$||
-n[S_{\frac{j}{n},\frac{j+1}{n}}-Id]g-B_\frac{j}{n}g
||_{C^0(\T^p)}   \tends 0   \eqno (6.14)$$
for every $g\in C^2(\T^p)$ uniformly for 
$\frac{j}{n}\in[\frac{-s}{n},-\e]$. Thus, the theorem reduces to proving this formula. The first equality below is the definition of 
$S_{\frac{j}{n},\frac{j+1}{n}}$. 
$$n[
S_{\frac{j}{n},\frac{j+1}{n}}-Id
]  g=
n\int_{\R^p}[
g(x-v)-g(x)
]  \g^\frac{1}{n}_\frac{j}{n}(x,v)\dr v=$$
$$n\int_{L(n^\frac{-1}{3})}[
g(x-v)-g(x)
]  \g^\frac{1}{n}_\frac{j}{n}(x,v)\dr v   +   \eqno (6.15)_a$$
$$n\int_{\R\setminus {L(n^\frac{-1}{3})}}[
g(x-v)-g(x)
]  \g^\frac{1}{n}_\frac{j}{n}(x,v)\dr v . \eqno (6.15)_b$$
By (6.4) and (6.5) of lemma 6.1 we get that 
$$(6.15)_b\tends 0   .   \eqno (6.16)$$
As for $(6.15)_a$, we want to substitute $\g^\frac{1}{n}_\frac{j}{n}$ with the Gaussian given by lemma 6.1. Indeed, since  $g$ is continuous, there is $\d_n\tends 0$ as 
$n\tends+\infty$ such that the first inequality below holds; the second one follows by (6.2) and (6.3), while the last one follows from the fact that the Gaussian has integral one.
$$\left\vert
\int_{L(n^\frac{-1}{3})}
[g(x-v)-g(x)]\cdot
\left[
N(\b_\frac{j}{n}(x),\frac{1}{n}Q_\frac{j}{n}(x))(v)-
\g^\frac{1}{n}_\frac{j}{n}(x,v)
\right]\dr v
\right\vert   \le  $$
$$\d_n\int_{L(n^\frac{-1}{3})}
\left\vert
N(\b_\frac{j}{n}(x),\frac{1}{n}Q_\frac{j}{n}(x))(v)-
\g^\frac{1}{n}_\frac{j}{n}(x,v)
\right\vert    \dr v\le$$
$$\d_n\int_{L(n^\frac{-1}{3})}
N(\b_\frac{j}{n}(x),\frac{1}{n}Q_\frac{j}{n}(x))(v)
\cdot[e^{\frac{1}{n}D_4\left( \frac{j}{n} \right)}  -1]   \dr v\le
\d_n\frac{1}{n} D_{15}\left( \frac{j}{n} \right)    .   $$
We multiply by $n$ and arrange the terms in a different way.
$$n\int_{L(n^\frac{-1}{3})}[
g(x-v)-g(x)
]  N\left(\b_\frac{j}{n}(x),\frac{1}{n}Q_\frac{j}{n}(x)\right) (v)\dr v
-\d_nD_{15}\left( \frac{j}{n} \right)\le$$
$$n\int_{L(n^\frac{-1}{3})}[
g(x-v)-g(x)
]   \g^\frac{1}{n}_{j}(x,v)\dr v\le$$
$$n\int_{L(n^\frac{-1}{3})}[
g(x-v)-g(x)
]  N\left(\b_\frac{j}{n}(x),\frac{1}{n}Q_\frac{j}{n}(x)\right) (v)\dr v+
\d_nD_{15}\left( \frac{j}{n} \right)  .  $$
By (6.1) and lemma 5.3, if $\frac{j}{n}\tends t$, then 
$n\b_\frac{j}{n}(x)\tends\partial_xu(t,x)$; on the other hand, 
$Q_\frac{j}{n}\tends Id$ by (5.12). Since $g\in C^2(\T^p)$, this implies in a standard way that the two integrals on the left and on the right of the formula above converge to 
$\D u-\partial_xu\cdot\partial_xg$; thus, 
$$(6.15)_a\tends\D g-\partial_xu\cdot\partial_xg  .  $$
Now (6.14) follows from (6.15), (6.16) and the last formula.

\fin 

This immediately calls for a definition.

\vskip 1pc

\noindent{\bf Definition.} Let $\mu_s$ be a limit as in proposition 6.2; in particular, it satisfies $\mu_{T}=\mu$. Then we say that $\mu_s$ is a limit minimal characteristic starting at $(T,\mu)$.

\prop{6.3} Let $T\in[-m,0)$ and let $s$ be the maximal integer such that $T\le\frac{-s}{n}$. Let 
$\{ \mu^\frac{1}{n}_\frac{j}{n} \}_j$ be a minimal 
$\{ \g^\frac{1}{n}_\frac{j}{n} \}_j$-sequence starting at 
$(\frac{-s}{n},\mu)$. Let 
$f^\frac{1}{n}_t$ be defined by (5.4); then, there is 
$f_0\in L^\infty(\T^p)$ such that the following holds. 

Up to subsequences, $f^\frac{1}{n}_t$ converges to a function $u$ which satisfies 
$(HJ)_{0,\bar\mu,f_0}$, where $\bar\mu_t$ is a limit minimal characteristic. The convergence is in 
$C([T,-\e],C^2(\T^p))$ for all $\e\in(0,T)$.

\proof {\bf Step 1.} In this step, we want to reduce to the situation of [14], i. e. to a problem where neither the potential nor the final condition depend on $n$.

As in proposition 6.2, we can interpolate the measures 
$\{ \mu^\frac{1}{n}_\frac{j}{n} \}_j$ by a curve of measures 
$\mu^\frac{1}{n}_t$. Taking subsequences, we can suppose that 
$f^\frac{1}{n}_t\tends u$ (lemma 5.3) and that 
$\mu^\frac{1}{n}_t\tends\mu$ (proposition 6.2). By point $ii$) of the definition of differentiability on densities, we can further refine our subsequence in order to have 
$e^{-f_0^\frac{1}{n}}\weak e^{-f_0}$ in $L^1(\T^p)$. 

For $\mu^\frac{1}{n}_t$ and its limit $\bar\mu_t$, we define as above
$$P^\frac{1}{n}(t,x)=V(t,x)+
W^{\mu_t^\frac{1}{n}}(x),\qquad
\bar P(t,x)=V(t,x)+ 
W^{\bar\mu_t}(x)  .  
$$
Note that the function $\bar P$ does not depend on $n$, since it is defined in terms of $\bar\mu_t$ which does not depend on $n$.

Since the curve of measures $\mu^\frac{1}{n}_t$ converges uniformly to $\bar\mu_t$ on $[T,-\e]$ for all $\e>0$, by the definition of $W^{\2\mu^\frac{1}{n}_t}$ we have that
$$\sup_{\frac{j}{n}\in[T,-\e]}||
P^\frac{1}{n}(\frac{j}{n},\cdot)-\bar P(\frac{j}{n},\cdot)
||_{C^4(\T^p)}       \le\d_n  \eqno (6.17)$$
with $\d_n\tends 0$ as $n\tends+\infty$.
We defined $\{ f_\frac{j}{n}^\frac{1}{n} \}_j$ as the linearized cost  for the problem with final condition 
$f_0^\frac{1}{n}=U^\prime(\mu_0^\frac{1}{n})$ and potential 
$P^\frac{1}{n}$; we let $\{ \bar f_\frac{j}{n}^\frac{1}{n} \}_j$ be the linearized cost for the problem with final condition 
$f_0$ and potential $\bar P$.

Since neither the potential nor the final condition for 
$\bar f^\frac{1}{n}_t$ depend on $n$, we are exactly in the case of [14]; by theorem 29 of [14],  $\bar f^\frac{1}{n}_t$  converges to a solution of $(HJ)_{0,\bar\r,f_0}$ as $n\tends+\infty$; thus, it suffices to show that, for all $\e>0$, 
$$\sup_{\frac{j}{n}\in[T,-\e]}
|| 
f^\frac{1}{n}_\frac{j}{n}-\bar f^\frac{1}{n}_\frac{j}{n}
||_{C^2(\T^p)}    \tends 0   \txt{as} n\tends+\infty   .   
\eqno (6.18)$$

\noindent {\bf Step 2.} Here we show that the Feynman-Kac formula (5.9) implies (6.18).

We define ${\cal P}_\frac{j}{n}$ as in formula (5.7); we define 
$\bar{\cal P}_\frac{j}{n}$ analogously, but for the potential 
$\bar P$. We set
$$c_\frac{j}{n}(x,y)\colon=
N(0,\frac{|j|}{n}Id)(x-y)
E_{0,0}\left[
{\cal P}_\frac{j}{n}(
x-a_{x-y}(\frac{j}{n})-\tilde w(\frac{j}{n}) , \dots,
x-a_{x-y}(-\frac{1}{n})-\tilde w(-\frac{1}{n})
)
\right] $$
and
$$\bar c_\frac{j}{n}(x,y)\colon=
N(0,\frac{|j|}{n}Id)(x-y)
E_{0,0}\left[
\bar{\cal P}_\frac{j}{n}(
x-a_{x-y}(\frac{j}{n})-\tilde w(\frac{j}{n}) , \dots,
x-a_{x-y}(-\frac{1}{n})-\tilde w(-\frac{1}{n})
)
\right]   .    $$
By (5.9) we get
$$e^{-f^\frac{1}{n}_\frac{j}{n}(x)}=
\int_{\R^p}c_\frac{j}{n}(x,y)e^{-f^\frac{1}{n}_0(y)}\dr y
\txt{and}
e^{-\bar f^\frac{1}{n}_\frac{j}{n}(x)}=
\int_{\R^p}\bar c_\frac{j}{n}(x,y)e^{-f_0(y)}\dr y  .  $$
Now, by the triangle inequality,
$$||e^{-f^\frac{1}{n}_\frac{j}{n}}-
e^{-\bar f^\frac{1}{n}_\frac{j}{n}}  ||_{C^2(\T^p)}\le  
\int_{\R^p}||
\bar c_\frac{j}{n}(\cdot,y)-c_\frac{j}{n}(\cdot,y)
||_{C^2(\R^p)}    e^{-f_0(y)}\dr y+$$
$$\left\vert\left\vert\int_{\R^p}
c_\frac{j}{n}(\cdot,y)
[ e^{-f_0(y)}-e^{-f_0^\frac{1}{n}(y)}]\dr y 
\right\vert\right\vert_{C^2(\R^p)}    .  $$
Thus, (6.18) follows if we prove that 
$$\sup_{-m\le\frac{j}{n}\le-\e}\left\vert\left\vert
\int_{\R^p}
c_\frac{j}{n}(\cdot,y)[e^{-f_0(y)}-e^{-f_0^\frac{1}{n}(y)}]\dr y
\right\vert\right\vert_{C^2(\R^p)}  \tends 0
\txt{as}  n\tends+\infty      \eqno (6.19)$$
and (recalling that $f_0$ is bounded by the definition of differentiability on densities)
$$\sup_{-m\le\frac{j}{n}\le-\e}
\int_{\R^p}
||c_\frac{j}{n}(\cdot,y)-\bar c_\frac{j}{n}(\cdot,y)||_{C^2(\R^p)}
\dr y
\tends 0   \txt{as} n\tends+\infty   .   \eqno (6.20)  $$
We begin with (6.19). For $L(a)$ defined as in lemma 6.1, we see that, by the periodicity of $f_0$ and $f_0^\frac{1}{n}$, 
$$\int_{\R^p}c_\frac{j}{n}(x,y)
[e^{-f_0(y)}-e^{-f_0^\frac{1}{n}(y)}]\dr y=
\int_{L(\2)}\left[
\sum_{k\in\Z^p}c_\frac{j}{n}(x,y+k)
\right]   \cdot
[e^{-f_0(y)}-e^{-f_0^\frac{1}{n}(y)}]\dr y   .   $$
Since we saw at the beginning of the proof that 
$e^{-f_0^\frac{1}{n}}\weak e^{-f_0}$ in $L^1(\T^p)$, (6.19) follows if we prove that the set of function of $y$ 
$$\left\{
\sum_{k\in\Z^p}\partial^l_x c_\frac{j}{n}(x,y+k)\st
x\in\T^p,\quad-m\le\frac{j}{n}\le-\e,\quad 0\le l\le 2,\quad x\in\T^p 
\right\}    $$
is relatively compact in $L^\infty(L(\2))$; since $c_\frac{j}{n}$ is the product of a Gaussian with variance greater than $\e$ and a periodic function bounded in $C^4$, this follows easily by Ascoli-Arzel\`a. 

We prove (6.20). We recall that, by the definition of $c_\frac{j}{n}$, 
$\bar c_\frac{j}{n}$ and (5.7), 
$$c_\frac{j}{n}(x,y)-\bar c_\frac{j}{n}(x,y)=
N\left( 0,\frac{|j|}{n}Id \right)(x-y)\cdot$$
$$E_{0,0}\left\{
\exp\left[
\frac{1}{n}\cdot\sum_{r=j}^0
P_\frac{r+1}{n}\left(
x-a_{x,y}\left(\frac{r}{n}\right)+\tilde w\left(\frac{r}{n}\right)
\right) 
\right]  -
\exp\left[\frac{1}{n}\cdot\sum_{r=j}^0
\bar P_\frac{r+1}{n}\left(
x-a_{x,y}\left(\frac{r}{n}\right)+\tilde w\left(\frac{r}{n}\right)
\right)
\right]
\right\}   .   $$
The first term in the product above is a Gaussian, which is bounded in $C^2$ if $\frac{|j|}{n}\ge\e$; as for the second one, it is easy to see that it tends to zero by (6.17). This implies (6.20).

\fin

\noindent{\bf End of the proof of theorem 1.} By proposition 6.3, the limit $u$ of the linearized value functions satisfies 
$(HJ)_{0,\mu,f_0}$, while the limit minimal characteristic satisfies 
$(FP)_{-m,-\partial_xu,\mu}$ by proposition 6.2; thus, the only things we have to prove are (1) and the semigroup property of 
$\Lambda^m$. As for the latter, it follows in a standard way from (1) (see for instance theorem 4 of [4]); thus, we shall skip its proof.

We prove 1). Let $\mu_t$ be a limit minimal characteristic starting at $(T,\mu)$; let us call $\r_t$ its density ($\mu_t$ has a density since the drift $-\partial_xu$ is regular by proposition 5.2 ) and let $u$ be the solution of the associated Hamilton-Jacobi equation; let us define the drift $Y$ as $Y=-\partial_x u$. Proposition 6.2 implies that 
$\mu_t$ is the push-forward of the Wiener measure by the map 
$\fun{}{\o}{\xi(t)(\o)}$, where $\xi$ solves the stochastic differential equation
$$\left\{
\eqalign{
\dr\xi(t)&= Y(t,\xi(s))\dr t+\dr w(t)\qquad t\in[T,0]\cr
\xi(T)&= X       
}     \right.    $$
and $X$ has distribution $\mu$. This implies that
$$\inf\{
E_w\int_T^0L_c^{\2\r}(t,\xi(t),Y(t,\xi(t)))\dr t+U(\r(0)\L^p)
\}  \le U(T,\mu)$$
where $E_w$ denotes expectation with respect to the Wiener measure. We must prove the opposite inequality.

Let $Y$ be a Lipschitz drift, and let 
$\g^\frac{1}{n}_\frac{j}{n}(x,\cdot)$ be the law of 
$\xi(\frac{j+1}{n})-x$, where $\xi$ solves
$$\left\{
\eqalign{
\dr \xi(t)&=Y(t,\xi(t))+\dr w(t)\cr
\xi(\frac{j}{n})&=x   .   
}
\right.    $$
Now we consider a $\{ \g_\frac{j}{n} \}_j$-sequence 
$\{ \mu_\frac{j}{n} \}_j$ starting at $(\frac{s}{n},\mu)$; by lemma 3.1, we have that
$$\sum_{j=s}^{-1}
\int_{\T^p\times\R^p}\left[
\frac{1}{n}L_c^{\2\mu_\frac{j}{n}}(t,x,nv)+
\log\g_\frac{j}{n}(x,v)
\right]   \g_\frac{j}{n}(x,v)\dr\mu_\frac{j}{n}\dr v   +
U(\mu_0)\ge   
\hat U(\frac{j}{n},\mu)  .  $$
It is easy to see that, if we let $n\tends+\infty$, the left hand side converges to 
$$E_w
\int_t^0L_c(t,\xi(t),Y(t,\xi(t)))\dr t+U(\r(0)\L^p)   .  $$
From this, the opposite inequality follows.

\fin

\vskip 2pc
\centerline{\bf Bibliography}

\noindent [1] L. Ambrosio, W. Gangbo, Hamiltonian ODE's in the Wasserstein space of probability measures, Communications on Pure and Applied Math., {\bf 61}, 18-53, 2008.

\noindent [2] L. Ambrosio, N. Gigli, G. Savar\'e, Gradient Flows, Birkhaeuser, Basel, 2005.

\noindent [3] L. Ambrosio, N. Gigli, G. Savar\'e, Heat flow and calculus on metric measure spaces with Ricci curvature bounded below - the compact case. Preprint 2012. 

\noindent [4] U. Bessi, Viscous Aubry-Mather theory and the Vlasov equation, Discrete and Continuous Dynamical Systems, 
{\bf 34}, 379-420, 2014.




\noindent [5] C. Dellacherie, P-A. Meyer, Probabilities and potential, Paris, 1978.

\noindent [6] I. Ekeland, R. Temam, Convex analysis and variational problems, Amsterdam, 1976.

\noindent [7] S. N. Ethier, T. G. Kurtz, Markov processes, Wiley, New York, 1986.

\noindent [8] J. Feng, T. Nguyen, Hamilton-Jacobi equations in space of measures associated with a system of conservation laws,  J. Math. Pures Appl., {\bf 97}, 318-390, 2012. 

\noindent [9] W. Gangbo, A. Tudorascu, Lagrangian dynamics on an infinite-dimensional torus; a weak KAM theorem, Adv. Math., 
{\bf 224}, 260-292, 2010.

\noindent [10] W. Gangbo, A. Tudorascu, Weak KAM theory on the Wasserstein torus with multi-dimensional underlying space, preprint.

\noindent [11] I. M. Gel'fand, A. M. Yaglom, Integration in functional spaces and its applications in Quantum Physics, J. Math. Phys. {\bf 1}, 48-69, 1960.

\noindent [12] W. Gangbo, T. Nguyen, A. Tudorascu, Hamilton-Jacobi equations in the Wasserstein space, Methods Appl. Anal., {\bf 15}, 155-183, 2008.

\noindent [13] D. Gomes, A stochastic analog of Aubry-Mather theory, Nonlinearity, {\bf 15}, 581-603, 2002.

\noindent [14] D. Gomes, E. Valdinoci, Entropy penalization method for Hamilton-Jacobi equations, Advances in Mathematics, {\bf 215}, 94-152, 2007.


\noindent [15] R. Jordan, D. Kinderlehrer, F. Otto, The variational formulation of the Fokker-Planck equation, SIAM J. Math. Anal., {\bf 29}, 1-17, 1998.

\noindent [16] C. Villani, Topics in optimal transpotation, Providence, R. I., 2003.

\end